\documentclass{article}
\pdfoutput=1
\usepackage{graphicx, amssymb, amsmath}
\usepackage{amsthm}
\usepackage{mathrsfs}
\usepackage{color}
\usepackage{caption}
\usepackage{algorithm}
\usepackage{algorithmic}
\usepackage{upgreek}
\usepackage{subfigure}
\usepackage{color,soul}
\usepackage[colorlinks,
            linkcolor=blue,
            anchorcolor=blue,
            citecolor=blue]{hyperref}
\usepackage{subfigure}
\usepackage{setspace}
\sethlcolor{blue}
\usepackage{float}
\usepackage{natbib}
\setcitestyle{numbers,square}

\def\k{\kern .5em}
\def\er{\kern .2em}

\begin{document}

\date{}
\author{}
\newcommand{\be}{\begin{equation}}
\newcommand{\ee}{\end{equation}}
\newcommand{\ba}{\begin{array}}
\newcommand{\ea}{\end{array}}
\newcommand{\beas}{\begin{eqnarray*}}
\newcommand{\eeas}{\end{eqnarray*}}
\newcommand{\bea}{\begin{eqnarray}}
\newcommand{\eea}{\end{eqnarray}}
\newcommand{\ome}{\Omega}

\newtheorem{theorem}{Theorem}[section]
\newtheorem{lemma}{Lemma}[section]
\newtheorem{remark}{Remark}[section]
\newtheorem{proposition}{Proposition}[section]
\newtheorem{definition}{Definition}[section]
\newtheorem{corollary}{Corollary}[section]

\newtheorem{theo}{Theorem}[section]
\newtheorem{lemm}{Lemma}[section]
\newcommand{\blem}{\begin{lemma}}
\newcommand{\elem}{\end{lemma}}
\newcommand{\bthe}{\begin{theorem}}
\newcommand{\ethe}{\end{theorem}}
\newtheorem{prop}{Proposition}[section]
\newcommand{\bprop}{\begin{proposition}}
\newcommand{\eprop}{\end{proposition}}
\newtheorem{defi}{Definition}[section]
\newtheorem{coro}{Corollary}[section]
\newtheorem{algo}{Algorithm}[section]
\newtheorem{rema}{Remark}[section]
\newtheorem{property}{Property}[section]
\newtheorem{assu}{Assumption}[section]
\newtheorem{exam}{Example}[section]

\renewcommand{\theequation}{\arabic{section}.\arabic{equation}}
\renewcommand{\thetheorem}{\arabic{section}.\arabic{theorem}}
\renewcommand{\thelemma}{\arabic{section}.\arabic{lemma}}
\renewcommand{\theproposition}{\arabic{section}.\arabic{proposition}}
\renewcommand{\thedefinition}{\arabic{section}.\arabic{definition}}
\renewcommand{\thecorollary}{\arabic{section}.\arabic{corollary}}
\renewcommand{\thealgorithm}{\arabic{section}.\arabic{algorithm}}
\newcommand{\lan}{\langle}
\newcommand{\curl}{{\bf curl \;}}
\newcommand{\rot}{{\rm curl}}
\newcommand{\grad}{{\bf grad \;}}
\newcommand{\dvg}{{\rm div \,}}
\newcommand{\ran}{\rangle}
\newcommand{\bR}{\mbox{\bf R}}
\newcommand{\bRn}{{\bf R}^3}
\newcommand{\Coinf}{C_0^{\infty}}
\newcommand{\disp}{\displaystyle}
\newcommand{\ra}{\rightarrow}
\newcommand{\Ra}{\Rightarrow}
\newcommand{\ud}{u_{\delta}}
\newcommand{\Ed}{E_{\delta}}
\newcommand{\Hd}{H_{\delta}}
\newcommand\varep{\varepsilon}
\newcommand{\RNum}[1]{\uppercase\expandafter{\romannumeral #1\relax}}
\newcommand{\tabincell}[2]{\begin{tabular}{@{}#1@{}}#2\end{tabular}}
\title{Error Analysis of Virtual Element Methods for the Time-dependent Poisson-Nernst-Planck Equations}
\author{*
	\and *
	\and *
}
\author{Ying Yang $^{1}$
	\and Ya Liu $^2$
	\and Shi Shu $^{3,*}$
     }
\footnotetext[1]
{School of  Mathematics and Computational Science, Guangxi Colleges and Universities Key Laboratory of Data Analysis and Computation, Key Laboratory of Intelligent Computing and Information Processing of Ministry of Education, Guilin University of Electronic Technology, Guilin, 541004, Guangxi, P.R. China. E-mail: yangying@lsec.cc.ac.cn }

\footnotetext[2]
{School of Mathematics and Computational Science, Guilin University of Electronic Technology, Guilin, 541004, Guangxi, P.R. China. E-mail: 979504374@qq.com }

%
%
\footnotetext[3]
{$^{*}$\textbf{Corresponding author.} School of Mathematics and Computational Science, Xiangtan University,
	Xiangtan, 411105, P.R. China. E-mail: shushi@xtu.edu.cn }

\maketitle

\noindent
{\bf Abstract}\quad We discuss and analyze the virtual element method on general polygonal meshes for the time-dependent Poisson-Nernst-Planck equations, which are a nonlinear coupled system widely used in semiconductors and ion channels. The spatial discretization is based on the elliptic projection and the $L^2$ projection operator, and for the temporal discretization, the backward Euler scheme is employed.  After presenting the semi and fully discrete schemes, we derive the a priori error estimates in the $L^2$ and $H^1$ norms. Finally, a numerical experiment  verifies the theoretical convergence results.



\section{Introduction}\label{sec1}
\noindent The virtual element method (VEM) could be seen as a deformation of the classical mimetic finite difference method,  
which was originally proposed in \cite{beirao2013basic}  as a generalization of the
finite element method (FEM). This method is applicable to general polygon/polyhedral grids, even including the multiply-connected or non-convex polygon grids, and hence has low requirements on grid quality.
The VEM, in comparison to the traditional FEM, does not require an explicit expression of the discrete basis functions. In addition to that, it only needs to define the appropriate degrees of freedom to convert the discrete formulation into the matrix form. Thanks to its applicability and simplicity, the VEM has been applied to many equations, for instance, the second-order
elliptic equation \cite{beirao2016virtual}, the parabolic equations \cite{adak2019convergence, vacca2015virtual} and hyperbolic equation \cite{vacca2017virtual}, the Stokes equations \cite{Antonietti2014A, Cangiani2016The, Loureno2015Divergence}, the elasticity problems \cite{da2013virtual,Gain2014On}, the plate bending problem \cite{brezzi2013virtual}, etc.

Here, we consider the  time-dependent Poisson-Nernst-Planck (PNP) equations.  The classic PNP equations are a coupled nonlinear system of partial differential equations, which consist of the electrostatic Poisson equation and the Nernst-Planck equation. The coupled nonlinear system was originally derived by W. Nernst \cite{nernst1889elektromotorische}  and M. Planck \cite{planck1890ueber} and has been
widely applied in semiconductors \cite{markowich1985stationary,Jerome1996Analysis}, biological ion channels \cite{eisenberg1998ionic,singer2009poisson} and electrochemical systems \cite{marcicki2012comparison,richardson2007time,van2010diffuse}.

Because of the high nonlinearity and strong coupling, it is difficult to find the analytic solution for PNP equations. Many numerical methods were developed to find the approximate solutions, for instance, finite volume methods \cite{bessemoulin2014study,chainais2004finite}, finite difference methods \cite{flavell2014conservative,he2016energy} and FEMs \cite{lu2010poisson,xie2020effective}, etc. 
 The FEM has been applied to PNP equations for many years and it is popular because of its flexibility and adaptability in dealing with the irregular interface. In recent years, some work on convergence analysis of FEM has emerged. We presented some error bounds in \cite{yang2013error} for a piecewise finite element approximation to the steady-state PNP equations describing the electrodiffusion of ions in a solvated biomolecular system. In \cite{sun2016error}, the authors discussed a priori error estimates of FEM for the time-dependent PNP  equations, where the optimal error estimates are obtained in $L^\infty(H^1)$ and $L^2(H^1)$ norms and the suboptimal error estimate is obtained in the $L^\infty(L^2)$ norm. 
 An optimal $L^2$ norm error estimate of the FEM to a linearized backward Euler scheme for the time-dependent PNP equations has been obtained in \cite{gao2017linearized}. Soon afterwards, the authors in \cite{shi2019superconvergence} gave the superconvergence analysis of the FEM for the time-dependent PNP equations. Based on the superconvergence result, an optimal $L^2$ norm error
estimate for the classic nonlinear backward Euler fully-discrete schemes was presented on the requirement of higher regularity of the solution. Recently, we presented a decoupling two-grid FEM for the time-dependent PNP equations in \cite{shen2019decoupling}.
This method costs less computational time and remains the same order of accuracy compared with the FEM combined with the Gummel iteration. The optimal $L^2$ error estimate for the classic nonlinear backward Euler scheme was also presented in \cite{shen2019decoupling} with a generic regularity assumption of the solution. In \cite{yang2020superconvergent}, we studied the superconvergent gradient recovery based on the finite element approximation for the strong nonlinear PNP equations and applied it to the realistic ion channel problem, and in \cite{shen2020gradient}, we presented the a posteriori error estimates for a class of steady-state PNP equations .
\par In this paper, the main purpose is to provide the a priori error analysis for the virtual element discretization of the time-dependent PNP equations. First, we design a suitable virtual element discretization scheme for the equations. Compared with the finite element discretization for PNP equations, it can be used on very general polygonal meshes and has lower requirements on grid quality. Hence it could be more suitable for PNP problems with extremely irregular interfaces, for example ion channel problems. Second, the optimal $H^1$ norm error estimates with $k$-th ($k\geq 1$) virtual element are presented for semidiscrete virtual element approximations. At last, the suboptimal $L^2$ norm error estimates with $k$-th ($k\geq 1$) virtual element are given for fully discrete virtual element approximations.

\par From a mathematical point of view, the PNP equations consist of a linear elliptic (Poisson) equation and two nonlinear parabolic (NP) equations. We follow the frame of convergence analysis in \cite{liu22} to present the error estimate for the elliptic equation. Some arguments in \cite{sun2016error} and \cite{vacca2015virtual} are used in the analysis of the nonlinear parabolic equation. Compared with these relevant work, we have some own characteristics in the analysis. For example,
different from the linear energy projection used in \cite{vacca2015virtual}, the energy projection applied in this paper is nonlinear, the error estimate of which is more complex. Compared with \cite{sun2016error}, although the definition of the the energy projection is similar, the analysis in this paper is improved than that presented in \cite{sun2016error}.
The existence or uniqueness of the nonlinear FEM energy projection is not discussed in \cite{sun2016error}, while the detailed proof of the existence and uniqueness of the nonlinear VEM energy projection are presented. The proof is not trivial, in which the error estimate of the $L^2$ projection in maximum norm needs to be discussed, see Lemma \ref{piinfty}, and more detailed discussion of the nonlinear form are also need to be given, see Lemma \ref{Bi}-\ref{lemmsup}. Our scheme was motivated by the work \cite{liu22} in which a VEM is provided for steady-state PNP equations. However, since the NP equation considered in this paper is a parabolic equation, while it is an elliptic one in \cite{liu22}, the frame of most of the analysis is quite different.
In addition, we also improved the error results from the following two aspects. In \cite{liu22}, the $H^1$ norm error estimates depend on the error estimate of  concentrations in $L^2$ norm, which is difficult to present. 
In this paper, we first present the suboptimal $L^2$ error estimates  of the virtual element solutions, see Theorem \ref{theorem3.1} and \ref{theorem3.3}, then the optimal error estimates of the $H^1$ norm are obtained,
see Theorem \ref{theorem_p_H1}. To get the error estimates of the virtual element method, the mathematical analysis in \cite{liu22} requires an assumption for the virtual element solution, e.g.  $\phi_h\in W^{1,\infty}$. In this paper, to avoid use this assumption, we show the boundness of the $L^2$ projection, see Lemma \ref{phiinfty}. Based on the result, the optimal error estimates of the virtual element approximation in the $H^1$ norm are obtained without the regularity assumption of
the discrete solution.

\par The rest of this paper is organised as follows. In Section 2, we introduce the time-dependent PNP equations and present the corresponding weak form.  In Section 3, the virtual element space and the corresponding semidiscrete virtual element approximation are discussed, and the convergence analysis in the $L^2$ and $H^1$ norms are derived. Section 4 introduces  the fully discrete scheme and  presents the convergence analysis in the $L^2$ norms. To confirm the efficiency of the proposed methods and the accuracy of the theoretical analysis, a numerical example is given in Section 5. Finally, some conclusions are made in Section 6.

\setcounter{lemma}{0}
\setcounter{theorem}{0}
\setcounter{corollary}{0}
\setcounter{equation}{0}
\setcounter{remark}{0}
\section{ The time-dependent Poisson-Nernst-Planck equations}

 Let $\Omega\subset \mathbb{R}^{2}$ be a bounded polygonal domain and $\partial\Omega$ be the Lipschitz continuous boundary of $\Omega$.
 We consider the following time-dependent PNP equations (cf. \cite{gao2017linearized,sun2016error}):
\begin{equation}\label{pnp equation}
\left\{\begin{array}{lr}
{p_t^{i}}-\nabla\cdot(\nabla{p^{i}}+q^{i}p^{i}\nabla\phi)=F^i, & \text { in } \Omega,~\text{for}~t\in(0, T], i=1,2, \vspace{1mm}\\
-\triangle\phi-\sum\limits^{2}_{i=1}q^{i}p^{i}=f, & \text { in } \Omega,~\text{for}~t\in(0, T],\vspace{1mm}\\
p^i(\cdot,0)=p^i_{0},~~\phi(\cdot,0)=\phi_0,~~f(\cdot,0)=f_0, ~~F^i(\cdot,0)=F_0^i,& \text { in } \Omega,
\end{array}\right.
\end{equation}
with the homogeneous Dirichlet boundary conditions
\begin{equation}\label{boundary}
\left\{
\begin{array}{c}
\;\phi=0,\;\text{on}\;\partial\Omega,~\text{for}~t\in(0, T],\vspace{0.5mm}\\
\;p^i=0,\;\text{on}\;\partial\Omega,~\text{for}~t\in(0, T],
\end{array}
\right.
\end{equation}
where  $p^i,~i=1,2$ denotes the concentration of the $i$-th ionic species, $p_t^{i}=\frac{\partial p^i}{\partial t}$, $\phi$ represents the electrostatic potential, the constant $q^i$ corresponds to the charge of the species $i$, $f$ and $F^{i}$ are the reaction source terms, $p^{i}_0,~\phi_{0},~f_{0}~\text{and}~F^i_0\in L^{2}(\Omega)$  denote initial data.\\

\par For any $u,v,\psi,u^1,u^2\in H _0^{1}(\Omega)$, define
\begin{equation}\label{2.4}
a(u,v)=(\nabla u,\nabla v),~ b_i(u,\psi,v)=(q^iu\nabla\psi,\nabla v),~ \tilde{b}(u^1,u^2,v)=(-\sum\limits_{i=1}^{2}q^iu^i,v).
\end{equation}
The weak formulation of (\ref{pnp equation})-(\ref{boundary}) is given as: find $p^i\in L^{2}\left(0, T; H_{0}^{1}(\Omega)\right), i=1,2,$ 
and $\phi\in L^{2}\left(0, T; H_{0}^{1}(\Omega)\right)$ such that
\begin{small}
\begin{equation}\label{weak formulation}
%
\left\{\begin{array}{lr}
(p_t^i(t),v)+a(p^i(t),v)+b_i(p^i(t),\phi(t),v)=(F^i(t),v),&~\forall v\in H _0^{1}(\Omega),~i=1,2,\vspace{2mm}\\
a(\phi(t),w)+\tilde{b}(p^1(t),p^2(t),w)=(f(t),w),&~\forall w\in H _0^{1}(\Omega),\vspace{2mm}\\
p^i(\cdot,0)=p^i_{0},~~\phi(\cdot,0)=\phi_0,~~f(\cdot,0)=f_0,~~F^i(\cdot,0)=F_0^i.\end{array}\right.
\end{equation}
\end{small}
The existence and uniqueness of  the solution for (\ref{weak formulation}) have been shown in \cite{gajewski1986basic} for $F_{i}=R\left(p^{1}, p^{2}\right)=r\left(p^{1}, p^{2}\right)\left(1-p^{1} p^{2}\right)$ with a Lipschitzian function $r: R_{+}^{2} \rightarrow R_{+}$. In this case, system (\ref{pnp equation})-(\ref{boundary}) describes the transport of mobile carriers in a semiconductor device,  $R\left(p^{1}, p^{2}\right)$ is the net recombination rate  and $p^1$ and $p^2$ represent the densities of mobile holes and electrons, respectively (see \cite{gajewski1986basic}).\\
\par In the rest of the paper, we follow the standard notations for Sobolev spaces $W^{s,p}(\Omega)$ with $\|\cdot\|_{s,p,\Omega}$ and $|\cdot|_{s,p,\Omega}$ denote the norm and the seminorm (cf.~\cite{adams1975sobolev, brenner1994mathematical}), respectively.
For $p=2$, the notations 
$H^{s}(\Omega)=W^{s,2}(\Omega), H_0^1(\Omega)=\{v\in H^{1}(\Omega):v|_{\partial\Omega}=0\}$. For simplicity, denote by $\|\cdot\|_{s}=\|\cdot\|_{W^{s,2}(\Omega)}$, $\|\cdot\|_{0}=\|\cdot\|_{L^2(\Omega)}$ and $\|\cdot\|_{0,\infty}=\|\cdot\|_{L^{\infty}(\Omega)}$.
 We adopt $(\cdot,\cdot)$ to denote the standard $L^{2}$-inner product.
\vspace{3mm}
\setcounter{lemma}{0}
\setcounter{theorem}{0}
\setcounter{corollary}{0}
\setcounter{equation}{0}
\section{The VEM semidiscrete scheme}
In this section, we present local and global virtual element space and outline the semi-discrete virtual element formulation of  (\ref{weak formulation}). Let $\left\{\mathcal{T}_{h}\right\}$ be a family of decompositions of $\Omega$ into elements $E$ with $h_E =diam(E)$
and $h=max\{h_E: E \in \mathcal{T}_h\}$. On each element $E$,  we suppose $a^{E}(\cdot, \cdot),~ b_{i}^{E}(\cdot, \cdot,\cdot) ~\text{and}~ \tilde{b}^{E}(\cdot, \cdot, \cdot)$ are the restrictions of the corresponding forms defined in (\ref{2.4}) on $E$.
Following \cite{ahmad2013equivalent, beirao2013basic}, we make the following assumption for the mesh $\mathcal{T}_h$.
\begin{assu}\label{assu}
Every element $E$ is star-shaped with respect to a ball of radius  greater than $\gamma h_E$ and the distance between any two vertices of $E$ is greater than or equal to $c h_E$, where $\gamma$ and $c$ are uniform positive constants.
\end{assu}

\par Next, let us begin with defining the virtual element space. For any integer $k \geq 0$, denote by $\mathbb{P}_k(D)$ the space of polynomial functions with the total degree up to $k$ living on $D$. Following \cite{beirao2016virtual}, for every element $E \in \mathcal{T}_h$, we  introduce a local space
\begin{eqnarray}\label{3.6}
\tilde{Q}_h^k(E)=\{v\in H^1(E):~ v|_{\partial E}\in C^0(\partial E), v|_e\in \mathbb{P}_k(e), \forall e\in \partial E, \triangle v\in \mathbb{P}_k(E)\}.\notag
\end{eqnarray}
The projection operator $\Pi_{k}^\nabla:\tilde{Q}_h^k(E)\rightarrow \mathbb{P}_k(E)$  is given by
\begin{eqnarray}\label{3.7}
(\nabla(\Pi_{k}^\nabla v_h-v_h),\nabla q)_E=0, ~~\int_{\partial E}(\Pi_{k}^\nabla v_h-v_h)ds=0,\quad\forall q\in \mathbb{P}_k(E).
\end{eqnarray}
\par Then we define the following local  virtual space
\begin{equation}\label{3.8}
Q_h^k(E)=\{v_h\in \tilde{Q}_h^k(E):~(v_h-\Pi_{k}^\nabla v_h,q)_E=0,\quad\forall q\in\big(\mathbb{P}_k/\mathbb{P}_{k-2}(E)\big)\},\notag
\end{equation}
\noindent where $\mathbb{P}_k/\mathbb{P}_{k-2}(E)$ denotes the polynomials in $\mathbb{P}_k(E)$ which are $L^2(E)$ orthogonal to $\mathbb{P}_{k-2}(E)$.
For $v_h \in {Q}_h^k(E)$, we define the following local degrees of freedom (cf.~\cite{ahmad2013equivalent}):\\
\vspace{-2mm}
\par (D1) The values of $v_h$ at the vertices of E.\\\vspace{-5mm}
\par(D2) For $k > 1$, the edge moments $\int_e v_hp_{k-2}ds,~p_{k-2}\in \mathbb{P}_{k-2}(e)$,~on each edge $e$ of $E$.\\\vspace{-5mm}
\par(D3) For $k > 1$, the internal moments $\int_Ev_hp_{k-2}dx, ~p_{k-2}\in\mathbb{P}_{k-2}(E)$.\\
\vspace{-7mm}\\
\par Finally, we can define the global virtual element space as follows:
\begin{equation}\label{3.9}
Q_h^k=\{v\in H_0^1(\Omega):~v|_E\in Q_h^k(E),\quad\forall E\in \mathcal{T}_h\}.\notag
\end{equation}
Let now $\Pi_{k}^0: Q_h^k(E)\rightarrow \mathbb{P}_k(E)$  be the $L^2$ projection operator defined by
\begin{equation}\label{L^2 projection}
(v_h-\Pi_{k}^0 v_h,q)_E=0,\quad\forall q\in \mathbb{P}_k(E).
\end{equation}
 Then the following approximation properties can be obtained (cf.~\cite{beirao2016virtual})
\begin{align}\label{3.11}&\|\Pi_{k}^0 v-v\|_{m,E}\leq Ch^{s-m}_E|v|_{s,E}, ~~m,s\in \mathbb{N},~~m\leq s\leq k+1,~~\forall v\in H^{s}(E),\\
\label{3.12}&\|\Pi_{k}^0 v\|_{m,E}\leq C\|v\|_{m,E},~~m\in \mathbb{N},~~m\leq k+1,~~\forall v\in H^m(E).
\end{align}
It is shown in \cite{chen2018} that under Assumption \ref{assu} for the mesh $\mathcal{T}_h$, for any element $E\in \mathcal{T}_h$, there is a "virtual triangulation" $\mathcal{T}_E$ of $E$ such that $\mathcal{T}_E$ is uniformly shape regular and quasi-uniform. Then,
from the inverse estimate and (\ref{3.12}), there holds (see \cite{liu22})
\begin{equation}\label{3.13}
\|\Pi_{k}^0 v\|_{0,\infty,E}\leq C\|v\|_{0,\infty,E}, ~~\forall v\in L^{\infty}(E),
\end{equation}
where we have used the element area $|E| \leq C h_{E}^{2}$.
\par Assume $v_I \in Q_h^k$ is the interpolant of $v$, which shares the value of the degrees of freedom with $v$. Then
$v_I$ satisfies (cf. \cite{beirao2017stability,brenner2017some})\\
\begin{equation}\label{3.14}
\|v-v_I\|_{1,E}\leq Ch_E\|v\|_{k+1,E},\quad\forall v\in H^{k+1}(E).
\end{equation}
And for any $v\in H^k(E)$ there exists a $v_{\pi}\in \mathbb{P}_k(E)$ such that (see \cite{beirao2013basic})
\begin{align}\label{vpi} \|v-v_{\pi}\|_{0,E}+h_E\|\nabla v-\nabla v_{\pi}\|_{0,E}\leq C h_E^k\|\nabla v\|_{k,E}.
\end{align}

\par For any $u_h, v_h, \psi_h, u_h^1, u_h^2$ on each element $E$, we construct the local forms
\begin{spacing}{2.0}\vspace{-13mm}
\begin{align}
&m_{h}^{E}\left(u_{h}, v_{h}\right)=\int_E[\Pi_{k}^{0} u_{h}]\cdot[\Pi_{ k}^{0} v_{h}]dx+S_m^E\big((I-\Pi_{ k}^{0}) u_{h},(I-\Pi_{ k}^{0}) v_{h}\big)_{E}~,\notag\\
&a_h^E(u_h,v_h)=\int_E[\Pi_{k-1}^0\nabla u_h]\cdot[\Pi_{k-1}^0\nabla v_h]dx+S_a^E\left((I-\Pi_k^\nabla)u_h, (I-\Pi_k^\nabla)v_h\right),\notag\\
&b_{i,h}^E(u_h,\psi_h,v_h)=\int_Eq^i[\Pi_{k-1}^0 u_h][\Pi_{k}^0 \nabla \psi_h]\cdot[\Pi_{k-1}^0\nabla v_h]dx,~~i=1,2,\notag\\
&\tilde{b}_{h}^E(u_h ^1,u_h
^2,v_h)=-\int_E\Big[\Pi_{k-1}^0\Big(\sum\limits_{i=1}^{2}q^iu_h ^i\Big)\Big][\Pi_{k-1}^0 v_h]dx,\notag\\
&(F_{h}^i,v)_E=\int_EF^i~\Pi_{k}^0 v_hdx,~~i=1,2,~~(f_h,v)_E=\int_Ef~\Pi_{k}^0 v_hdx,\notag
\end{align}
\end{spacing}\vspace{-5mm}
\noindent where the symmetric bilinear forms $S_m^E(\cdot,\cdot)$ and $S_a^E(\cdot,\cdot):Q_h^k(E)\times Q_h^k(E)\rightarrow \mathbb{R}$ satisfy
\begin{align*}
&\alpha_1(v_h,v_h)_E\leq S_m^E(v_h, v_h)\leq \alpha_2(v_h,v_h)_E,\quad\forall v_h\in Q_h^k(E)~~\text{with}~~\Pi_k^0 v_h
=0,
\end{align*}
and
\begin{align}\label{sae}
&\beta_1a^E(v_h,v_h)\leq S_a^E(v_h, v_h)\leq \beta_2a^E(v_h,v_h),\quad\forall v_h\in Q_h^k(E)~~\text{with}~~\Pi_k^\nabla v_h=0,
\end{align}
respectively, for four positive constants $\alpha_1$, $\alpha_2$, $\beta_1$ and $\beta_2$.
Moreover, define the discrete forms
\begin{spacing}{2.0}\vspace{-13mm}
\begin{align}
\label{ahbh}&a_{h}\left(u_{h}, v_{h}\right):=\sum\limits_{E}a_{h}^{E}\left(u_{h}, v_{h}\right), \quad b_{i, h}\left(u_{h}, \psi_{h}, v_{h}\right):=\sum\limits_{E} b_{i, h}^{E}\left(u_{h}, \psi_{h}, v_{h}\right), \quad i=1,2, \vspace{1mm}\\
&m_{h}\left(u_{h}, v_{h}\right):=\sum\limits_{E}m_{h}^{E}\left(u_{h},v_{h}\right),\quad\tilde{b}_{h}\left(u_{h}^{1}, u_{h}^{2}, v\right):=\sum\limits_{E} \tilde{b}_{h}^{E}\left(u_{h}^{1}, u_{h}^{2}, v_{h}\right), \vspace{1mm}\notag\\
&\left(f_{h}, v_{h}\right):=\sum\limits_{E}\left(f_{h}, v_{h}\right)_{E},\quad\left(F_{h}^i, v_{h}\right):=\sum\limits_{E}\left(F_{h}^i, v_{h}\right)_{E},\quad i=1,2.\notag
\end{align}
\end{spacing}\vspace{-8mm}
The  semidiscrete virtual element formulation corresponding to (\ref{weak formulation}) reads as: find $p_{h}^{i} \in L^{2}\left(0, T, Q_{h}^{k}\right) \\\text { with } p_{h, t}^{i} \in L^{2}\left(0, T, Q_{h}^{k}\right),~i=1,2$ and $\phi_{h} \in L^{2}\left(0, T, Q_{h}^{k}\right)$ such~that
\begin{equation}\label{semi form}
\left\{\begin{array}{lr}
m_h(p_{h,t}^i,v_h)+a_{h}\left(p_{h}^{i}, v_{h}\right)+b_{i, h}\left(p_{h}^{i}, \phi_{h}, v_{h}\right)=\left(F_{h}^i, v_{h}\right), &\quad \forall v_{h} \in Q_{h}^{k}~\text{for}~\text{a.e.}~t~\text{in}~(0, T), \vspace{2mm}\\
{a}_{h}\left(\phi_{h}, w_{h}\right)+\tilde{b}_{h}\left(p_{h}^{1}, p_{h}^{2}, w_{h}\right)=\left(f_{h}, w_{h}\right), &\quad \forall w_{h} \in Q_{h}^{k},~\text{for}~\text{a.e.}~t~\text{in}~(0, T), \vspace{2mm}\\
p_h^i(0)=p^i_{h,0},~\phi_h(0)=\phi_{h,0}.
\end{array}\right.
\end{equation}
In order to derive the error  estimates in the $H^1$ norm for the semidiscrete case, we require the following lemmas.
\begin{lemma}(\cite{beirao2016virtual}) (Stability) There exist two positive constants $C_0$ and $C_1$ independent of $h$ and $E$, such that
 \begin{align}\label{sta} C_0 a^E(v_h,v_h)\leq a_h^E(v_h,v_h)\leq C_1 a^E(v_h,v_h),~~\forall v_h\in Q_h^k(E).
 \end{align}
\end{lemma}
From (\ref{sta}), it is easy to obtain
\begin{equation}\label{continuity a}
a_{h}^E\left(u_{h}, v_{h}\right) \leq C\left\|\nabla u_{h}\right\|_{0,E}\left\|\nabla v_{h}\right\|_{0,E},~~\forall u_h,~~v_h\in Q_{h}^{E}.
\end{equation}
and
\begin{equation}\label{continuity u}
a_{h}^E\left(u_{h}, u_{h}\right) \geq C\left\|\nabla u_{h}\right\|_{0,E}^2.
\end{equation}

\begin{lemma} \label{cons a}(K-Consistency) For any $q \in \mathbb{P}_{k}(E)$ and $v_{h} \in Q_{h}^{k}(E)$, such that
	\begin{equation}
	\begin{aligned}\label{ahcon}
	{a}_{h}^{E}\left(q, v_{h}\right)-a^{E}\left(q, v_{h}\right)=0.
	\end{aligned}
	\end{equation}
\end{lemma}
\noindent\textbf{Proof.}~It is easy to get the result from the definitions of $\Pi_{k}^{\nabla}$ and $\Pi_{k-1}^{0}$.
$\hfill\Box$\\
\par From \cite{vacca2015virtual}, the bilinear ${m}_{h}^{E}(\cdot, \cdot)$ satisfies the following consistency property and stability property.
\vspace{2mm}
\begin{lemma}(\cite{vacca2015virtual}) (K-Consistency) For all $\chi \in \mathbb{P}_{k}(E)$ and $v_{h} \in Q_{h}^{k}(E)$, there holds
\begin{equation}\label{cons}
\begin{aligned}
{m}_{h}^{E}\left(\chi, v_{h}\right)=\left(\chi, v_{h}\right)_E.
\end{aligned}
\end{equation}
(Stability) There exits two positive constants $C_*$ and $C^*$ independent of h and E, such that
\begin{equation}\label{stab}
C_* \left(v_{h}, v_{h}\right)_E \leq {m}_{h}^{E}\left(v_{h}, v_{h}\right) \leq C^*\left(v_{h}, v_{h}\right)_E,\quad \forall v_{h} \in Q_{h}^{k}(E).
\end{equation}
\end{lemma}
From the stability conditions, it is easy to get the continuity of $m_h$:
\begin{equation}\label{continuity}
m_h\left(u_{h}, v_{h}\right) \leq C\left\|u_{h}\right\|_{0}\|v_{h}\|_{0},\quad \forall u_{h}, v_{h} \in Q_{h}^{k}(E).
\end{equation}
\par In order to estimate the difference between the continuous form and the discrete bilinear form, we need the following lemma.
\vspace{2mm}
\begin{lemma}(\cite{beirao2016virtual})\label{lemma3.3} For any $u,v\in H^1(E)$, if $\kappa \in L^{\infty}(E)$ and $\lambda\in\left[L^{\infty}(E)\right]^{2}$, then we have the estimate
\begin{equation}\label{3.23}
\begin{aligned}
|(\lambda u, &\nabla v)_{E}-(\lambda \Pi_{k-1}^{0} u, \Pi_{k-1}^{0} \nabla v)_{E}| \leq\|\lambda u-\Pi_{k-1}^{0}(\lambda u)\|_{0, E}\|\nabla v-\Pi_{k-1}^{0} \nabla v\|_{0, E}\\
&+\|\lambda \cdot \nabla v-\Pi_{k-1}^{0}(\lambda \cdot \nabla v)\|_{0, E}\|u-\Pi_{k-1}^{0} u\|_{0, E} +C_{\lambda}\|u-\Pi_{k-1}^{0} u\|_{0, E}\|\nabla v-\Pi_{k-1}^{0} \nabla v\|_{0, E}
\end{aligned}
\end{equation}
and
\begin{equation}\label{3.24}
\begin{aligned}
|(\kappa u, v)_{E}-(&\kappa \Pi_{k-1}^{0} u, \Pi_{k-1}^{0} v)_{E}| \leq\|\kappa u-\Pi_{k-1}^{0}(\kappa u)\|_{0, E}\|v-\Pi_{k-1}^{0} v\|_{0, E} \\
&+\|\kappa v-\Pi_{k-1}^{0}(\kappa v)\|_{0, E}\|u-\Pi_{k-1}^{0} u\|_{0, E}+C_{\kappa}\|u-\Pi_{k-1}^{0} u\|_{0, E}\|v-\Pi_{k-1}^{0} v\|_{0, E}.
\end{aligned}
\end{equation}
\end{lemma}

\vspace{2mm}
The error between the continuous form $b_{i}^E(\cdot,\cdot,\cdot)$ and the discrete form $b_{i,h}^E(\cdot,\cdot,\cdot)$ is estimated as follows.
\begin{lemma} \label{bi-bhi}Suppose $w|_{E} \in W^{k+1, \infty}(E)$ and $w_{h}\in Q_{h}^{k}(E)$ satisfying $\|\Pi_k^0\nabla w_h\|_{0,\infty,E}\leq C$. For any  $u \in L^{\infty}(E) \cap H^{k}(E)$ and $u_{h} \in Q_{h}^{k}(E)$, there holds 
\begin{align} b_i^E(u,w,v_h)-b_{i,h}^E(u_h,w_h,v_h)
\leq &C\Big(h_E^k(\|w\|_{k+1,E}+\|u\|_{k,E})+\|\nabla w_h-\nabla w\|_{0,E}\notag\\
&+\|\Pi_k^0\nabla w_h\|_{0,\infty,E}\|u_h-u\|_{0,E}\Big)
\|\nabla v_h\|_{0,E},~~~\forall v_{h} \in Q_{h}^{k}(E).\notag
\end{align}
\end{lemma}
\textbf{Proof.}~ First, for any $v_{h} \in Q_{h}^{k}(E)$, we get
\begin{align}\label{b-}
&b_i^E(u,w,v_h)-b_{i,h}^E(u_h,w_h,v_h)\notag\\
&=(q^iu\nabla w,\nabla v_h)_E-\big(q^i[\Pi_{k-1}^0u_h][\Pi_{k}^0\nabla w_h],[\Pi_{k-1}^0\nabla
v_h]\big)_E\notag\\
&=\Big((q^iu\nabla w,\nabla v_h)_E-(q^i\Pi_{k-1}^0u\Pi_{k}^0\nabla w,\Pi_{k-1}^0\nabla v_h)_E\Big)\notag\\
&\;\;\;\;-\big(q^i\Pi_{k-1}^0(u_h-u)\Pi_{k}^0\nabla w_h,\Pi_{k-1}^0\nabla v_h\big)_E
-\big(q^i\Pi_{k-1}^0u\Pi_{k}^0(\nabla w_h-\nabla w),\Pi_{k-1}^0\nabla v_h\big)_E\notag\\
&=\RNum{1}_1+\RNum{1}_2+\RNum{1}_3.
\end{align}
\noindent For the term $\RNum{1}_1$, it holds
\begin{align}\label{I1}
\RNum{1}_1&=(q^iu\nabla w,\nabla v_h)_E-(q^i\Pi_{k-1}^0u\Pi_{k}^0\nabla w,\Pi_{k-1}^0\nabla v_h)_E\notag\\
&=\Big{(}q^iu(\nabla w-\Pi_{k}^0\nabla w),\nabla v_h\Big{)}_E+(q^iu\Pi_{k}^0\nabla w,\nabla v_h)_E-(q^i\Pi_{k-1}^0u\Pi_{k}^0\nabla w,\Pi_{k-1}^0\nabla v_h)_E\notag\\
&\leq Ch_E^k\|w\|_{k+1,E}\|\nabla v_h\|_{0,E}+\big((q^iu\Pi_{k}^0\nabla w,\nabla v_h)_E-(q^i\Pi_{k-1}^0u\Pi_{k}^0\nabla w,\Pi_{k-1}^0\nabla v_h)_E\big).
\end{align}
For the second term, setting $\beta=q^i\Pi_{k}^0\nabla w$, then from (\ref{3.11})-(\ref{3.13}) and (\ref{3.23}), we obtain
\begin{align}\label{I11}
&(q^iu\Pi_{k}^0\nabla w,\nabla v_h)_E-(q^i\Pi_{k-1}^0u\Pi_{k}^0\nabla w,\Pi_{k-1}^0\nabla v_h)_E\notag\\
&=(\beta u,\nabla v_h)_E-(\beta\Pi_{k-1}^0u,\Pi_{k-1}^0\nabla v_h)_E\notag\\
&\leq \|\beta u-\Pi_{k-1}^0(\beta u)\|_{0,E}\|\nabla v_h-\Pi_{k-1}^0\nabla v_h\|_{0,E}\notag\\
&\;\;+\|\beta\cdot\nabla v_h-\Pi_{k-1}^0(\beta\cdot\nabla v_h)\|_{0,E}\|u-\Pi_{k-1}^0u\|_{0,E}\notag\\
&\;\;+C_{\beta}\|u-\Pi_{k-1}^0u\|_{0,E}\|\nabla v_h-\Pi_{k-1}^0\nabla v_h\|_{0,E}\notag\\
&\leq C\big(\|\beta u-\Pi_{k-1}^0(\beta u)\|_{0,E}+h_E^k\|u\|_{k,E}\big)\|\nabla v_h\|_{0,E}.
\end{align}
where we have used $w|_{E} \in W^{1, \infty}(E)$. Denote by $\hat{\beta}=q^i\nabla w$, then we have
\begin{align} \label{hatu-u}\|\beta u-\Pi_{k-1}^0(\beta u)\|_{0,E}&\leq \|(\beta-\hat{\beta}) u\|_{0,E}+\|\hat{\beta} u-\Pi_{k-1}^0(\hat{\beta}u)\|_{0,E}
+\|\Pi_{k-1}^0(\hat{\beta}u)-\Pi_{k-1}^0({\beta}u)\|_{0,E}\notag\\
&\leq  \|(\beta-\hat{\beta}) u\|_{0,E}+Ch_E^k\|\nabla wu\|_{k,E}\notag\\
&\leq Ch_E^k(\|w\|_{k+1,E}+\|\nabla wu\|_{k,E}).
\end{align}
Inserting (\ref{I11}) and (\ref{hatu-u}) into (\ref{I1}),  it follows that
\begin{eqnarray}\label{RNum{1}}
\begin{split}
\RNum{1}_1&\leq Ch_E^k(\|w\|_{k+1,E}+\|u\|_{k,E})\|\nabla v_h\|_{0,E},
\end{split}
\end{eqnarray}
where we have used $w|_{E} \in W^{k+1, \infty}(E)$.
For the term $\RNum{1}_2$, it yields
\begin{align}\label{RNum{2}}
\RNum{1}_2&=\big(q^i\Pi_{k-1}^0(u_h-u)\Pi_{k}^0\nabla w_h,\Pi_{k-1}^0\nabla v_h\big)_E\notag\\
&\leq |q^i|\|\Pi_{k}^0\nabla w_h\|_{0,\infty,E}\|\Pi_{k-1}^0(u_h-u)\|_{0,E}\|\Pi_{k-1}^0\nabla v_h\|_{0,E}\notag\\
&\leq C\|\Pi_{k}^0\nabla w_h\|_{0,\infty,E}\|u_h-u\|_{0,E}\|\nabla v_h\|_{0,E}.
\end{align}
For the term $\RNum{1}_3$, we deduce that
\begin{align}\label{RNum{3}}
\RNum{1}_3&=\big(q^i\Pi_{k-1}^0u\Pi_{k}^0(\nabla w_h-\nabla w),\Pi_{k-1}^0\nabla v_h\big)_E\notag\\
&\leq |q^i|\|\Pi_{k-1}^0u\|_{0,\infty,E}\|\Pi_{k}^0(\nabla w_h-\nabla w)\|_{0,E}\|\Pi_{k-1}^0\nabla v_h\|_{0,E}\notag\\
&\leq C\|\nabla w_h-\nabla w\|_{0,E}\|\nabla v_h\|_{0,E}.
\end{align}
Then inserting (\ref{RNum{1}})-(\ref{RNum{3}}) into (\ref{b-}), it yields
\begin{align}
&b_i^E(u,w,v_h)-b_{i,h}^E(u_h,w_h,v_h)\notag\\
&\leq C\Big(h_E^k(\|w\|_{k+1,E}+\|u\|_{k,E})+\|\nabla w_h-\nabla w\|_{0,E}+\|\Pi_{k}^0\nabla w_h\|_{0,\infty,E}\|u_h-u\|_{0,E}\Big)
\|\nabla v_h\|_{0,E}.\notag
\end{align}
This completes the proof.$\hfill\Box$
\vspace{2mm}
\begin{lemma}\label{lemma4.2} Suppose $v_h~\text{and}~u_h^{i}\in Q_{h}^{k}(E),~i=1,2$. There holds
\begin{equation}
\Big|\tilde{b}_{h}(u_h^{1}, u_h^{2},v_h)-\tilde{b}(u^{1},u^{2},v_h)\Big|
\leq C(h^{k+1}\sum\limits^{2}_{i=1}\|u^{i}\|_{k}|v_h|_1+\sum\limits^{2}_{i=1}\|u_h^{i}-u^{i}\|_{0}\|v_h\|_0)
.\notag
\end{equation}
\end{lemma}
\textbf{Proof.}~ From (\ref{3.11}), (\ref{3.12}) and (\ref{3.24}), we get
\begin{align}\label{tilde{b}_{h}}
&\Big|\tilde{b}_{h}(u_h^{1},u_h^{2},v_h)-\tilde{b}(u^{1},u^{2},v_h)\Big|\notag\\
&=\Big|\sum_E\Big\{(\sum\limits_{i=1}^{2}q^iu^{i},v_h)_E-\big(\Pi_{k-1}^0(\sum\limits_{i=1}^{2}q^iu_h^{i}),\Pi_{k-1}^0v_h\big)_E\Big\}\Big|\notag\\
&\leq\sum_E\Big\{\Big|(\sum\limits^{2}_{i=1}q^iu^{i},v_h)_E-\big(\Pi_{k-1}^0(\sum\limits_{i=1}^{2}q^iu^{i}),\Pi_{k-1}^0v_h\big)_E\Big|
+\Big|\Big(\Pi_{k-1}^0(\sum\limits_{i=1}^{2}q^iu^{i}-\sum\limits_{i=1}^{2}q^iu_h^{i}),\Pi_{k-1}^0v_h\Big)_E\Big|\Big\}\notag\\
&\leq\sum_E\Big\{\|\sum\limits_{i=1}^{2}q^iu^{i}-\Pi_{k-1}^0(\sum\limits_{i=1}^{2}q^iu^{i})\|_{0,E}\|v_h-\Pi_{k-1}^0v_h\|_{0,E}+
C\sum\limits_{i=1}^{2}\|u_h^{i}-u^{i}\|_{0,E}\|v_h\|_{0,E}\Big\}\notag\\
&\leq C\sum_E\Big\{h_E^{k+1}\sum\limits^{2}_{i=1}\|u^{i}\|_{k,E}|v_h|_{1,E}+\sum\limits_{i=1}^{2}\|u_h^{i}-u^{i}\|_{0,E}\|v_h\|_{0,E}\Big\}\notag\\
&\leq C(h^{k+1}\sum\limits^{2}_{i=1}\|u^{i}\|_{k}|v_h|_1+\sum\limits^{2}_{i=1}\|u_h^{i}-u^{i}\|_{0}\|v_h\|_0).
\end{align}
We complete the proof of this Lemma. $\hfill\Box$\\
\setcounter{lemma}{0}
\setcounter{theorem}{0}
\setcounter{corollary}{0}
\section{Error estimates for semidiscrete case}\setcounter{equation}{0}
In this section,  we present the a priori error estimates in the $L^2$ and $H^1$ norms for the semidiscrete system (\ref{semi form}). First, we deduce the error estimate in the $L^2$ norm. Then, we give some error estimates of an $H^1$ projection. At last, we present the priori error estimates in the $H^1$ norms.
\vspace{3mm}
\subsection{Error estimates for semidiscrete case in the $L^2$ norm}
In this subsection,  we present the a priori error estimates in the $L^2$ norm for the semidiscrete system (\ref{semi form}).
First, we develop the error analysis of $\phi(t)-\phi_h(t)$ in the $H^1$ norm.
\vspace{3mm}
\begin{lemma}\label{theorem3.2}
Let $(\phi,p^{i})$  and $(\phi_h,p_h^{i})$ be the solutions of (\ref{weak formulation}) and (\ref{semi form}),
respectively.
 Then for all $t\in (0,T]$, the following estimation holds
\begin{equation}
\|\phi(t)-\phi_h(t)\|_{1}\leq C\Big(h^k(||f||_{k}+\sum\limits^{2}_{i=1}\|p^{i}\|_k+\|\phi\|_{k+1})+\sum\limits^{2}_{i=1}\|p_h^{i}-p^{i}\|_0\Big).\notag
\end{equation}
\end{lemma}
\textbf{Proof.}~ Let $\;e_h=\phi_h(t)-\phi_I(t)$, where $\phi_I(t)\in Q_h^k$ is the interpolant of $\phi(t)$. Now, we eatimate $||\phi_h(t)-\phi_{I}(t)\|_1$. First,
\begin{align}\label{a1}
{a}_h(e_h,e_h)&={a}_h(\phi_h(t),e_h)-{a}_h(\phi_I(t),e_h)\notag\\
&=(f_h(t),e_h)-\widetilde{b}_h(p_h^{1}(t),p_h^{2}(t),e_h)-{a}_h(\phi_I(t),e_h).
\end{align}
From (\ref{ahcon}), it yields
\begin{align}\label{a2}
{a}_h(\phi_I(t),e_h)&=\sum_{E}\Big({a}_h^{E}(\phi_I(t)-\Pi_k^0{\phi(t)},e_h)+{a}_h^{E}(\Pi_k^0{\phi(t)},e_h)\Big)\notag\\
&=\sum_{E}\Big({a}_h^{E}(\phi_I(t)-\Pi_k^0{\phi(t)},e_h)+a^{E}(\Pi_k^0{\phi(t)},e_h)\Big)\notag\\
&=\sum_{E}\Big({a}_h^{E}(\phi_I(t)-\Pi_k^0{\phi(t)},e_h)+a^{E}(\Pi_k^0{\phi(t)}-\phi(t),e_h)+a^{E}(\phi(t),e_h)\Big)\notag\\
&=\sum_{E}\Big({a}_h^{E}(\phi_I(t)-\Pi_k^0{\phi(t)},e_h)+a^{E}(\Pi_k^0{\phi(t)}-\phi(t),e_h)\Big)+(f(t),e_h)-\widetilde{b}(p^{1}(t),p^{2}(t),e_h).
\end{align}
Then inserting (\ref{a2}) into (\ref{a1}), we deduce that
\begin{align}
{a}_h(e_h,e_h)&=\big(f_h(t)-f(t),e_h\big)-\sum_{E}\Big({a}_h^{E}\big(\phi_I(t)-\Pi_k^0{\phi(t)},e_h\big)+a^{E}\big(\Pi_k^0{\phi(t)}-\phi(t),e_h\big)\Big)\notag\\
&~~~~+\widetilde{b}\big(p^{1}(t),p^{2}(t),e_h\big)-\widetilde{b}_h\big(p_h^{1}(t),p_h^{2}(t),e_h\big)\notag\\
&\leq C\Big(h^k\left\|f(t)\right\|_k\left\|e_h\right\|_0
+\sum_{E}\big(\|\phi_I(t)-\phi(t)\|_{1,E}+\|\phi(t)-\Pi_k^0{\phi(t)}\|_{1,E}\big)\left\|e_h\right\|_{1}\notag\\
&\;\;\;\;+\big(h^{k+1}\sum\limits^{2}_{i=1}\|p^{i}(t)\|_k+\sum\limits^{2}_{i=1}\|p_h^{i}(t)-p^{i}(t)\|_0\big)\|e_h\|_1\Big)
~~(by ~~(\ref{continuity a}) ~~and~~ Lemma~~ \ref{lemma4.2})\notag\\
&\leq C\Big(h^k(\|f(t)\|_k+\sum\limits^{2}_{i=1}\|p^{i}(t)\|_k+\|\phi(t)\|_{k+1})+\sum\limits^{2}_{i=1}\|p_h^{i}(t)-p^{i}(t)\|_0\Big)
\|e_h\|_1,~~(by~~(\ref{3.14}))\notag
\end{align}
At last, from (\ref{continuity u}), we obtain
\begin{align}
C_0|e_h|_1^2&=C_0a(e_h,e_h)\leq{a}_h(e_h,e_h)\notag\\
&\leq C\Big(h^k(\|f(t)\|_k+\sum\limits^{2}_{i=1}\|p^{i}(t)\|_k+\|\phi(t)\|_{k+1})+\sum\limits^{2}_{i=1}\|p_h^{i}(t)-p^{i}(t)\|_0\Big)\|e_h\|_1,\notag
\end{align}
which together with the interpolation error estimate (\ref{3.14}) and Poincar$\acute{e}$ inequality, we prove the assertion of the theorem
\begin{equation}
\|\phi(t)-\phi_h(t)\|_{1}\leq C\Big(h^k(||f||_{k}+\sum\limits^{2}_{i=1}\|p^{i}\|_k+||\phi||_{k+1})+\sum\limits^{2}_{i=1}\|p_h^{i}-p^{i}\|_0\Big).\notag
\end{equation}
 This completes the proof.$\hfill\Box$\\
\vspace{3mm}

\par The following lemma which shall be used to present the error estimate in the $L^2$ norm.
\begin{lemma}\label{phiinfty} Suppose $(\phi,p^{i})$  and $(\phi_h,p_h^{i})$ are the solutions of (\ref{weak formulation}) and (\ref{semi form}), respectively. If $\mathcal{T}_h$ is quasi-uniform, then there holds
\begin{align*} \|\Pi_k^0\nabla\phi_h\|_{0,\infty,E}\leq Ch^{-1}\sum_{i=1}^2\|p^i-p^i_h\|_0.
\end{align*}
\end{lemma}
\textbf{Proof.}~ For any $q_k\in \mathbb{P}_k(E)$, from the standard inverse inequality, we have (see Section 2.3 in \cite{liu22})
  \begin{align*} \|q_k\|_{0,\infty,E}\leq Ch_E^{-1}\|q_k\|_{0,E}.
 \end{align*}
 Hence, it follows that
 \begin{align*} \|\Pi_k^0\nabla \phi_h\|_{0,\infty,E}&\leq \|\Pi_k^0\nabla \phi_h-\Pi_k^0\nabla \phi\|_{0,\infty,E}+\|\Pi_k^0\nabla \phi\|_{0,\infty,E}\notag\\
 &\leq Ch^{-1}_E\|\Pi_k^0\nabla \phi_h-\Pi_k^0\nabla \phi\|_{0,E}+C_1\notag\\
 &\leq C h^{-1}_E\|\nabla\phi-\nabla\phi_h\|_{0,E}.
 \end{align*}
Since $\mathcal{T}_h$ is quasi-uniform, we have $h\leq C h_E,~~\forall E\in \mathbb{T}_h$. Then from Lemma \ref{theorem3.2}, we get
\begin{align*}\|\Pi_k^0\nabla \phi_h\|_{0,\infty,E}&\leq Ch^{-1}\|\nabla\phi-\nabla\phi_h\|_0\notag\\
&\leq Ch^{-1}(h^k+\sum_{i=1}^2\|p^i-p^i_h\|_0)\notag\\
 &\leq Ch^{-1}\sum_{i=1}^2\|p^i-p^i_h\|_0.
\end{align*}
This completes the proof of this lemma.$\hfill\Box$\\

Next, we derive error estimates for $p^i$ in the $L^2$ norm. Assume
\begin{align}\label{reassu} &p^i\in L^{\infty}(0,T;H^{k+1}(\ome)\cap L^{\infty}(\ome)),~~p^i\in L^{\infty}(0,T;H^{k+1}(\ome))~~i=1,2 \notag\\
~~&and~~
\phi\in L^{\infty}(0,T;H^{k+1}(\ome)\cap W^{k+1,\infty}(\ome)).
\end{align}
We also suppose
\begin{align}\label{assuf} f\in L^{\infty}(0,T;H^{k}(\ome))~~and~~ F^i\in L^{\infty}(0,T;H^{k+1}(\ome)).
\end{align}
\begin{theorem}\label{theorem3.1}
 Suppose $\mathcal{T}_h$ is quasi-uniform. Let $(\phi,p^{i})$  and $(\phi_h,p_h^{i})$ be the solutions of (\ref{weak formulation}) and (\ref{semi form}), respectively, and set $p_{h,0}^i:=(p_0^i)_i$, the interpolant function of the initial value of $p^i_0$ in $Q_h^k$.
 Assume (\ref{reassu}) and (\ref{assuf}) holds, then for $t\in (0,T]$, there holds
\begin{align}
\sum\limits_{i=1}^{2}\|p_{h}^i( t)-p^i( t)\|_{0} \leq Ch^k,~~k\geq 1.\notag
\end{align}
\end{theorem}
\textbf{Proof.}~ Decompose the error as follows
\begin{eqnarray}
p_h^{i}(t)-p^{i}(t)=\big(p_h^{i}(t)-\Pi_{k}^0p^{i}(t)\big)+\big(\Pi_{k}^0p^{i}(t)-p^{i}(t)\big)=:\upsilon^{i}(t)+\varrho^{i}(t),\notag
\end{eqnarray}
which are then estimated separately. For the second term $\varrho^{i}(t)$, it is easy to get
\begin{align}\label{3.50}
\|\varrho^{i}(t)\|_0&=\|\Pi_{k}^0p^{i}(t)-p^{i}(t)\|_0\leq Ch^{k+1}|p^{i}(t)|_{k+1}.
\end{align}
Now, we proceed the eatimate for $\upsilon^{i}(t)$. For any $v_{h}\in Q_h^k$, an application of the first equation in (\ref{weak formulation}) together with (\ref{semi form}) yields
\vspace{-15mm}
\begin{spacing}{2.0}
\begin{align}\label{smh+ah}
&m_{h}\left(\upsilon_t^{i}(t), v_{h}\right)+a_{h}\left(\upsilon^{i}(t), v_{h}\right)\notag\\
&=\big(m_h(p_h^i(t),v_h)+a_h(p_h^i(t),v_h)\big)-m_{h}\left(\frac{d}{d t} \Pi_{k}^0 p^i( t), v_{h}\right)
-a_{h}\left(\Pi_{k}^0 p^i(t), v_{h}\right)\notag\\
&~~+a(p^i( t), v_{h})-a(p^i( t), v_{h})\notag\\
&=\big((F_{h}^i(t), v_{h})-b_{i,h}(p_h^i(t),\phi_h(t),v_h)\big)-m_{h}\left(\Pi_{k}^0 p^i_{t}( t), v_{h}\right)
-a_{h}\left(\Pi_{k}^0 p^i(t), v_{h}\right)\notag\\
&~~+a(p^i( t), v_{h}) -\big((F^i(t),v_h)-(p_t^i(t),v_h)-b_i(p^i(t),\phi(t),v_h)\big)\notag\\
&=\left(F_{h}^i(t)-F^i(t), v_{h}\right)+\big((p^i_{t}(t), v_{h})-m_{h}(\Pi_{k}^0 p^i_{t}( t), v_{h})\big)
+\left(a(p^i( t), v_{h})-a_h(\Pi_{k}^0p^i( t), v_{h})\right)\notag\\
&~~+\left(b_i(p^i(t),\phi(t),v_h)-b_{i,h}(p_h^i(t),\phi_h(t),v_h)\right)\notag\\
&:= H_1+H_2+H_3+H_4.
\end{align}
\end{spacing}
\vspace{-8mm}
\noindent The first term can be estimated as follows
\begin{eqnarray}\label{sh1}
\begin{split}
H_{1}=(F_{h}^i(t)-F^i(t),v_h)=(\Pi_{k}^0F^i(t)-F^i(t),v_h)\leq Ch^{k+1}\|F^i\|_{k+1}\|v_h\|_0.
\end{split}
\end{eqnarray}
The second term can be bounded by the consistency property (\ref{cons})
\begin{align}\label{sh2}
H_{2}=\left(p^i_{t}(t), v_{h}\right)-m_{h}\left(\Pi_{k}^0 p^i_{t}( t), v_{h}\right)&=\sum_{E}\Big(\left(p^i_{t}(t), v_{h}\right)_E-m_{h}^E\left(\Pi_{k}^0 p^i_{t}( t), v_{h}\right)\Big) \notag\\
&=\sum_{E}\Big(\left(p^i_{t}(t), v_{h}\right)_E-\left(\Pi_{k}^0 p^i_{t}( t), v_{h}\right)\Big) \notag\\
&\leq C \sum_{E}\|p^i_{t}( t)-\Pi_{k}^0 p^i_{t}( t)\|_{0,E}\|\| v_{h} \|_{0,E}  \notag\\
&\leq C h^{k+1}\|p_t^i(t)\|_{k+1}\| v_{h} \|_{0}.
\end{align}
\noindent From Lemma \ref{cons a}, we can express the third term $a(p^{i}(t),v_h)-a_h(\Pi_{k}^0p^{i}(t),v_h)$ as
\begin{align}\label{sh3}
H_{3}&=\sum_{E}\Big(a^E\big(p^i( t), v_{h}\big)-a_h^E\big(\Pi_{k}^0p^i( t), v_{h}\big)\Big)\notag \\
&=\sum_{E}\Big(a^E\big(p^i( t), v_{h}\big)-a^E\big(\Pi_{k}^0p^i( t), v_{h}\big)\Big)\notag\\
&\leq Ch^k\|p^i(t)\|_{k+1}\|\nabla v_h\|_{0}.
\end{align}
For the fourth term, from Lemmas \ref{bi-bhi}, \ref{theorem3.2} and \ref{phiinfty}, we have
\begin{align}\label{sh4}
H_{4}&=b_i(p^i(t),\phi(t),v_h)-b_{i,h}(p_h^i(t),\phi_h(t),v_h)
=\sum_{E}\big(b_i^E(p^i(t),\phi(t),v_h)-b_{i,h}^E(p_h^i(t),\phi_h(t),v_h)\big)\notag\\
&\leq \sum_{E}\Big(h_E^k+\|\nabla\phi(t)-\nabla\phi_h(t)\|_{0,E}+\|\Pi_k^0\nabla \phi_h\|_{0,\infty,E}\|p^i(t)-p_h^i(t)\|_{0,E}\Big)\|\nabla v_h\|_{0,E}\notag\\
&\leq C\Big(h^k+\sum\limits^{2}_{i=1}\|p_h^{i}(t)-p^{i}(t)\|_{0}+h^{-1}(\sum\limits^{2}_{i=1}\|p_h^{i}(t)-p^{i}(t)\|_{0})^2\Big)
\|\nabla v_h\|_{0}
\end{align}
Setting $v_h=\upsilon^{i}(t)$ in (\ref{smh+ah}) and using (\ref{sh1})-(\ref{sh4}), we get
\vspace{-2mm}
\begin{spacing}{1.0}
\begin{align}
&m_{h}\left(\upsilon^{i}_t(t),\upsilon^{i}(t)\right)+a_{h}(\upsilon^{i}(t), \upsilon^{i}(t))\notag\\
& \leq C h^{k+1}\|\upsilon^{i}(t)\|_{0}+C \Big(h^{k}+\sum\limits^{2}_{i=1}\|p_h^{i}(t)-p^{i}(t)\|_{0}+h^{-1}(\sum\limits^{2}_{i=1}\|p_h^{i}(t)-p^{i}(t)\|_{0})^2\Big)\|\nabla\upsilon^{i}(t)\|_{0}.\notag
\end{align}\end{spacing}
\vspace{0mm}
\noindent By introducing $\|v\|_{h}^{2}=m_{h}(v, v)~ \text{for all}~ v \in Q_{h}^k$, using (\ref{continuity u}) and the Poincar$\acute{e}$ inequality, we infer that\vspace{-13mm}
\begin{spacing}{2.0}
\begin{align}
&\frac{1}{2}\frac{d}{dt}\|\upsilon^{i}(t)\|_h^2+C_0\|\nabla\upsilon^{i}(t)\|_{0}^2\notag\\
& \leq C h^{k+1}\|\upsilon^{i}(t)\|_{0}+C\Big(h^{k}+\sum\limits^{2}_{i=1}\|p_h^{i}(t)-p^{i}(t)\|_{0}+h^{-1}(\sum\limits^{2}_{i=1}\|p_h^{i}(t)-p^{i}(t)\|_{0})^2\Big)\|\nabla\upsilon^{i}(t)\|_{0}.\notag\\
& \leq C \Big(h^{k}+\sum\limits^{2}_{i=1}\|\upsilon^i(t)\|_{0}+(h^{-1}\sum\limits^{2}_{i=1}\|\upsilon^i(t)\|_{0})\sum\limits^{2}_{i=1}\|\upsilon^i(t)\|_{0}\Big)\|\nabla\upsilon^{i}(t)\|_{0}.\notag\\
&\leq C \Big (h^{2k}+(\sum\limits^{2}_{i=1}\|\upsilon^i(t)\|_{0})^2+h^{-2}(\sum\limits^{2}_{i=1}\|\upsilon^i(t)\|_{0})^2\sum\limits^{2}_{i=1}\|\upsilon^i(t)\|_{0}^2\Big)+\frac{C_0}{2}\|\nabla\upsilon^{i}(t)\|_{0}^2.\notag
\end{align}
\end{spacing}
\vspace{-8mm}
 Hence
\vspace{0mm}
\begin{align}\label{svin}
\frac{d}{dt}\|\upsilon^{i}(t)\|_0^2&\leq C\Big(h^{2k}+(\sum\limits^{2}_{i=1}\|\upsilon^i(t)\|_{0})^2+h^{-2}(\sum\limits^{2}_{i=1}\|\upsilon^i(t)\|_{0})^2\sum\limits^{2}_{i=1}\|\upsilon^i(t)\|_{0}^2\Big).
\end{align}
\vspace{0mm}
Next, following the arguments in \cite{sun2016error} we shall show $h^{-1}\|\upsilon^i(t)\|_{0}\leq C$ by induction. First, from (\ref{3.14}) we have
\begin{align*}h^{-1}\|\upsilon^i(0)\|_{0}&=h^{-1}\|p_{h,0}^i-\Pi_{k}^0p^i_0\|_0\notag\\
&\leq h^{-1}(\|p_{h,0}^i-p^i_0\|_0+\|\varrho^{i}(0)\|_0)\leq Ch^{k+1}\|p^i_0\|_{k+1} \leq C.
\end{align*}
Then, assume $h^{-1}\|\upsilon^i(t)\|_{0}\leq C$ holds for $t\in [0,T_0],~~T_0< T$. From (\ref{svin}), we get
\begin{align*}
\frac{d}{dt}\|\upsilon^{i}(t)\|_0^2&\leq C h^{2k}+C\sum\limits^{2}_{i=1}\|\upsilon^i(t)\|_{0}^2.
\end{align*}
Integrating the above from 0 to $t$, the following inequality holds

\vspace{-2mm}
\begin{align}
\|\upsilon^{i}(t)\|_0^2&\leq \|\upsilon^{i}(0)\|_0^2
+Ch^{2k}+C\int_0^t\sum\limits^{2}_{i=1}\|\upsilon^i(s)\|_{0}^2ds.\notag
\end{align}
\vspace{-2mm}
Summing up for the index $i$, we get
\vspace{3mm}
\begin{align}
\sum\limits^{2}_{i=1}\|\upsilon^{i}(t)\|_0^2&\leq \sum\limits^{2}_{i=1}\|\upsilon^{i}(0)\|_0^2
+Ch^{2k}+C\int_0^t\sum\limits^{2}_{i=1}\|\upsilon^i(s)\|_{0}^2ds.\notag
\end{align}
\vspace{3mm}
By using Gronwall's inequality, we deduce that
\begin{align}
\sum\limits^{2}_{i=1}\|\upsilon^{i}(t)\|_0^2&\leq \sum\limits^{2}_{i=1}\|\upsilon^{i}(0)\|_0^2
+Ch^{2k}.\notag
\end{align}
That is
\begin{align}
\sum\limits^{2}_{i=1}\|\upsilon^{i}(t)\|_0&\leq \sum\limits^{2}_{i=1}\|\upsilon^{i}(0)\|_0
+Ch^{k}.\notag
\end{align}
From (\ref{3.14}), the term $\sum\limits^{2}_{i=1}\|\upsilon^{i}(0)\|_0$ can be estimated as follows
\begin{align}
\sum\limits^{2}_{i=1}\|\upsilon^{i}(0)\|_0&=\sum\limits^{2}_{i=1}\|p_{h,0}^i-\Pi_{k}^0p^i_0\|_0\leq\sum\limits^{2}_{i=1}(\|p_{h,0}^i-p^i_0\|_0+\|p^i_0-\Pi_{k}^0p^i_0\|_0)\notag\\
&\leq \sum\limits^{2}_{i=1}(\|p_{h,0}^i-p^i_0\|_0+C h^{k+1}|p^i_0|_{k+1})\leq C h^{k+1}|p^i_0|_{k+1}.\notag
\end{align}
Hence, we get if the assumption $h^{-1}\|\upsilon^i(t)\|_{0}\leq C$ holds for $t\in [0,T_0],~~T_0< T$ holds, then
\begin{align}\label{uphk}\sum\limits^{2}_{i=1}\|\upsilon^{i}(t)\|_0\leq C h^k,~~t\in [0,T_0].
\end{align}
It yields
\begin{align*}h^{-1}\|\upsilon^{i}(t)\|_0\leq C h^{k-1}\leq C,~~ for~~ k\geq 1~~ and~~ t\in [0,T_0].
\end{align*}
Since $h^{-1}\|\upsilon^{i}(t)\|_0$ is a continuous function with respect to $t\in [0,T]$, from the uniform continuity with time, there is $\delta$ such that $h^{-1}\|\upsilon^{i}(t)\|_0\leq C$ holds for $t\in [0, T_0+\delta]$. Since $[0,T]$ is a finite interval, we have
\begin{align*}h^{-1}\|\upsilon^{i}(t)\|_0\leq C,~~for~~ t\in [0,T].
\end{align*}
Thus, from (\ref{uphk}), we have
\begin{align} \sum\limits^{2}_{i=1}\|\upsilon^{i}(t)\|_0\leq C h^k,~~for~~ t\in [0,T].
\end{align}
Combining the estimates for $\upsilon^{i}(t)$ and $\varrho^{i}(t)$, we get
\begin{align}
\sum\limits_{i=1}^{2}\|p_{h}^i( t)-p^i( t)\|_{0}\leq \sum\limits^{2}_{i=1}(\|\varrho^{i}(t)\|_0+\|\upsilon^{i}(t)\|_0)\leq Ch^k.\notag
\end{align}
We complete the proof. $\hfill\Box$\\

Note that a suboptimal error estimate in the $L^2$ norm is presented for the virtual element solution $p_h^i$ in Theorem \ref{theorem3.1}. Since the PNP equations contain the nonlinear coupled term $p^i\nabla\phi$, it is not easy to present the optimal error estimate for the VEM. The similar difficulty needs be dealt with in the error analysis for the FEM (cf. \cite{sun2016error}), in which a suboptimal error estimate is obtained in the $L^2$ norm by using the standard arguments of analysis.
\begin{theorem}\label{theorem3.3}
Let $(\phi,p^{i})$ and $(\phi_h,p_h^{i})$ be the solutions of (\ref{weak formulation}) and (\ref{semi form}), respectively. Under the assumptions of Theorem \ref{theorem3.1}, for all $t\in(0,T]$, the following estimate holds
\begin{align}
\|\phi_h(t)-&\phi(t)\|_{0}\leq Ch^k.\notag
\end{align}
\end{theorem}
\textbf{Proof.}~  Let $\psi(t) \in H_{0}^{2}(\Omega)$ be the solution of following equation
\begin{equation}\label{phiw}
a(w,\psi(t))=(\phi(t)-\phi_h(t),w),\quad w\in Q_h^k(E).
\end{equation}
Then, the following regularity result can be obtained
\begin{equation}\label{spsi 1}
\|\psi(t)\|_{2} \leq C\left\|\phi(t)-\phi_h(t)\right\|_{0},
\end{equation}
where $C$ is a constant depending on the domain $\Omega$. 
\begin{spacing}{2.0}\vspace{-13mm}
\begin{align}\label{sphi}
\left\|\phi(t)-\phi_h(t)\right\|_{0}^{2} &=a\big(\phi(t)-\phi_h(t), \psi(t)\big) \notag\\
&=a\big(\phi(t)-\phi_h(t), \psi(t)-\psi_{I}(t)\big)+a\big(\phi(t), \psi_{I}(t)\big)-a\big(\phi_h(t),\psi_{I}(t)\big) \notag\\
&=a\big(\phi(t)-\phi_h(t), \psi(t)-\psi_{I}(t)\big)+\big(f(t), \psi_{I}(t)\big)-\tilde{b}(p^{1}(t), p^{2}(t), \psi_{I}(t))-a\big(\phi_h(t), \psi_{I}(t)\big)\notag\\
&~~+{a}_{h}\big(\phi_h(t), \psi_{I}(t)\big)+\tilde{b}_h\big(p_h^{1}(t), p_h^{2}(t), \psi_{I}(t)\big)
-\big(f_h(t),\psi_{I}(t)\big) \notag\\
&=a\big(\phi(t)-\phi_h(t),\psi(t)-\psi_{I}(t)\big)+\big(f(t)-f_h(t),\psi_{I}(t)\big)
+{a}_{h}\big(\phi_h(t),\psi_{I}(t)\big)\notag\\
&~~-a\big(\phi_h(t),\psi_{I}(t)\big)
+\tilde{b}_h\big(p_h^{1}(t),p_h^{2}(t),\psi_{I}(t)\big)-\tilde{b}\big(p^{1}(t),p^{2}(t),\psi_{I}(t)\big)\notag\\
&=a\big(\phi(t)-\phi_h(t),\psi(t)-\psi_{I}(t)\big)+\big(f(t)-f_h(t),\psi_{I}(t)-\Pi_{k}^0\psi_{I}(t)\big)\notag\\
&~~+\Big({a}_{h}\left(\phi_h(t)-\Pi_{k}^0\phi(t),\psi_{I}(t)\right)-a\left(\phi_h(t)-\Pi_{k}^0\phi(t),\psi_{I}(t)\right)\Big)\notag\\
&~~+\Big({a}_{h}\left(\Pi_{k}^0\phi(t),\psi_{I}(t)\right)-a\left(\Pi_{k}^0\phi(t),\psi_{I}(t)\right)\Big)\notag\\
&~~+\Big(\tilde{b}_h\left(p_h^{1}(t),p_h^{2}(t),\psi_{I}(t)\right)-\tilde{b}\left(p^{1}(t),p^{2}(t),\psi_{I}(t)\right)\Big)\notag\\
&:=M_1+M_2+M_3+M_4+M_5.
\end{align}\end{spacing}\vspace{-8mm}
\noindent From (\ref{3.14}), Lemma \ref{theorem3.2} and (\ref{spsi 1}), the first term can be estimated as follows\begin{spacing}{2.0}\vspace{-13mm}
\begin{align}\label{sm1}
M_1=a\big(\phi(t)-\phi_h(t),\psi(t)-\psi_{I}(t)\big)&\leq \|\phi(t)-\phi_h(t)\|_1\|\psi(t)-\psi_{I}(t)\|_1\notag\\
&\leq C \|\phi(t)-\phi_h(t)\|_1h|\psi(t)|_2\notag\\
&\leq C\Big(h^{k+1}\big(||f(t)||_{k}+\sum\limits^{2}_{i=1}\|p^{i}(t)\|_k+||\phi(t)||_{k+1}\big)\notag\\
&~~~~+h\sum\limits^{2}_{i=1}\|p_h^{i}(t)-p^{i}(t)\|_0\Big)\|\phi(t)-\phi_h(t)\|_0.
\end{align}\end{spacing}\vspace{-8mm}
\noindent For the second term, from (\ref{3.14}) and (\ref{spsi 1}), we have
\begin{align}\label{sm2}
M_2=\left(f(t)-f_h(t),\psi_{I}(t)-\Pi_{k}^0\psi_{I}(t)\right)&\leq \|f(t)-\Pi_{k}^0f(t)\|_0\|\psi_{I}-\Pi_{k}^0\psi_{I}(t)\|_0\notag\\
&\leq Ch^k|f(t)|_k h|\psi_{I}(t)|_1\notag\\
&\leq Ch^{k+1}|f(t)|_k\big(|\psi(t)-\psi_{I}(t)|_1+|\psi(t)|_1\big)\notag\\
&\leq Ch^{k+1}|f(t)|_k\|\phi(t)-\phi_h(t)\|_0.
\end{align}
To estimate $M_3$, setting $\bar{u}=\phi_h(t)-\Pi_{k}^0\phi(t)$ and $\bar{v}=\psi_{I}(t)$ for simplicity, the third term can be estimated with
\begin{align}\label{m3}
M_3&={a}_{h}(\bar{u},\bar{v})-a(\bar{u},\bar{v})\notag\\
&=\sum_E\Big((\Pi_{k-1}^0 \nabla\bar{u},\Pi_{k-1}^0\nabla \bar{v})_E+S^E_a\big((I-\Pi_k^\nabla)\bar{u},(I-\Pi_k^\nabla)\bar{v}\big)-(\nabla \bar{u}, \nabla \bar{v})_E\Big)\notag\\
&=\sum_E\Big(\big((\Pi_{k-1}^0 \nabla\bar{u},\Pi_{k-1}^0\nabla\bar{v})_E-(\nabla \bar{u}, \nabla \bar{v})_E\big)+S^E_a\big((I-\Pi_k^\nabla)\bar{u},(I-\Pi_k^\nabla)\bar{v}\big)\Big)\notag\\
&:=\sum_E(\tilde{M}_1+\tilde{M}_2).
\end{align}
Then
\begin{align}\label{tm1}
\tilde{M}_1&=(\Pi_{k-1}^0 \nabla\bar{u},\Pi_{k-1}^0\nabla\bar{v})_E-(\nabla \bar{u}, \nabla \bar{v})_E
=(\nabla \bar{u},\Pi_{k-1}^0 \nabla \bar{v})_E-(\nabla \bar{u}, \nabla \bar{v})_E\notag\\
&\leq \|\nabla \bar{u}\|_{0,E}\|\Pi_{k-1}^0 \nabla\bar{v}-\nabla \bar{v}\|_{0,E}\notag\\
&=\|\nabla \bar{u}\|_{0,E}\|\Pi_{k-1}^0 \nabla\psi_{I}(t)-\Pi_{k-1}^0 \nabla\psi(t)+\Pi_{k-1}^0 \nabla\psi(t)-\nabla\psi(t)+\nabla\psi(t)-\nabla\psi_{I}(t)\|_{0,E}\notag\\
&\leq Ch_E|\psi|_{2,E}\|\nabla\bar{u}\|_{0,E}\notag\\
&= Ch_E\|\nabla \big(\phi_h(t)-\Pi_{k}^0\phi(t)\big)\|_{0,E}|\psi(t)|_{2,E}.
\end{align}
And by using (\ref{sae}) and the standard approximation estimates, it holds (see (5.26) in \cite{beirao2016virtual})
\begin{align}\label{tm2}
\tilde{M}_2\leq Ch_E\|\nabla \bar{u}\|_{0,E}|\bar{v}|_{2,E}=Ch_E\|\nabla \big(\phi_h(t)-\Pi_{k}^0\phi(t)\big)\|_{0,E}|\psi(t)|_{2,E}.
\end{align}
Inserting (\ref{tm1}) and (\ref{tm2}) into (\ref{m3}), and from Lemma \ref{theorem3.2} and (\ref{spsi 1}), it follows that\vspace{3mm}
\begin{align}\label{sm3}
{M}_3&\leq Ch\sum_E\|\nabla \big(\phi_h(t)-\Pi_{k}^0\phi(t)\big)\|_{0,E}|\psi(t)|_{2,E}\notag\\
&=Ch\sum_E\|\nabla\big(\phi_h(t)-\phi(t)+\phi(t)-\Pi_{k}^0\phi(t)\big)\|_{0,E}\|\phi(t)-\phi_h(t)\|_{0,E}\notag\\
&\leq Ch\big(\|\phi(t)-\phi_h(t)\|_{1}+h^k\|\phi(t)\|_{k+1}\big)\|\phi(t)-\phi_h(t)\|_0\notag\\
&\leq C\Big(h^{k+1}\big(\|f(t)\|_k+\sum\limits_ {i=1}^{2}\|p^{i}(t)\|_k+\|\phi(t)\|_{k+1}\big)+h\sum\limits_ {i=1}^{2}\|p^{i}(t)-p_h^{i}(t)\|_0\Big)\|\phi(t)-\phi_h(t)\|_0.
\end{align}
From Lemma \ref{cons a}, we have
\begin{align}\label{sm4}
M_4={a}_{h}\left(\Pi_{k}^0\phi(t),\psi_{I}(t)\right)-a\left(\Pi_{k}^0\phi(t),\psi_{I}(t)\right)=0.
\end{align}
From (\ref{3.14}) and Lemma \ref{lemma4.2}, we can get
\begin{align}\label{sm5}
M_5&=\tilde{b}_h\left(p_h^{1}(t),p_h^{2}(t),\psi_{I}(t)\right)-\tilde{b}\left(p^{1}(t),p^{2}(t),\psi_{I}(t)\right)\notag\\
&\leq C\big(h^{k+1}\sum\limits^{2}_{i=1}\|p^{i}(t)\|_{k}|\psi_I(t)|_1+\sum\limits^{2}_{i=1}\|p^{i}(t)-p_h^{i}(t)\|_{0}\|\psi_I(t)\|_0\big)\notag\\
&\leq C\sum\limits^{2}_{i=1}\big(h^{k+1}\|p^{i}(t)\|_{k}+\|p^{i}(t)-p_h^{i}(t)\|_{0}\big)\|\phi(t)-\phi_{h}(t)\|_{0}.
\end{align}
By collecting (\ref{sm1})-(\ref{sm2}) and (\ref{sm3})-(\ref{sm5}) in (\ref{sphi}),
we get
\begin{equation}\label{sphi-sphih}
\|\phi(t)-\phi_{h}(t)\|_{0} \leq C\Big(h^{k+1}\big(\|f(t)\|_{k}+\sum_{i=1}^{2}\|p^{i}(t)\|_{k}+\|\phi(t)\|_{k+1}\big)+\sum_{i=1}^{2}\|p_{h}^{i}(t)-p^{i}(t)\|_{0}\Big).
\end{equation}
From Theorem \ref{theorem3.1}, we complete the proof of the theorem.$\hfill\Box$

\vspace{3mm}

\subsection{Error estimates of the $H^1$ projection}
In this subsection,  we deduce the error estimates of the $H^1$ projection, which shall be used in the a priori error estimates in the $H^1$ norm for the semidiscrete system (\ref{semi form}) in the next subsection.
Define the $H^1$ projection $R_h: H_0^1(\ome)\rightarrow Q_h^k$ satisfying: for any $u\in H_0^1(\ome)$ and at any given time $t\in [0,T]$,
\begin{align}\label{defrh}
a_h(R_hu,v_h)+b_{i,h}(R_hu,\phi(t),v_h)=a(u,v_h)+b_i(u,\phi(t),v_h),~~\forall v_h\in Q_h^k,
\end{align}
where $a(\cdot,\cdot),~~b_i(\cdot,\cdot,\cdot)$ and $a_h(\cdot,\cdot),~~b_{i,h}(\cdot,\cdot,\cdot)$ are defined in (\ref{2.4}) and (\ref{ahbh}), respectively.
 In order to present the existence and uniqueness of the solution of (\ref{defrh}), we need to show some lemmas, see Lemma \ref{piinfty}-\ref{lemmsup}. First the error estimate of the $L^2$ projection in the $L^{\infty}$ norm is presented as follows, which shall be used in Lemma \ref{Bih} later.
 \vspace{2mm}
\begin{lemma} \label{piinfty}If $w\in W^{k,\infty}(E)\cap H^{k+1}(E)$, then there holds
 \begin{align*} \|\Pi_k^0w-w\|_{0,\infty,E}\leq C h^k_E(\|w\|_{k,\infty,E}+\|w\|_{k+1,E}).
 \end{align*}
 \end{lemma}
 \textbf{Proof.}~ Setting $I_h:~H^1(E)\rightarrow S^h(E)$ is a piecewise polynomial interpolant, where $S^h(E)$ is the $k$-$th$ degree  finite element space defined on the element $E$. As mentioned in Section 3, for each $E\in \mathcal{T}^h$, there exists a virtual triangulation $\mathcal{T}_E$ of $E$ such that $\mathcal{T}_E$ is uniformly quasi-uniform (see \cite{chen2018}). Then by using the standard inverse inequality, for all $\tau\in \mathcal{T}_E$, there holds (cf...)
 \begin{align*} \|q\|_{0,\infty,\tau}\leq Ch^{-1}_E\|q\|_{0,\tau},~~\forall q\in \mathcal{P}_k(E),
 \end{align*}
 which implies
 \begin{align*} \|q\|_{0,\infty,E}\leq Ch^{-1}_E\|q\|_{0,E},~~\forall q\in \mathcal{P}_k(E).
 \end{align*}
 Hence,
 \begin{align*} \|\Pi_k^0w-I_hw\|_{0,\infty,E}&\leq C h^{-1}_E (\|\Pi_k^0w-w\|_{0,E}+\|w-I_hw\|_{0,E})\notag\\
 &\leq C h_{E}^{k}\|w\|_{k+1,E}.
 \end{align*}
 Then
 \begin{align*}\|\Pi_k^0w-w\|_{0,\infty,E}&\leq \|\Pi_k^0w-I_hw\|_{0,\infty,E}+\|I_hw-w\|_{0,\infty,E}\notag\\
 &\leq Ch^{k}_E(\|w\|_{k+1,E}+\|w\|_{k,\infty,E}),
 \end{align*}
  which finishes the proof of this lemma. $\hfill\Box$\\

For any $u,w,v\in H^1_0(\ome)$, set
 \begin{align*} B_i(u,w,v)=a(u,v)+b_i(u,w,v),
 \end{align*}
 where $a(\cdot,\cdot),~~b_i(\cdot,\cdot,\cdot)$ are defined as (\ref{2.4}).

  \begin{lemma} \label{Bi}Suppose $w\in W^{2,\infty}(\ome)$. There holds
  \begin{align}\label{supB} \sup_{v\in H_0^1} \frac {B_i(u,w,v)}{\|v\|_1}\geq C\|u\|_1,~i=1,2,~~\forall u\in H_0^1(\ome),
 \end{align}
 and
 \begin{align}\label{supB2} \sup_{u\in H_0^1} \frac {B_i(u,w,v)}{\|u\|_1}\geq C\|v\|_1,~i=1,2,~~\forall v\in H_0^1(\ome).
 \end{align}
 \end{lemma}
 \textbf{Proof.}~Consider the following problem
   \begin{equation}\label{steadpnp}
\left\{\begin{array}{lr}
Lu=-\nabla\cdot(\nabla{u}+q^i\nabla w u)=f,  ~~\text {in}~~ \Omega, \vspace{1mm}\\
u=0, \text { on} ~~\partial\Omega.\vspace{1mm}\\
\end{array}\right.
\end{equation}
If $w\in W^{2,\infty}(\Omega)$, then there exists a unique solution to (\ref{steadpnp}) and there holds (cf. Chapter 5 in \cite{sayas2019variational}).
\begin{align} \label{reg}\|u\|_1\leq C\|f\|_{-1}.
 \end{align}
By using (\ref{reg}), we have (cf.\cite{beirao2016virtual, liu22})
 \begin{align*}\sup_{v\in H_0^1(\ome)} \frac {B_i(u,w,v)}{\|v\|_1}=\sup_{v\in H_0^1(\ome)} \frac {<Lu,v>}{\|v\|_1}=\|Lu\|_{-1}=\|f\|_{-1}\geq C\|u\|_1,~~\forall u\in H_0^1(\ome).
 \end{align*}
 Similarly, if $w\in W^{2,\infty}(\ome)$, then there exists a unique solution to the following adjoint problem (see Chapter 5 in \cite{sayas2019variational})
 \begin{equation*}
\left\{\begin{array}{lr}
Lv=\nabla\cdot(\nabla{v}+q^{i}\nabla w\cdot\nabla v)=\tilde{f},  ~~\text {in}~~ \Omega, \vspace{1mm}\\
v=0, ~~ \text {on} ~~\partial\Omega,\vspace{1mm}\\
\end{array}\right.
\end{equation*}
 and
\begin{align*} \|v\|_1\leq C\|\tilde{f}\|_{-1}.
 \end{align*}
 Hence
 \begin{align*}\sup_{u\in H_0^1(\ome)} \frac {B_i(u,w,v)}{\|u\|_1}=\sup_{u\in H_0^1(\ome)} \frac {<Lv,u>}{\|u\|_1}=\|Lv\|_{-1}=\|\tilde{f}\|_{-1}\geq C\|v\|_1,~i=1,2,~~\forall v\in H_0^1(\ome).
 \end{align*}
 This completes the proof.$\hfill\Box$\\

For any $u\in H_0^1(\ome)$, $u_h$ and $w_h\in Q_h^k$, set
 \begin{align*}B_{i,h}(u_h,w,v_h)=a_h(u_h,v_h)+b_{i,h}(u_h,w,v_h),
 \end{align*}
 where $a_h(\cdot,\cdot)$ and $b_{i,h}(\cdot,\cdot,\cdot)$ are defined in (\ref{ahbh}). Then the definition of the $H^1$ projection (\ref{defrh}) can be written as
 \begin{align}\label{rhdef} B_{i,h}(R_hu,\phi(t),v_h)=B_i(u,\phi(t),v_h),~~\forall v_h\in Q_h^k.
 \end{align}
 If $\phi(t)|_E\in W^{1,\infty}(E)$, then the form $B_{i,h}$ is bounded, i.e.
 \begin{align}\label{boundness}B_{i,h}(u_h,\phi(t),v_h)\leq C \|u_h\|_1\|v_h\|_1.
 \end{align}

 From Lemma \ref{Bi}, we can get the following two lemmas for the discrete form $B_{i,h}$.
\vspace{2mm}
\begin{lemma} \label{Bih}Suppose $\phi\in W^{2,\infty}(E)\cap H^3(E)\cap Q_h^k$. There holds
 \begin{align*} \sup_{v_h\in Q_h^k} \frac {B_{i,h}(u_h,\phi(t),v_h)}{\|v_h\|_1} \geq C\|u_h\|_1,~~\forall u_h\in Q_h^k.
 \end{align*}
 \end{lemma}
 \textbf{Proof.}~  From (\ref{supB}), there exists $v\in H_0^1(\ome)$, such that
 \begin{align}\label{B_bound} B_i(u,\phi,v)\geq {C}_B \|u\|_1\|v\|_1,~~\forall u\in H_0^1(\ome).
 \end{align}
 From Lemma 5.6 in \cite{beirao2016virtual}, for any $v\in H_0^1(\ome)$, there exists a $v_h\in Q_h^k$, such that
 \begin{align}\label{ah} a_h(v_h,u_h)=a(v,u_h),~~\forall u_h\in Q_h^k,
 \end{align}
 and
 \begin{align}\label{v-vh} h\|v_h-v\|_1+\|v_h-v\|_0\leq Ch\|v\|_1,~~\|v_h\|_1\leq C\|v\|_1.
 \end{align}
 Then we have
 \begin{align}\label{bih2}
 B_{i,h}(u_h,\phi,v_h)&=a_h(u_h,v_h)+b_{i,h}(u_h,\phi,v_h)\notag\\
                      &=a(u_h,v)+b_{i,h}(u_h,\phi,v_h)-b_i(u_h,\phi,v)+b_i(u_h,\phi,v)\notag\\
                      &=B_i(u_h,\phi,v)+b_{i,h}(u_h,\phi,v_h)-b_i(u_h,\phi,v)\notag\\
                      &\geq {C}_B\|u_h\|_1\|v\|_1-|b_{i,h}(u_h,\phi,v_h)-b_i(u_h,\phi,v)|~~(by~~(\ref{B_bound}))\notag\\
                      &={C}_B\|u_h\|_1\|v\|_1-|b_{i,h}(u_h,\phi,v_h)-b_i(u_h,\phi,v_h)+b_i(u_h,\phi,v_h-v)|.
 \end{align}
 Next, we estimate $b_{i,h}(u_h,\phi,v_h)-b_i(u_h,\phi,v_h)$ and $b_i(u_h,\phi,v_h-v)$, respectively.
 \begin{align}\label{bih1} b_{i,h}(u_h,\phi,v_h)&-b_i(u_h,\phi,v_h)=\sum_E\big{\{} (q^i\Pi_{k-1}^0u_h\Pi_k^0\nabla\phi,\Pi_{k-1}^0\nabla v_h)_E-(q^iu_h\nabla\phi,\nabla v_h)_E\big{\}}\notag\\
 &=\sum_E\{(q^iu_h(\Pi_k^0\nabla\phi-\nabla\phi,\nabla v_h)_E+\big{(}(q^i\Pi_{k-1}^0u_h\Pi_k^0\nabla\phi,\Pi_{k-1}^0\nabla v_h)_E-(q^iu_h\Pi_k^0\nabla\phi,\nabla v_h)_E\big{)}\notag\\
 &=\sum_E\{I_1+I_2\}.
 \end{align}
 From Lemma \ref{piinfty}, we have
 \begin{align} \label{diff1}I_1&=(q^iu_h(\Pi_k^0\nabla\phi-\nabla\phi,\nabla v_h)_E\notag\\
        &=(q^i(\Pi_{k-1}^0u_h-u_h)\nabla\phi+\Pi_{k-1}^0u_h(\Pi_k^0\nabla\phi-\nabla\phi)+(u_h-\Pi_{k-1}^0u_h)\Pi_k^0\nabla\phi,\nabla v_h)_E\notag\\
        &\leq C(\|\Pi_{k-1}^0u_h-u_h\|_{0,E}\|\nabla\phi\|_{0,\infty,E}+\|\Pi_{k-1}^0u_h\|_{0,E}\|\Pi_k^0\nabla\phi-\nabla\phi\|_{0,\infty,E}
        )\|v_h\|_{1,E}\notag\\
        &\leq Ch_E\|u_h\|_{1,E}\|v_h\|_{1,E}.
 \end{align}
 Following the deduction of (\ref{I11}), we get
 \begin{align*}\label{diff2}I_2&=(q^i\Pi_{k-1}^0u_h\Pi_k^0\nabla\phi,\Pi_{k-1}^0\nabla v_h)_E-(q^iu_h\Pi_k^0\nabla\phi,\nabla v_h)_E\notag\\
& \leq C(\|\beta u_h-\Pi_{k-1}^0(\beta u_h)\|_{0,E}+h_E\|u_h\|_{1,E})\|v_h\|_{1,E},
 \end{align*}
 where $\beta=q^i\Pi_{k}^0\nabla\phi$. Taking $u=u_h$ and $\hat{\beta}=q^i\nabla\phi$ in (\ref{hatu-u}), and using Lemma \ref{piinfty}, then we deduce that
 \begin{align}\|\beta u_h-\Pi_{k-1}^0(\beta u_h)\|_{0,E}&\leq C\|(\beta-\hat{\beta})u_h\|_{0,E}+Ch_E\|\nabla\phi u_h\|_{1,E}~~(By~~following ~~(\ref{hatu-u}))\notag\\
 &\leq \|\beta-\hat{\beta}\|_{0,\infty,E}\|u_h\|_{1,E}+Ch_E\|\nabla\phi u_h\|_{1,E}\notag\\
 &\leq Ch_E\|u_h\|_{1,E}.
 \end{align}
 Hence
 \begin{align}\label{diff2}I_2\leq h_E\|u_h\|_{1,E}\|v_h\|_{1,E}.
 \end{align}
 Inserting (\ref{diff1}) and (\ref{diff2}) into (\ref{bih1}), we get
 \begin{align}\label{diff3} b_{i,h}(u_h,\phi,v_h)-b_i(u_h,\phi,v_h)\leq C_{1} h\|u_h\|_1\|v_h\|_1.
 \end{align}
 Now, we estimate the term $b_i(u_h,\phi,v_h-v)$. There holds
 \begin{align}\label{diff4} b_i(u_h,\phi,v_h-v)&=\big{(}q^iu_h\nabla\phi,\nabla(v_h-v)\big{)}\notag\\
 &=-\big{(}div(q^iu_h\nabla\phi),v_h-v\big{)}+\int_{\partial\ome} q^i u_h\nabla\phi\cdot n(v_h-v)\notag\\
 &\leq C \|div(q^iu_h\nabla\phi)\|_0\|v-v_h\|_0\notag\\
&\leq C_{2} h\|u_h\|_1\|v\|_1~~(by~(\ref{v-vh})).
 \end{align}
 Substituting (\ref{diff3}) and (\ref{diff4}) into (\ref{bih2}), it yields
 \begin{align*} B_{i,h}(u_h,\phi,v_h)&\geq \tilde{C}_B \|u_h\|_1\|v\|_1\geq {C} \|u_h\|_1\|v_h\|_1,
 \end{align*}
 where $h(C_{1}+C_{2})\leq \frac {\tilde{C}_B} 2$. This completes the proof.$\hfill\Box$\\

\vspace{2mm}
\begin{lemma}\label{lemmsup} Assume $\phi\in W^{2,\infty}(E)\cap H^3(E)\cap Q_h^k$. For any $v_h\in Q_{h}^{k}$, we have
\begin{align*} \sup_{u_h\in Q_h^k}|B_{i,h}(u_h,\phi,v_h)|>0,~~\forall v_h\neq 0.
\end{align*}
\end{lemma}
\textbf{Proof.}~ We shall follow Lemma \ref{Bih} to present the result of this lemma. From (\ref{supB2}), we know that there exists a $u\in H_0^1(\ome)$, such that
\begin{align}\label{bi2} B_i(u,\phi,v)\geq C \|u\|_1\|v\|_1,~~\forall v\in H_0^1(\ome).
\end{align}
From (\ref{ah}), we get for any $u\in H_0^1(\ome)$, there exists a $u_h\in Q_h^k$ such that
 \begin{align} \label{ah2}a_h(u_h,v_h)=a(u,v_h),~~\forall v_h\in Q_h^k,
 \end{align}
 and
\begin{align}\label{uh-u} h\|u_h-u\|_1+\|u_h-u\|_0\leq Ch\|u\|_1,~~\|u_h\|_1\leq C\|u\|_1.
 \end{align}
 By following the arguments of (\ref{bih2}) and using (\ref{bi2})-(\ref{ah2}), for $u\in H_0^1(\ome)$, there exists a $u_h\in Q_h^k$, satisfying
 \begin{align}\label{Bih2}B_{i,h}(u_h,\phi,v_h)&=a_h(u_h,v_h)+b_{i,h}(u_h,\phi,v_h)\notag\\
 &=a(u,v_h)+b_{i,h}(u_h,\phi,v_h)-b_i(u,\phi,v_h)+b_i(u,\phi,v_h)\notag\\
 &=B_i(u,\phi,v_h)+b_{i,h}(u_h,\phi,v_h)-b_i(u,\phi,v_h)\notag\\
 &\geq C \|u\|_1\|v_h\|_1-|b_{i,h}(u_h,\phi,v_h)-b_i(u_h,\phi,v_h)+b_i(u_h-u,\phi,v_h)|,~~\forall v_h\in Q_h^k.
 \end{align}
 A difference between (\ref{bih2}) and the above is that $u_h\in Q_h^k$ is an arbitrary function in (\ref{bih2}) while it is
 a function dependent on $u$ in the above. Then from (\ref{diff3}), we have
 \begin{align*}b_{i,h}(u_h,\phi,v_h)-b_i(u_h,\phi,v_h)\leq Ch\|u_h\|_1\|v_h\|_1,
 \end{align*}
 and by (\ref{uh-u}),
 \begin{align*}b_i(u_h-u,\phi,v_h)&\leq \|u_h-u\|_0\|\phi\|_{1,\infty}\|v_h\|_1\notag\\
 & \leq Ch \|u\|_1\|v_h\|_1,
 \end{align*}
 which combining with (\ref{Bih2}) yields
 \begin{align} \label{Bih3}B_{i,h}(u_h,\phi,v_h)\geq C \|u\|_1\|v_h\|_1\geq C \|u_h\|_1\|v_h\|,~~\forall v_h\in Q_h^k,
 \end{align}
where we have used (\ref{uh-u}). Since $u_h$ is the solution of (\ref{ah2}), it can not be a trivial solution. Thus, from (\ref{Bih3}), we have there exists $u_h\neq 0$, such that
 \begin{align*}
 B_{i,h}(u_h,\phi,v_h)> 0,~~\forall v_h\neq 0,
 \end{align*}
 which leads
 \begin{align*}
 \sup_{u_h\in Q_h^k}|B_{i,h}(u_h,\phi,v_h)|>0,~~\forall v_h\neq 0.
 \end{align*}
 This completes the proof of this lemma. $\hfill\Box$\\

Now we can show the existence and uniqueness of the solution of (\ref{defrh}) as follows.
\begin{lemma} \label{exi}If $\phi\in W^{2,\infty}(E)\cap H^3(E)\cap Q_h^k$, then there exists a unique solution $R_hu$ satisfying (\ref{defrh}).
\end{lemma}
\textbf{Proof.}~ It is easy to get the result of Lemma \ref{exi} by using generalized Lax-Milgramm Lemma (see \cite{ga14}), since the discrete form $B_{i,h}$ satisfies all the conditions of generalized Lax-Milgramm Lemma from (\ref{boundness}) and Lemmas \ref{Bih}-\ref{lemmsup}. $\hfill\Box$\\

Next, we shall presented the error estimates for the $H^1$ projection $R_h$.
\begin{lemma}\label{rhl2err}Suppose $u\in H^{k+1}(\ome)$ and $\phi\in W^{2,\infty}(\ome)\cap H^{k+2}(\ome)$. Then the following estimate for the projection $R_h$ holds.
\begin{align*}
\|R_hu-u\|_0+\|\nabla(R_hu-u)\|_0\leq Ch^k.
\end{align*}
\end{lemma}
\textbf{Proof.}~ Suppose $u_I\in Q_h^k$ is the the interpolant to $u\in H^1_0(\ome)$ and $u_{\pi}\in P_k(E)$. From (\ref{sta}) and (\ref{defrh}), we have
\begin{align}\label{b12}
C_0\|\nabla(R_hu-u_I)\|_0^2&=C_0a(R_hu-u_I,R_hu-u_I)\leq a_h(R_hu-u_I,R_hu-u_I)\notag\\
&=a_h(R_hu,R_hu-u_I)-a_h(u_I,R_hu-u_I)\notag\\
&=a(u,R_hu-u_I)+b_i(u,\phi,R_hu-u_I)-b_{i,h}(R_hu,\phi,R_hu-u_I)\notag\\
&~~~~-\sum_E\big{(}a_h^E(u_I-u_{\pi},R_hu-u_I)+a_h^E(u_{\pi},R_hu-u_I)\big{)}\notag\\
&=-\sum_E\big{(}a_h^E(u_I-u_{\pi},R_hu-u_I)+a^E(u_{\pi}-u,R_hu-u_I)\big{)} \notag\\
&~~+\Big(b_i(u,\phi,R_hu-u_I)-b_{i,h}(R_hu,\phi,R_hu-u_I)\Big)\notag\\
&:=B_1+B_2.
\end{align}
Using (\ref{3.14}), (\ref{vpi}) and (\ref{continuity a}), it yields
\begin{align}\label{b1}
B_1&=-\sum_E\big{(}a_h^E(u_I-u_{\pi},R_hu-u_I)+a^E(u_{\pi}-u,R_hu-u_I)\big{)}\notag\\
&\leq C\sum_E(\|\nabla(u_I-u)\|_{0,E}+\|\nabla(u_{\pi}-u)\|_{0,E})\|\nabla(R_hu-u_I)\|_{0,E}\notag\\
&\leq C h^{k}_E\|\nabla(R_hu-u_I)\|_{0}.
\end{align}

To estimate $B_2$, from (\ref{I1}) and the first inequality in (\ref{I11}), for any $w,~v_h\in H^1(\ome)$ and $u_h\in Q_h^k$ we have
\begin{align}\label{bi-bih2}
& b_i^E(u_h,w,v_h)-b_{i,h}^E(u,w,v_h)=(q^iu_h\nabla w,\nabla v_h)_E-(q^i\Pi_{k-1}^0u_h\Pi_{k}^0\nabla w,\Pi_{k-1}^0\nabla v_h)_E\notag\\
&=(q^iu_h(\nabla w-\Pi_k^0\nabla w),\nabla v_h)_E+\Big((q^iu_h\Pi_k^0\nabla w,\nabla v_h)_E-(q^i\Pi_{k-1}^0u_h\Pi_{k}^0\nabla w,\Pi_{k-1}^0\nabla v_h)_E\Big)\notag\\
&=\tilde{B}_{21}+\tilde{B}_{22}.
\end{align}
For $w\in W^{1,\infty}(E)$ and $u\in L^{\infty}(E)$, there holds
\begin{align}\label{b21}\tilde{B}_{21}&=(q^iu_h(\nabla w-\Pi_k^0\nabla w),\nabla v_h)_E\notag\\
&=\Big( q^i\big( (u_h-u)\nabla w+u(\nabla w-\Pi_k^0\nabla w)+(u-u_h)\Pi_k^0\nabla w\big),\nabla v_h\Big)_E\notag\\
&\leq C\Big(\|u-u_h\|_{0,E}\|\nabla w\|_{0,\infty,E}+\|u\|_{0,\infty,E}\|\nabla w-\Pi_k^0\nabla w\|_{0,\infty,E}\Big)\|\nabla v_h\|_{0,E}\notag\\
&\leq C\Big(\|u-u_h\|_{0,E}+h^k\|w\|_{k+1,E}\Big)\|\nabla v_h\|_{0,E}.
\end{align}
And setting $\beta=q^i\Pi_{k}^0\nabla w$, from (\ref{3.23}) we have
\begin{align}\label{b22}\tilde{B}_{22}&= (q^iu_h\Pi_{k-1}^0\nabla w,\nabla v_h)_E-(q^i\Pi_{k-1}^0u_h\Pi_{k-1}^0\nabla w,\Pi_{k-1}^0\nabla v_h)_E\notag\\
&\leq \|\beta u_h-\Pi_{k-1}^0(\beta u_h)\|_{0,E}\|\nabla v_h-\Pi_{k-1}^0\nabla v_h\|_{0,E}\notag\\
&\;\;+\|\beta\cdot\nabla v_h-\Pi_{k-1}^0(\beta\cdot\nabla v_h)\|_{0,E}\|u_h-\Pi_{k-1}^0u_h\|_{0,E}\notag\\
&\;\;+C_{\beta}\|u_h-\Pi_{k-1}^0u_h\|_{0,E}\|\nabla v_h-\Pi_{k-1}^0\nabla v_h\|_{0,E}\notag\\
& \leq C(\|u_h-\Pi_{k-1}^0u_h\|_{0,E}+\|\beta u_h-\Pi_{k-1}^0(\beta u_h)\|_{0,E})\|\nabla v_h\|_{0,E},
\end{align}
Combining (\ref{b21})-(\ref{b22}) with (\ref{bi-bih2}), we have
\begin{align}\label{bi-bih3}b_i^E(u_h,w,v_h)-b_{i,h}^E(u,w,v_h)\leq &\Big(\|u-u_h\|_{0,E}+\|u_h-\Pi_{k-1}^0u_h\|_{0,E}\notag\\
&~~~~+\|\beta u_h-\Pi_{k-1}^0(\beta u_h)\|_{0,E}+h^k\|w\|_{k+1,E}\Big)\|\nabla v_h\|_{0,E}
\end{align}
Taking $u_h=R_hu$, $w=\phi$ and $v_h=R_hu-u_I$ in (\ref{bi-bih3}), we have
\begin{align}\label{3.49}
&b_i(R_hu,\phi,R_hu-u_I)-b_{i,h}(R_hu,\phi,R_hu-u_I)\notag\\
&\leq \sum_E C\Big(h_E^k+\|u-R_hu\|_{0,E}+\|R_hu-\Pi_{k-1}^0R_hu\|_{0,E}+\|\beta R_hu-\Pi_{k-1}^0(\beta R_hu)\|_{0,E}\Big)\|\nabla R_hu-u_I\|_{0,E}\notag\\
&\leq C(\|u-R_hu\|_0+h^k)\|\nabla(R_hu-u_I)\|_0.
\end{align}
Hence,
\begin{align}\label{b2}
B_2&=b_i(u,\phi,R_hu-u_I)-b_{i,h}(R_hu,\phi,R_hu-u_I)\notag\\
&=b_i(u,\phi,R_hu-u_I)-b_i(R_hu,\phi,R_hu-u_I)+b_i(R_hu,\phi,R_hu-u_I)-b_{i,h}(R_hu,\phi,R_hu-u_I)\notag\\
&\leq C(\|u-R_hu\|_0+h^k)\|\nabla(R_hu-u_I)\|_0.
\end{align}
Substituting (\ref{b1}) and (\ref{b2}) into (\ref{b12}), it yields
\begin{align}\label{rh-ui}
\|\nabla(R_hu-u_I)\|_0\leq C(\|u-R_hu\|_0+h^k).
\end{align}
\par Next, we present the $L^2$ estimate for $R_h-u$. Define the adjoint problem as follows:
\begin{small}
\begin{equation}\label{ad-prob}
\left\{\begin{array}{lr}
-\Delta w^i+q^i\nabla\phi\cdot\nabla w^i=u-R_hu,~~x\in \ome,\vspace{2mm}\\
w^i=0,~~x\in \partial\ome.\end{array}\right.
\end{equation}
\end{small}
If $\phi\in W^{1,\infty}(\ome)$, then the regularity result holds (cf. \cite{sun2016error})
\begin{align}\label{reg-adj}\|w^i\|_2\leq C\|u-R_hu\|_0.
\end{align}
From (\ref{defrh}) and (\ref{ad-prob}), we have
\begin{align}\label{d1234}\|u-R_hu\|_0^2&=(-\Delta w^i,u-R_hu)+(q^i\nabla\phi\cdot\nabla w^i,u-R_hu)\notag\\
&=a(u-R_hu,w^i)+b_i(u-R_hu,\phi,w^i)\notag\\
&=a(u-R_hu,w^i-w_I^i)+a(u-R_hu,w^i_I)+b_i(u-R_hu,\phi,w^i-w^i_I)\notag\\
&~~+b_i(u-R_hu,\phi,w^i_I)+a_h(R_hu,w^i_I)-a_h(R_hu,w^i_I)\notag\\
&~~+b_{i,h}(R_hu,\phi,w^i_I)-b_{i,h}(R_hu,\phi,w^i_I)\notag\\
&=a(u-R_hu,w^i-w_I^i)+\big{(}a_h(R_hu,w^i_I)-a(R_hu,w^i_I)\big{)}\notag\\
&~~+b_i(u-R_hu,\phi,w^i-w^i_I)+\big{(}b_{i,h}(R_hu,\phi,w^i_I)-b_{i}(R_hu,\phi,w^i_I)\big{)}~~(by~(\ref{defrh}))\notag\\
&:=D_1+D_2+D_3+D_4.
\end{align}
From (\ref{reg-adj}), we get
\begin{align}\label{d1}
D_1=a(u-R_hu,p^i-p_I^i)&\leq C\|\nabla(u-R_hu)\|_0\|\nabla(w^i-w_I^i)\|_0\notag\\
&\leq Ch\|\nabla(u-R_hu)\|_0\|w^i\|_2\notag\\
&\leq Ch\|\nabla(u-R_hu)\|_0\|u-R_hu\|_0.
\end{align}
From Lemma 3.2, we have
\begin{align}
D_2=&a_h(R_hu,w^i_I)-a(R_hu,w^i_I)\notag\\
&=\sum_E\{ a_h(R_hu-u_{\pi},w^i_I-\Pi_k^0w^i)-a(R_hu-u_{\pi},w^i_I-\Pi_k^0w^i)\}\notag\\
&\leq \sum_E\|\nabla(R_hu-u_{\pi})\|_{0,E}\|\nabla(w^i_I-\Pi_k^0w^i)\|_{0,E}\notag\\
&\leq C(h^k\|u\|_{k+1}+\|\nabla(R_hu-u)\|_{0})h\|w^i\|_2\notag\\
&\leq C(h^{k+1}+h\|\nabla(R_hu-u)\|_{0})\|R_hu-u\|_0.
\end{align}
and
\begin{align}\label{d3}
D_3=b_i(u-R_hu,\phi,w^i-w^i_I)&\leq C\|\nabla(R_hu-u)\|_{0}\|\nabla(w^i-w^i_I)\|_0\notag\\
&\leq C\|\nabla(R_hu-u)\|_{0}h\|w^i\|_2\notag\\
&\leq Ch\|\nabla(R_hu-u)\|_{0}\|R_hu-u\|_0.
\end{align}
It remains to estimate $D_4$. There holds
\begin{align}\label{d4}
D_4&=b_{i,h}(R_hu,\phi,w^i_I)-b_{i}(R_hu,\phi,w^i_I)\notag\\
&=\sum_E\Big(q^i(\Pi_{k-1}^0(R_hu)\Pi_k^0\nabla\phi,\Pi_{k-1}^0\nabla w_I^i)_{0,E}-q^i(R_hu\nabla\phi,\nabla w^i_I)_{0,E}\Big)\notag\\
&=\sum_E\Big(\big(q^i(\Pi_{k-1}^0(R_hu)\Pi_k^0\nabla\phi,\Pi_{k-1}^0\nabla w_I^i)_{E}-q^i(R_hu\Pi_k^0\nabla\phi,\nabla w^i_I)_{0,E}\big)\notag\\
&~~~~+q^i(R_hu(\Pi_k^0\nabla\phi-\nabla\phi),\nabla w_I^i)_{E}\Big)\notag\\
&:=\sum_E(D_{41}+D_{42}).
\end{align}
Setting $\tilde{\beta}=q^i\Pi_k^0\phi$ and $u_h=R_hu$, from (\ref{3.23}) it yields
\begin{align}\label{d41}D_{41}&=\sum_E\Big(q^i(\Pi_{k-1}^0(R_hu)\Pi_k^0\nabla\phi,\Pi_{k-1}^0\nabla w_I^i)_{0,E}-q^i(R_hu\Pi_k^0\nabla\phi,\nabla w^i_I)_{0,E}\Big)\notag\\
&=(\tilde{\beta}\Pi_{k-1}^0u_h,\Pi_{k-1}^0\nabla w_I^i)_E-(\tilde{\beta}u_h,\nabla w_I^i)_E\notag\\
&\leq \|\tilde{\beta}\cdot\nabla w_I^i-\Pi_{k-1}^0(\tilde{\beta}\cdot\nabla w_I^i)\|_{0,E}\|u_h-\Pi_{k-1}^0u_h\|_{0,E}+\|\nabla w_I^i-\Pi_{k-1}^0(\nabla w_I^i)\|_{0,E}\|\tilde{\beta}u_h-\Pi_{k-1}^0(\tilde{\beta}u_h)\|_{0,E}\notag\\
&~~+C\|\nabla w_I^i-\Pi_{k-1}^0\nabla w_I^i\|_{0,E}\|u_h-\Pi_{k-1}^0u_h\|_{0,E}.
\end{align}
Note that
\begin{align}\label{d41-1}
\|u_h-\Pi_{k-1}^0u_h\|_{0,E}&=\|R_hu-\Pi_{k-1}^0(R_hu)\|_{0,E}\notag\\
&=\|R_hu-u\|_{0,E}+\|u-\Pi_{k-1}^0u\|_{0,E}+\|\Pi_{k-1}^0u-\Pi_{k-1}^0R_hu\|_{0,E}\notag\\
&\leq \|R_hu-u\|_{0,E}+Ch^k\|u\|_{k,E},
\end{align}
\begin{align}\label{d41-2}
\|\tilde{\beta}\cdot\nabla w_I^i-\Pi_{k-1}^0(\tilde{\beta}\cdot\nabla w_I^i)\|_{0,E}&\leq \|\tilde{\beta}\cdot(\nabla w_I^i-\nabla w^i)\|_{0,E}
+\|\tilde{\beta}\nabla w^i-\Pi_{k-1}^0(\tilde{\beta}\cdot\nabla w^i)\|_{0,E}\notag\\
&~~~~+\|\Pi_{k-1}^0\tilde{\beta}\cdot(\nabla w^i-\nabla w_I^i)\|_{0,E}\notag\\
&\leq Ch\|w^i\|_{2,E},
\end{align}
and
\begin{align}\label{d41-3}
\|\tilde{\beta}u_h-\Pi_{k-1}^0(\tilde{\beta}u_h)\|_{0,E}&=\|\tilde{\beta}R_hu-\Pi_{k-1}^0(\tilde{\beta}R_hu)\|_{0,E}\notag\\
&\leq \|\tilde{\beta}(R_hu-u)\|_{0,E}+\|\tilde{\beta}u-\Pi_{k-1}^0(\tilde{\beta}u)\|_{0,E}+\|\Pi_{k-1}^0\tilde{\beta} u-\tilde{\beta}R_hu)\|_{0,E}\notag\\
&\leq C(\|R_hu-u\|_{0,E}+h^k\|u\|_{k,E}). ~~(by~~(\ref{hatu-u})~~ with~~ w=\phi)
\end{align}
Combining (\ref{d41-1})-(\ref{d41-3}) with (\ref{d41}), we get
\begin{align}\label{d41-4}
D_{41}\leq Ch\|w^i\|_{2,E}(\|R_hu-u\|_{0,E}+h^k\|u\|_{k,E})\leq C(h\|R_hu-u\|_{0,E}^2+h^{k+1}\|R_hu-u\|_{0,E}).
\end{align}
Now, we estimate $D_{42}$. From Lemma \ref{piinfty}, it follows that
\begin{align}\label{d42}
D_{42}&=q^i(R_hu(\Pi_k^0\nabla\phi-\nabla\phi),\nabla w_I^i)_{0,E}\notag\\
&=q^i\Big((R_hu-u)(\Pi_k^0\nabla\phi-\nabla\phi),\nabla w_I^i)_{0,E}+q^i\Big(u(\Pi_k^0\nabla\phi-\nabla\phi),\nabla w_I^i)_{0,E}\notag\\
&\leq C\Big((\|R_hu-u\|_{0,E}\|\Pi_k^0\nabla\phi-\nabla\phi\|_{0,\infty,E}\|\nabla w_I^i\|_{0,E}+
\|u\|_{0,\infty,E}\|\Pi_k^0\nabla\phi-\nabla\phi\|_{0,E})\|\nabla w_I^i\|_{0,E}\Big)\notag\\
& \leq C(h^k \|R_hu-u\|_{0,E}^2+h^k\|R_hu-u\|_{0,E}).
\end{align}
Substituting (\ref{d41-4}) and (\ref{d42}) into (\ref{d4}), we deduce
\begin{align}\label{d4-1}
D_4&=b_{i,h}(R_hu,\phi,w_I^i)-b_i(R_hu,\phi,w_I^i)\notag\\
&\leq C(h\|R_hu-u\|_{0,E}^2+h^k\|R_hu-u\|_{0,E}).
\end{align}

Inserting (\ref{d1})-(\ref{d4-1}) to (\ref{d1234}) and using (\ref{rh-ui}), we have
\begin{align}
\|u-R_hu\|_0^2&\leq C(h^k\|R_hu-u\|_0+h\|R_hu-u\|_0^2).
\end{align}
Hence, if $h$ is small enough, then it follows that
\begin{align}\label{rhl2}
\|u-R_hu\|_0\leq C h^k.
\end{align}
Combining the above inequality with (\ref{rh-ui}), it yields
\begin{align}
\|\nabla(u-R_hu)\|_0\leq C h^k,
\end{align}
which finishes the proof of this lemma. $\hfill\Box$

\begin{lemma} \label{rhh1}Suppose $u,~\partial_t u\in H^{k+1}(\ome)$, $\phi\in W^{2,\infty}(\ome)\cap H^{k+2}(\ome)$ and $\partial_t\phi(t)\in W^{1,\infty}(\ome)$. There holds
\begin{align}
\|\nabla\partial_t(R_hu-u)\|_0\leq C(h^k+\|\partial_t(R_hu-u)\|_0).
\end{align}
\end{lemma}
\textbf{Proof.}~  Since
\begin{align*}\|\nabla\partial_t(R_hu-u)\|_0&\leq \|\nabla\partial_t(R_hu-\Pi_k^0u)\|_0+\|\nabla\partial_t(\Pi_k^0u-u)\|_0\\
&\leq \|\nabla\partial_t(R_hu-\Pi_k^0u)\|_0+Ch^k\|\partial_t u\|_{k+1},
\end{align*}
it suffices to present $\|\nabla\partial_t(R_hu-\Pi_k^0u)\|_0$. Integrating on both sides of (\ref{defrh}) with respecting to $t$, we have
 \begin{align}\label{3.74}&a(\partial_tu,v_h)+a(u,\partial_tv_h)+q^i(\partial_t(u\nabla\phi),\nabla v_h)+q^i(u\nabla\phi,\partial_t\nabla v_h)\notag\\
 &=a_h(\partial_tR_hu,v_h)+a_h(R_hu,\partial_tv_h)+q^i(\partial_t(\Pi_{k-1}^0R_hu\Pi_{k}^0\nabla\phi),\Pi_{k-1}^0\nabla v_h)\notag\\
 &+q^i(\Pi_{k-1}^0R_hu\Pi_{k}^0\nabla\phi,\partial_t\Pi_{k-1}^0\nabla v_h).
 \end{align}
 Setting $v_h=\partial_tv_h$ in (\ref{defrh}), it follows that
 \begin{align} \label{3.75} a(u,\partial_tv_h)+q^i(u\nabla\phi,\nabla\partial_tv_h)=a_h(R_hu,\partial_tv_h)+q^i(\Pi_{k-1}^0R_hu\Pi_{k}^0\nabla\phi,
 \nabla\partial_t v_h).
 \end{align}
 Combining (\ref{3.74}) with (\ref{3.75}), we get
\begin{align}\label{3.76} a(\partial_tu,v_h)+q^i(\partial_t(u\nabla\phi),\nabla v_h)=a_h(\partial_tR_hu,v_h)+q^i(\partial_t(\Pi_{k-1}^0R_hu\Pi_{k}^0\nabla\phi),
 \Pi_{k-1}^0\nabla v_h).
 \end{align}
 To estimate $\|\nabla\partial_t(R_hu-\Pi_k^0u)\|_0$, for simplicity, we set $\psi=R_h u-\Pi_k^0 u$. Then from (\ref{3.76}), we have
 \begin{align}\label{gamma12}
 C_0\|\nabla\partial_t\psi\|_0^2&=C_0a(\partial_t\psi,\partial_t\psi)\leq a_h(\partial_t\psi,\partial_t\psi)\notag\\
 &=a_h(\partial_tR_hu,\partial_t\psi)-a_h(\partial_t\Pi_k^0 u,\partial_t\psi)\notag\\
 &=\{a(\partial_tu,\partial_t\psi)-a_h(\partial_t\Pi_k^0 u,\partial_t\psi)\}\notag\\
 &~~~~+\{q^i(\partial_t(u\nabla\phi),\nabla \partial_t\psi)-q^i(\partial_t(\Pi_{k-1}^0R_hu\Pi_{k}^0\nabla\phi),
 \Pi_{k-1}^0\nabla \partial_t\psi)\}~~(by~~(\ref{3.76}))\notag\\
 &:=\Gamma_1+\Gamma_2.
 \end{align}
 Using (\ref{ahcon}), it is easy to deduce that
\begin{align}\label{gamma1}
\Gamma_1&=a(\partial_tu,\partial_t\psi)-a_h(\partial_t\Pi_k^0 u,\partial_t\psi)=\sum_E a^E(\partial_tu-\partial_t\Pi_k^0 u,\partial_t\psi)\notag\\
&\leq \sum_E\|\nabla(\partial_tu-\partial_t\Pi_k^0 u)\|_{0,E}\|\nabla\partial_t\psi\|_{0,E}\leq Ch^k\|\partial_tu\|_{k+1}\|\nabla\partial_t\psi\|_{0}\notag\\
&\leq Ch^{2k}\|\partial_tu\|_{k+1}^2+\epsilon\|\nabla\partial_t\psi\|_{0}^2.
\end{align}
There holds
\begin{align}\label{gamma2}
\Gamma_2&=q^i(\partial_t(u\nabla\phi),\nabla \partial_t\psi)-q^i(\partial_t(\Pi_{k-1}^0R_hu\Pi_{k}^0\nabla\phi),
 \Pi_{k-1}^0\nabla \partial_t\psi)\notag\\
 &=\{q^i(\partial_t(u\nabla\phi-R_hu\nabla\phi),\nabla \partial_t\psi)\}\notag\\
&~~~~+\{q^i(\partial_t(R_hu\nabla\phi),\nabla \partial_t\psi)-q^i(\partial_t(\Pi_{k-1}^0R_hu\Pi_{k}^0\nabla\phi),
 \Pi_{k-1}^0\nabla \partial_t\psi)\}\notag\\
&\leq C(\|\partial_t(u-R_hu)\|_0+\|u-R_hu\|_0)\|\nabla \partial_t\psi\|_0.\notag\\
&~~~~+\{q^i(\partial_t(R_hu\nabla\phi),\nabla \partial_t\psi)-q^i(\partial_t(\Pi_{k-1}^0R_hu\Pi_{k}^0\nabla\phi),
 \Pi_{k-1}^0\nabla \partial_t\psi)\}.
\end{align}
Next, we estimate the last term.
\begin{align}\label{3.82}
&q^i(\partial_t(R_hu\nabla\phi),\nabla \partial_t\psi)-q^i(\partial_t(\Pi_{k-1}^0R_hu\Pi_{k}^0\nabla\phi),
 \Pi_{k-1}^0\nabla \partial_t\psi)\notag\\
&~~~~=q^i(\partial_t(R_hu)\nabla\phi,\nabla \partial_t\psi)+q^i(R_hu\partial_t\nabla\phi),\nabla \partial_t\psi)\notag\\
&~~~~~~-q^i(\partial_t(\Pi_{k-1}^0R_hu)\Pi_{k}^0\nabla\phi,
 \Pi_{k-1}^0\nabla \partial_t\psi)-q^i(\Pi_{k-1}^0R_hu\partial_t(\Pi_{k}^0\nabla\phi), \Pi_{k-1}^0\nabla \partial_t\psi)\notag\\
&=\sum_E \{b_i^E(\partial_t(R_hu),\phi,\partial_t\psi)-b_{i,h}^E(\partial_t(R_hu),\phi,\partial_t\psi)\notag\\
&~~~~+b_i^E(R_hu,\partial_t\phi,\partial_t\psi)-b_{i,h}^E(R_hu,\partial_t\phi,\partial_t\psi)\}.
\end{align}
Similarly as the deduction of (\ref{3.49}), by taking $u=\partial_tR_hu,~w=\phi,~v_h=\partial_t\psi$ and $u=R_hu,~w=\partial_t\phi,~v_h=\partial_t\psi$  in (\ref{bi-bih3}), respectively, we have
\begin{align*} &b_i^E(\partial_t(R_hu),\phi,\partial_t\psi)-b_{i,h}^E(\partial_t(R_hu),\phi,\partial_t\psi)\leq C(\|\partial_tu-\partial_tR_hu\|_{0,E}+h^k_E\|\phi\|_{k+1,E})\|\nabla\partial_t\psi\|_{0,E},
\end{align*}
and
\begin{align*} &b_i^E(R_hu,\partial_t\phi,\partial_t\psi)-b_{i,h}^E(R_hu,\partial_t\phi,\partial_t\psi)\leq C(\|u-R_hu\|_{0,E}+h^k_E\|\partial_t\phi\|_{k+1,E})\|\nabla\partial_t\psi\|_{0,E}.
\end{align*}
Inserting the above two inequality into (\ref{3.82}) and using (\ref{rhl2}), we get
\begin{align}\label{3.83}&q^i(\partial_t(R_hu\nabla\phi),\nabla \partial_t\psi)-q^i(\partial_t(\Pi_{k-1}^0R_hu\Pi_{k}^0\nabla\phi),
 \Pi_{k-1}^0\nabla \partial_t\psi)\notag\\
&\leq C \sum_E\{(\|\partial_tu-\partial_tR_hu\|_{0,E}+\|u-R_hu\|_{0,E}+h^k)\|\nabla\partial_t\psi\|_{0,E}\}\notag\\
&\leq C(h^{2k}+\|\partial_t(u-R_hu)\|_{0}^2)+\epsilon\|\nabla\partial_t\psi\|_{0}^2.
\end{align}
Combining (\ref{gamma2}) with (\ref{3.83}), and using (\ref{rhl2}) again, it yields
\begin{align}\label{gamma2-2}
\Gamma_2\leq C(h^{2k}+\|\partial_t(u-R_hu)\|_{0}^2)+\epsilon\|\nabla\partial_t\psi\|_{0}^2.
\end{align}
At last, substituting (\ref{gamma1}) and (\ref{gamma2-2}) into (\ref{gamma12}),  we have
\begin{align*}
\|\nabla\partial_t(R_hu-\Pi_k^0u)\|_0=\|\nabla\partial_t\psi\|_0\leq C(h^k+\|\partial_t(u-R_hu)\|_{0}).
\end{align*}
Hence,
\begin{align*}
\|\nabla\partial_t(R_hu-u)\|_0\leq C(h^k+\|\partial_t(u-R_hu)\|_{0}).
\end{align*}
This completes the proof of the lemma. $\hfill\Box$
\begin{lemma}\label{trhl2} Suppose $u,~\partial_t u\in H^{k+1}(\ome)$, $\phi\in W^{2,\infty}(\ome)\cap H^{k+2}(\ome)$ and $\partial_t\phi(t)\in W^{1,\infty}(\ome)$. Then we have
\begin{align*}
\|\partial_t(R_hu-u)\|_0\leq Ch^k.
\end{align*}
\end{lemma}
\textbf{Proof.}~  Similar as the proof of (\ref{rhl2}), we fist define an adjoint problem as follows:
\begin{small}
\begin{align}\label{ad-prob2}
\left\{\begin{array}{lr}
-\Delta w^i+q^i\nabla\phi\cdot\nabla w^i=\partial_t(u-R_hu),~~x\in \ome,\vspace{2mm}\\
w^i=0,~~x\in \partial\ome.\end{array}\right.
\end{align}
\end{small}
If $\phi\in W^{1,\infty}$,  then we have (cf. \cite{sun2016error})
\begin{align}\label{regu2}
\|w^i\|_2\leq C\|\partial_t(u-R_hu)\|_0.
\end{align}
Denote by $\xi=u-R_hu$. Then from (\ref{ad-prob2}), we get
\begin{align}\label{f1-5}
\|\partial_t\xi\|_0^2&=(\partial_t\xi,\partial_t\xi)=(\nabla\partial_t\xi,\nabla w^i)+(q^i\nabla\phi\nabla w^i,\partial_t\xi)\notag\\
&=a(\partial_t\xi,w^i)+b_i(\partial_t\xi,\phi,w^i)\notag\\
&=a(\partial_t\xi,w^i-w^i_I)+a(\partial_t\xi,w_I^i)+b_i(\partial_t\xi,\phi,w^i-w^i_I)+b_i(\partial_t\xi,\phi,w_I^i)\notag\\
&~~~~+a_h(\partial_tR_h u,w^i_I)-a_h(\partial_tR_hu,w_I^i)+b_{i,h}(\partial_tR_hu,\phi,w_I^i)-b_{i,h}(\partial_tR_hu,\phi,w_I^i)\notag\\
&=a(\partial_t\xi,w^i-w^i_I)+\{a_h(\partial_tR_h u,w^i_I)-a(\partial_tR_h u,w^i_I)\}\notag\\
&~~~~+b_i(\partial_t\xi,\phi,w^i-w^i_I)+\{b_{i,h}(\partial_tR_hu,\phi,w_I^i)-b_{i}(\partial_tR_hu,\phi,w_I^i)\}\notag\\
&~~~~+\{a(\partial_t u,w_I^i)-a_h(\partial_tR_hu,w_I^i)+b_{i}(\partial_t u,\phi,w_I^i)-b_{i,h}(\partial_tR_hu,\phi,w_I^i)\}\notag\\
&:=F_1+F_2+F_3+F_4+F_5.
\end{align}
Similar as the estimate for $D_1-D_4$ in (\ref{d1})-(\ref{d3}) and (\ref{d4-1}), respectively, by using (\ref{regu2}) and Lemma \ref{rhh1}, we have
\begin{align}\label{f1}
F_1&=a(\partial_t\xi,w^i-w_I^i)\leq C\|\nabla\partial_t\xi\|_0h\|\partial_t\xi\|_0~~(by~~(\ref{regu2})\notag\\
&\leq C(h^k+\|\partial_t\xi\|_0)h\|\partial_t\xi\|_0~~(by~~Lemma~\ref{rhh1})\notag\\
&\leq C(h^{2k+2}+h^2\|\partial_t\xi\|_0^2)+\epsilon\|\partial_t\xi\|_0^2,
\end{align}
\begin{align}
F_2&=a_h(\partial_tR_h u,w^i_I)-a(\partial_tR_h u,w^i_I)\leq C(h^{k+1}+h\|\nabla\partial_t\xi\|_0)\|\partial_t\xi\|_0\notag\\
&\leq C(h^{2k+2}+h^2\|\partial_t\xi\|_0^2)+\epsilon\|\partial_t\xi\|_0^2,
\end{align}
\begin{align}
F_3&=b_i(\partial_t\xi,\phi,w^i-w^i_I)\leq C\|\nabla\partial_t\xi\|_0h\|\partial_t\xi\|_0\notag\\
&\leq C(h^{2k+2}+h^2\|\partial_t\xi\|_0^2)+\epsilon\|\partial_t\xi\|_0^2~~(by~~Lemma~\ref{rhh1}),
\end{align}
and
\begin{align}\label{f4}
F_4&=b_{i,h}(\partial_tR_hu,\phi,w_I^i)-b_{i}(\partial_tR_hu,\phi,w_I^i)\leq C(h\|\partial_t\xi\|_0^2+h^k \|\partial_t\xi\|_0)\notag\\
&\leq C(h^{2k}+h\|\partial_t\xi\|_0^2)+\epsilon\|\partial_t\xi\|_0^2.
\end{align}
Next, we show the estimate for $F_5$. First, for any $v_h\in Q_h^k$, (\ref{3.76}) can be written as
\begin{align*}
&a(\partial_t u,v_h)-a_h(\partial_t R_h u,v_h)+b_i(\partial_t u,\phi,v_h)-b_{i,h}(\partial_t(R_hu),\phi,v_h)\notag\\
&~~~~=b_{i,h}(R_hu,\partial_t\phi,v_h)-b_i(u,\partial_t\phi,v_h).
\end{align*}
Taking $v_h=w^i_I$, then
\begin{align}\label{tf5}
F_5&=a(\partial_t u,w^i_I)-a_h(\partial_t R_h u,w^i_I)+b_i(\partial_t u,\phi,w^i_I)-b_{i,h}(\partial_t(R_hu),\phi,w^i_I)\notag\\
&=b_{i,h}(R_hu,\partial_t\phi,w^i_I)-b_i(u,\partial_t\phi,w^i_I)\notag\\
&=b_{i,h}(R_hu,\partial_t\phi,w^i_I-w^i)+b_{i,h}(R_hu,\partial_t\phi,w^i)-b_{i}(u,\partial_t\phi,w^i_I-w^i)-b_{i}(u,\partial_t\phi,w^i)\notag\\
&=\{b_{i,h}(R_hu,\partial_t\phi,w^i_I-w^i)-b_{i}(R_hu,\partial_t\phi,w^i_I-w^i)\}\notag\\
&+\{b_{i,h}(R_hu,\partial_t\phi,w^i)-b_{i}(R_hu,\partial_t\phi,w^i)\}\notag\\
&+\{b_{i}(R_hu,\partial_t\phi,w^i_I-w^i)-b_{i}(u,\partial_t\phi,w_I^i-w^i)
+\{b_i(R_hu,\partial_t\phi,w^i)-b_i(u,\partial_t\phi,w^i)\}\notag\\
&:=\tilde{F_1}+\tilde{F_2}+\tilde{F_3}+\tilde{F_4}.
\end{align}
Similarly as the deduction of (\ref{3.49}), by taking $u=R_hu,~w=\partial_t\phi,~v_h=p_I^i-p^i$ in (\ref{bi-bih3}), respectively, we have
\begin{align*} \tilde{F_1}=&\sum_E\{b_{i,h}^E(R_hu,\partial_t\phi,w_I^i-w^i)-b_{i}^E(R_hu,\partial_t\phi,w_I^i-w^i)\}\notag\\
&\leq C\sum_E(\|\xi\|_{0,E}+h^k_E)\|\nabla(w_I^i-w^i)\|_{0,E}\notag\\
&\leq Ch^{2k}+\epsilon\|\partial_t\xi\|_0^2~~(by~~(\ref{rhl2})~~and~~(\ref{regu2})).
\end{align*}
Similarly, 
from (\ref{rhl2}), we deduce
\begin{align*} \tilde{F_2}=&\sum_E\{b_{i,h}^E(R_hu,\partial_t\phi,w^i)-b_{i}^E(R_hu,\partial_t\phi,w^i)\}\notag\\
&\leq C\sum_E(\|\xi\|_{0,E}+h^k_E)\|\nabla w^i\|_{0,E}\notag\\
&\leq Ch^{2k}+\epsilon\|\partial_t\xi\|_0^2
\end{align*}
and
\begin{align*}\tilde{F_3}+\tilde{F_4}&=\{b_{i}(R_hu-u,\partial_t\phi,w^i_I-w^i)\}+\{b_{i}(R_hu-u,\partial_t\phi,w^i)\}\notag\\
&\leq C\|\xi\|_0\|(\nabla(w^i_I-w^i)\|_0+\|\nabla w^i\|_0)\notag\\
&\leq Ch^{2k}+\epsilon\|\partial_t\xi\|_0^2.
\end{align*}
Inserting the estimates for $\tilde{F_1}-\tilde{F_4}$ into (\ref{tf5}), it yields
\begin{align} \label{f5-2}F_5\leq C h^{2k}+\epsilon\|\partial_t\xi\|_0^2.\end{align}
Combining (\ref{f1-5}), (\ref{f1})-(\ref{f4}) and (\ref{f5-2}), we get
\begin{align*} \|\partial_t\xi\|_0^2\leq C(h^{2k}+h\|\partial_t\xi\|_0^2)+\epsilon\|\partial_t\xi\|_0^2.
\end{align*}
Hence, if $h$ is small enough, we have
\begin{align*}
\|\partial_t(u-R_hu)\|_0=\|\partial_t\xi\|_0\leq Ch^k,
\end{align*}
which finishes the proof of this lemma. $\hfill\Box$

\vspace{3mm}
\subsection{Error estimates for semidiscrete case in the $H^1$ norm}

In this subsection, we give a priori error estimate in the $H^1$ norm for the semi-discrete solution.  Assume
\begin{align}\label{reassu2} &p^i,~~p^i_t\in L^{\infty}(0,T;H^{k+1}(\ome)\cap L^{\infty}(\ome)),~~i=1,2,~~\notag\\
&\phi\in L^{\infty}(0,T;H^{k+2}(\ome)\cap W^{2,\infty}(\ome)),~~and~~\partial_t\phi(t)\in L^{\infty}(0,T;W^{1,\infty}(\ome)).
\end{align}
We also suppose
\begin{align}\label{assuf2} f\in L^{\infty}(0,T;H^{k}(\ome))~~and~~F^i\in L^{\infty}(0,T;H^{k+1}(\ome)).
\end{align}
In the following lemma, we present the error estimate of $\nabla\phi-\nabla\phi_h$ in the $H^1$ norm.
\begin{lemma}\label{tphi} Let $(\phi,p^{i})$ and $(\phi_h,p_h^{i})$ be the solutions of (\ref{weak formulation}) and (\ref{semi form}), respectively. Assume (\ref{reassu2}) and (\ref{assuf2}) hold. Then we have
\begin{align*}
\|\partial_t(\nabla\phi-\nabla\phi_h)\|_0\leq C(h^k+\sum_{i=1}^2\|\partial_t(R_hp^i-p_h^i)\|_0).
\end{align*}
\end{lemma}
\textbf{Proof.}~For any $t\in (0,T)$, suppose $\tilde{R}_h\phi \in Q_h^k$ is the $H^1$ projection to $\phi(t)$, satisfying
\begin{align}\label{deftrh}a_h(\tilde{R}_h\phi,w_h)=a(\phi,w_h),~~\forall w_h\in Q_h^k.\end{align}

From (\ref{weak formulation}) and (\ref{semi form}), we have
\begin{align*} a(\phi,w_h)-a_h(\phi_h,w_h)+\tilde{b}(p^1,p^2,w_h)-\tilde{b}_h(p^1_h,p^2_h,w_h)=(f-f_h,w_h).
\end{align*}
Then, from (\ref{deftrh}), it follows that
\begin{align} \label{3.103}a_h(\tilde{R}_h\phi-\phi_h,w_h)=(f-f_h,w_h)+\tilde{b}_h(p^1_h,p^2_h,w_h)-\tilde{b}(p^1,p^2,w_h).
\end{align}
Integrating with respect to $t$ on both side of the above and setting $\eta=\tilde{R}_h\phi-\phi_h$, then it yields
\begin{align}\label{3.105} a_h(\partial_t\eta,w_h)+a_h(\eta,\partial_t w_h)&=(\partial_t(f-f_h),w_h)+(f-f_h,\partial_tw_h)\notag\\
&~~~~+\tilde{b}_h(\partial_tp^1_h,\partial_tp^2_h,w_h)+\tilde{b}_h(p^1_h,p^2_h,\partial_tw_h)\notag\\
&~~~~-\tilde{b}(\partial_tp^1,\partial_tp^2,w_h)-\tilde{b}(p^1,p^2,\partial_tw_h).
\end{align}
Setting $w_h=\partial_tw_h$ in (\ref{3.103}), then we get
\begin{align*}a_h(\eta,\partial_tw_h)=(f-f_h,\partial_tw_h)+\tilde{b}_h(p^1_h,p^2_h,\partial_tw_h)-\tilde{b}(p^1,p^2,\partial_tw_h).
\end{align*}
Inserting the above into (\ref{3.105}), it follows
\begin{align*}
a_h(\partial_t\eta,w_h)=(\partial_t(f-f_h),w_h)+\tilde{b}_h(\partial_tp^1_h,\partial_tp^2_h,w_h)-\tilde{b}(\partial_tp^1,\partial_tp^2,w_h).
\end{align*}
Taking $w_h=\partial_t\eta$ and using Lemma \ref{lemma4.2}, we have
\begin{align}a_h(\partial_t\eta,\partial_t\eta)&=(\partial_t(f-f_h),\partial_t\eta)+\tilde{b}_h(\partial_tp^1_h,\partial_tp^2_h,\partial_t\eta)
-\tilde{b}(\partial_tp^1,\partial_tp^2,\partial_t\eta)\notag\\
&\leq \|\partial_t(f-f_h)\|_0\|\partial_t\eta\|_0+|\tilde{b}_h(\partial_tp^1_h,\partial_tp^2_h,\partial_t\eta)
-\tilde{b}(\partial_tp^1,\partial_tp^2,\partial_t\eta)|\notag\\
&\leq C(h^k+\sum_{i=1}^2\|\partial_t(p^i-p_h^i\|_0)\|\nabla\partial_t\eta\|_0
\end{align}
Hence,
\begin{align*} C_0\|\partial_t\nabla \eta\|_0^2&=C_0\|\nabla\partial_t\eta\|_0^2\leq a_h(\partial_t\eta,\partial_t\eta)\notag\\
&\leq C(h^k+\sum_{i=1}^2\|\partial_t(p^i-p_h^i\|_0)\|\nabla\partial_t\eta\|_0,
\end{align*}
which deduce  that
\begin{align}\label{trhphi}\|\partial_t(\nabla\tilde{R}_h\phi-\nabla\phi_h)\|_0\leq C(h^k+\sum_{i=1}^2\|\partial_t(p^i-p_h^i)\|_0).
\end{align}
Next, we estimate $\|\partial_t(\nabla\tilde{R}_h\phi-\nabla\phi)\|_0$. First, integrating with respect to $t$ on both sides of
(\ref{deftrh}), we have
\begin{align*}
a_h(\partial_t\tilde{R}_h\phi,w_h)+a_h(\tilde{R}_h\phi,\partial_tw_h)=a(\partial_t\phi,w_h)+a(\phi,\partial_tw_h).
\end{align*}
Taking $w_h=\partial_tw_h$ in (\ref{deftrh}) and inserting the resulted equation into the above, we get
\begin{align*}a_h(\partial_t\tilde{R}_h\phi,w_h)=a(\partial_t\phi,w_h).
\end{align*}
Then setting $w_h=\partial_t\tilde{\eta}=\partial_t(\nabla\tilde{R}_h\phi-\Pi_k^0\nabla\phi)$, it follows that
\begin{align*}C_0\|\nabla\partial_t\tilde{\eta}\|_0^2&\leq C_0a_h(\partial_t\tilde{\eta},\partial_t\tilde{\eta})= a_h(\partial_t\tilde{R}_h\phi,\partial_t\tilde{\eta})-a_h(\partial_t\Pi_k^0\phi,\partial_t\tilde{\eta})\notag\\
&= \sum_E\big{(}a^E(\partial_t\phi,\partial_t\tilde{\eta})-a^E(\partial_t\Pi_k^0\phi,\partial_t\tilde{\eta})\big{)}\leq C\sum_E\|\nabla\partial_t\phi-\nabla\Pi_k^0\nabla\partial_t\phi\|_0\|\nabla\partial_t\tilde{\eta}\|_0\notag\\
&\leq Ch^k\|\nabla\partial_t\tilde{\eta}\|_0.
\end{align*}
Hence,
\begin{align*}
\|\partial_t(\nabla\tilde{R}_h\phi-\Pi_k^0\nabla\phi)\|_0\leq Ch^k.
\end{align*}
Then
\begin{align*}\|\partial_t(\nabla\phi-\nabla\tilde{R}_h\phi)\|_0&\leq \|\partial_t(\nabla\phi-\Pi_k^0\nabla\phi)\|_0
+\|\partial_t(\Pi_k^0\nabla\phi-\nabla\tilde{R}_h\phi)\|_0\notag\\
&\leq C h^k.
\end{align*}
Combing the above with (\ref{trhphi}) and using Lemma \ref{trhl2}, we get
\begin{align*}
\|\partial_t(\nabla\phi-\nabla\phi_h)\|_0\leq C(h^k+\sum_{i=1}^2\|\partial_t({R}_hp^i-p_h^i)\|_0).
\end{align*}
This completes the proof of this lemma.$\hfill\Box$

\begin{theorem}\label{theorem_p_H1}
Let $(\phi,p^{i})$ and $(\phi_h,p_h^{i})$ be the solutions of (\ref{weak formulation}) and (\ref{semi form}), respectively. Suppose (\ref{reassu2}) and (\ref{assuf}) hold. For all $t\in(0,T]$, the following estimate holds
\begin{align}
\|p^i(t)-p_h^i(t)\|_{1}+\|\phi(t)-\phi_h(t)\|_1\leq C h^k.
\end{align}
\end{theorem}
\textbf{Proof.}~  From (\ref{weak formulation}) and (\ref{semi form}), for any $v_h\in Q_h^k$,we have
\begin{align}\label{eqn3321}
(p_t^i,v_h)-m_h(p_{h,t}^i&,v_h)+a(p^i,v_h)-a_h(p_h^i,v_h)\notag\\
&+b_i(p^i,\phi_h,v_h)-b_{i,h}(p_h^i,\phi_h,v_h)=(F^i,v_h)-(F^i_h,v_h).
\end{align}
Then from (\ref{defrh}), it follows that
\begin{align*}
(p_t^i-(R_hp^i)_t&+(R_hp^i)_t,v_h)+m_h\big{(}(R_hp^i)_t-p_{h,t}^i,v_h\big{)}-m_h\big{(}(R_hp^i)_t,v_h\big{)}\\
&+a_h(R_hp^i-p_h^i,v_h)+b_{i,h}(R_hp^i,\phi,v_h)-b_{i,h}(p_h^i,\phi_h,v_h)=(F^i-F_h^i,v_h).\notag\\
\end{align*}
Let $\theta^i=R_hp^i-p_h^i,~\eta^i=p^i-R_hp^i$ and take $v_h=\theta_t^i:=(R_hp^i-p_h^i)_t$ in the above equation. Then it follows that
\begin{align}\label{a1234}
&m_h(\theta_t^i,\theta_t^i)+a_h(\theta^i,\theta_t^i)\notag\\
&=-(\eta_t^i,\theta_t^i)+\big{\{}m_h\big{(}(R_hp^i)_t,\theta_t^i)-\big{(}(R_hp^i)_t,\theta_t^i\big{)}\big{\}}\notag\\
& +(F^i-F_h^i,\theta_t^i)+\big{\{}b_{i,h}(p_h^i,\phi_h,\theta_t^i)-b_{i,h}(R_hp^i,\phi,\theta_t^i)\big{\}}\notag\\
&:=A_1+A_2+A_3+A_4.
\end{align}
Next, we shall estimate $A_i,~i=1,2,3,4$ respectively.
First, from Lemma \ref{trhl2}, we get
\begin{align}
A_1&=-(\eta_t^i,\theta_t^i)\leq \|p_t^i-(R_hp^i)_t\|_0\|\theta_t^i\|_0\notag\\
&\leq Ch^{k}\|\theta_t^i\|_0\leq Ch^{2k}+\epsilon\|\theta_t^i\|_0^2.
\end{align}

Using (\ref{cons}), (\ref{continuity}), Lemmas \ref{rhl2err} and \ref{trhl2}, it yields
\begin{align}
A_2&=m_h\big{(}(R_hp^i)_t,\theta_t^i)-\big{(}(R_hp^i)_t,\theta_t^i\big{)}\notag\\
&=m_h\big{(}(R_hp^i)_t-\Pi_k^0(R_hp^i)_t,\theta_t^i\big{)}+m_h\big{(}(\Pi_k^0(R_hp^i)_t,\theta_t^i\big{)}-\big{(}(R_hp^i)_t,\theta_t^i\big{)}\notag\\
&\leq C\|(R_hp^i)_t-\Pi_k^0(R_hp^i)_t\|_0\|\theta_t^i\|_0~~\big(by~~(\ref{cons})~~and~~(\ref{continuity})\big)\notag\\
&\leq C h^k\|\theta_t^i\|_0~~\big(by~~Lemmas~~ \ref{rhl2err} ~~and~~ \ref{trhl2}\big)\notag\\
&\leq Ch^{2k}+\epsilon \|\theta_t^i\|_0^2.
\end{align}

For the term $A_3$, we have
\begin{align}\label{a3} A_3&=(F^i-F_h^i,\theta_t^i)\leq Ch^{k+1}\|F^i\|_{k+1}\|\theta_t^i\|_0\notag\\
&\leq Ch^{2k+2}+\epsilon\|\theta_t^i\|_0^2.
\end{align}

To estimate the term $A_4$, first we get
\begin{align}\label{a4}A_4&=b_{i,h}(p_h^i,\phi_h,\theta_t^i)-b_{i,h}(R_hp^i,\phi,\theta_t^i)\notag\\
&=b_{i,h}(p_h^i-R_hp^i,\phi,\theta_t^i)+b_{i,h}(p_h^i-p^i,\phi_h-\phi,\theta_t^i)+b_{i,h}(p^i,\phi_h-\phi,\theta_t^i)\notag\\
&=A_{41}+A_{42}+A_{43}.
\end{align}
Next, we shall estimate the terms $A_{41},~A_{42},~A_{41}$, respectively.
\par From Lemma \ref{trhl2} and Theorem \ref{theorem3.1}, it yields
\begin{align}\label{a41}&A_{41}=b_{i,h}(p_h^i-R_hp^i,\phi,\theta_t^i)=\sum_Eq^i(\Pi_{k-1}^0\theta^i\Pi_{k}^0\nabla\phi,\Pi_{k-1}^0
\nabla\theta_t^i)_E\notag\\
&=\sum_Eq^i\frac {\partial}{\partial_t}(\Pi_{k-1}^0\theta^i\Pi_{k}^0\nabla\phi,\Pi_{k-1}^0\nabla\theta^i)_E-\sum_E\Big\{q^i\big((\ \Pi_{k-1}^0\theta^i)_t\Pi_{k}^0\nabla\phi,\Pi_{k-1}^0\nabla\theta^i\big)_E\notag\\
& +q^i\big{(}\Pi_{k-1}^0\theta^i(\Pi_{k}^0\nabla\phi)_t,\Pi_{k-1}^0\nabla\theta^i\big{)}_E\Big\}\notag\\
&\leq \frac{\partial}{\partial_t}b_{i,h}(\theta^i,\phi,\theta^i)+\sum_E\Big{\{}\|(\Pi_{k-1}^0\theta^i)_t\|_{0,E}\|\Pi_{k}^0\nabla\phi\|_{0,\infty,E}\|\Pi_{k-1}^0\nabla\theta^i\|_{0,E}\notag\\
& +\|\Pi_{k-1}^0\theta^i\|_{0,E}\|(\Pi_{k}^0\nabla\phi)_t\|_{0,\infty,E}\|\Pi_{k-1}^0\nabla\theta^i\|_{0,E}\Big{\}}\notag\\
&\leq \frac{\partial}{\partial_t}b_{i,h}(\theta^i,\phi,\theta^i)+C(\|\theta^i\|_0^2+\|\nabla\theta^i\|_0^2)+\epsilon\|\theta_t^i\|_0^2\notag\\
&\leq \frac{\partial}{\partial_t}b_{i,h}(\theta^i,\phi,\theta^i)+C(h^{2k}+\|\nabla\theta^i\|_0^2)+\epsilon\|\theta_t^i\|_0^2~~(by~~Lemma~~ \ref{rhl2err}~~ and~~ Theorem ~~\ref{theorem3.1}).
\end{align}

Similarly, we get
\begin{align} A_{42}&=b_{i,h}(p_h^i-p^i,\phi_h-\phi,\theta_t^i)=\sum_Eq^i(\Pi_{k-1}^0(p^i_h-p^i)\Pi_{k}^0\nabla(\phi_h-\phi),\Pi_{k-1}^0
\nabla\theta_t^i)_E\notag\\
&=\sum_Eq^i\frac {\partial}{\partial_t}(\Pi_{k-1}^0(p^i_h-p^i)\Pi_{k}^0\nabla(\phi_h-\phi),\Pi_{k-1}^0\nabla\theta^i)_E\notag\\
& -\sum_E\big{\{}q^i\big{(}(\ \Pi_{k-1}^0(p^i_h-p^i))_t\Pi_{k}^0\nabla(\phi_h-\phi),\Pi_{k-1}^0\nabla\theta^i\big{)}_E
+q^i\big{(}\Pi_{k-1}^0(p^i_h-p^i)(\Pi_{k}^0\nabla(\phi_h-\phi))_t,\Pi_{k-1}^0\nabla\theta^i\big{)}_E\big{\}}\notag\\
&\leq \frac{\partial}{\partial_t}b_{i,h}(p^i_h-p^i,\phi_h-\phi,\theta^i)+\sum_E\big{\{}\|(\Pi_{k-1}^0(p^i-p_h^i))_t\|_{0,E}\|\Pi_{k}^0\nabla(\phi_h-\phi)\|_{0,\infty,E}\|\Pi_{k-1}^0\nabla\theta^i\|_{0,E}\notag\\
& +\|\Pi_{k-1}^0(p^i-p_h^i)\|_{0,E}\|(\Pi_{k}^0\nabla(\phi_h-\phi))_t\|_{0,\infty,E}\|\Pi_{k-1}^0\nabla\theta^i\|_{0,E}\big{\}}\notag\\
&\leq \frac{\partial}{\partial_t}b_{i,h}(p^i_h-p^i,\phi_h-\phi,\theta^i)+C(\|p^i-p_h^i\|_0^2+\|\nabla\theta^i\|_0^2)+\epsilon\|(p ^i-p_h^i)_t\|_0^2\notag\\
&\leq \frac{\partial}{\partial_t}b_{i,h}(p^i_h-p^i,\phi_h-\phi,\theta^i)+C(h^{2k}+\|\nabla\theta^i\|_0^2)+\epsilon\|\theta_t^i\|_0^2,
~~(by~~Lemma~~ \ref{rhl2err}~~ and~~ Theorem \ref{theorem3.1})
\end{align}
and
\begin{align}\label{a43}
A_{43}&=b_{i,h}(p^i,\phi_h-\phi,\theta_t^i)=\sum_Eq^i(\Pi_{k-1}^0p^i\Pi_{k}^0\nabla(\phi_h-\phi),\Pi_{k-1}^0
\nabla\theta_t^i)_E\notag\\
&=\sum_Eq^i\frac {\partial}{\partial_t}(\Pi_{k-1}^0p^i\Pi_{k}^0\nabla(\phi_h-\phi),\Pi_{k-1}^0\nabla\theta^i)_E\notag\\
& -\sum_E\big{\{}q^i\big{(}(\ \Pi_{k-1}^0p^i)_t\Pi_{k}^0\nabla(\phi_h-\phi),\Pi_{k-1}^0\nabla\theta^i\big{)}_E
+q^i\big{(}\Pi_{k-1}^0p^i(\Pi_{k}^0\nabla(\phi_h-\phi))_t,\Pi_{k-1}^0\nabla\theta^i\big{)}_E\big{\}}\notag\\
&\leq \frac{\partial}{\partial_t}b_{i,h}(p^i,\phi_h-\phi,\theta^i)+
\sum_E\big{\{}\|(\Pi_{k-1}^0p^i)_t\|_{0,\infty,E}\|\Pi_{k}^0\nabla(\phi_h-\phi)\|_{0,E}\|\Pi_{k-1}^0\nabla\theta^i\|_{0,E}\notag\\
& +\|\Pi_{k-1}^0p^i\|_{0,\infty,E}\|(\Pi_{k}^0\nabla(\phi_h-\phi))_t\|_{0,E}\|\Pi_{k-1}^0\nabla\theta^i\|_{0,E}\big{\}}\notag\\
&\leq \frac{\partial}{\partial_t}b_{i,h}(p^i,\phi_h-\phi,\theta^i)+
C(\|\nabla(\phi_h-\phi)\|_0+\|\nabla(\phi_h-\phi)_t\|_0)\|\nabla\theta^i\|_0\notag\\
&\leq \frac{\partial}{\partial_t}b_{i,h}(p^i,\phi_h-\phi,\theta^i)+C(h^{2k}+
\|\nabla\theta^i\|_0^2)+\epsilon\sum_{i=1}^2\|\theta_t^i\|_0^2,
\end{align}
where Lemma~~ \ref{trhl2},~~Theorem \ref{theorem3.1}
~~and~~Lemma \ref{tphi} are used in the last inequality. Inserting (\ref{a41})-(\ref{a43}) into (\ref{a4}), it follows that
\begin{align}\label{a4-2}A_4&\leq \frac{\partial}{\partial_t}\big{(}b_{i,h}(\theta^i,\phi,\theta^i)+b_{i,h}(p_h^i-p^i,\phi_h-\phi,\theta^i)
+b_{i,h}(p^i,\phi_h-\phi,\theta^i)\big{)}\notag\\
& +C(\|\nabla\theta^i\|_0^2+h^{2k})+\sum_{i=1}^2\epsilon\|\theta_t^i\|_0^2.
\end{align}
Thus combing (\ref{a1234})-(\ref{a3}) and (\ref{a4-2}), we have
\begin{align*}\tilde{C}_*\|\theta_t^i\|_0^2+\frac 1 2\|\nabla\theta_t^i\|_0^2&\leq m_h(\theta_t^i,\theta_t^i)+a_h(\theta^i,\theta_t^i)\notag\\
&\leq \frac{\partial}{\partial_t}\big{(}b_{i,h}(\theta^i,\phi,\theta^i)+b_{i,h}(p_h^i-p^i,\phi_h-\phi,\theta^i)
+b_{i,h}(p^i,\phi_h-\phi,\theta^i)\big{)}\notag\\
& +C(\|\nabla\theta^i\|_0^2+h^{2k})+\epsilon\sum_{i=1}^2\|\theta_t^i\|_0^2.
\end{align*}
Then integrating with respect to $t$ on both sides, it yields
\begin{align*}
\|\nabla\theta\|_0^2& \leq C\big{(}b_{i,h}(\theta^i,\phi,\theta^i)+b_{i,h}(p_h^i-p^i,\phi_h-\phi,\theta^i)
+b_{i,h}(p^i,\phi_h-\phi,\theta^i)+\int_0^t\|\nabla\theta\|_0^2dt+h^{2k}\big{)}\notag\\
& \leq C\big{(}(\|\theta^i\|_0\|\nabla\phi\|_{0,\infty}+\|p_h^i-p^i\|_0\|\nabla\phi-\nabla\phi_h\|_{0,\infty}
+\|p^i\|_{0,\infty}\|\nabla\phi-\nabla\phi_h\|_{0})\|\nabla\theta^i\|_0\notag\\
&~~~~+\int_0^t\|\nabla\theta\|_0^2dt+h^{2k}\big{)}\notag\\
 &\leq C \big{(}\int_0^t\|\nabla\theta\|_0^2dt+h^{2k}+\epsilon\|\nabla\theta^i\|_0^2\big{)}.
 ~~(by~~Lemmas~~ \ref{theorem3.2}~~and~~\ref{rhl2err},~~and~~Theorem~~ \ref{theorem3.1})
 \end{align*}
 Thus, we have
 \begin{align*}
 \|\nabla\theta^i\|_0^2\leq C\big{(}\int_0^t\|\nabla\theta\|_0^2dt+h^{2k}\big{)}.
 \end{align*}
 By using the Gronwall inequlity, we get
 \begin{align*}
 \|\nabla\theta^i\|_0^2\leq C h^{2k}.
 \end{align*}
 Then using Lemma \ref{rhl2err}, we have
 \begin{align*}\|p^i-p_h^i\|_1\leq Ch^k,
 \end{align*}
 which combining with Lemma \ref{theorem3.2} show
 \begin{align*} \|\phi-\phi_h\|_1\leq Ch^k.
 \end{align*}
 This completes the proof of the lemma. $\hfill\Box$

\vspace{3mm}
\setcounter{lemma}{0}
\setcounter{theorem}{0}
\setcounter{corollary}{0}
\setcounter{equation}{0}
\section{Error estimates for fully discrete case}

we first construct the fully discrete virtual element formulation of (\ref{semi form}) 
and then derive the a priori error estimates for the scheme in the $L^2$ norm.
\par We employ the backward Euler scheme for the approximation of time derivative. Let $p_h^n:= p_h(\cdot, t_n)$, $n = 0, 1, . . . , N $ with $t_n=n\tau,~\tau=T/N.$ The fully discrete form corresponding to (\ref{semi form}) reads as: find $p_h^{i,n}$, $\phi_{h}^n\in Q_h^k$ such that
\begin{small}
\begin{equation}\label{full}
\left\{\begin{array}{lr}
m_{h}\left(\frac{p_h^{i,n}-p_h^{i,n-1}}{\tau}, v_{h}\right)+a_{h}\left(p_h^{i,n}, v_{h}\right)+b_{i,h}(p_h^{i,n},\phi_h^n, v_{h})=\left(F_{h}^{i,n}, v_{h}\right)~~\forall v_{h} \in Q_{h}^k,~\text { for } n=1, \ldots, N,\vspace{3mm}\\
{a}_{h}\left(\phi_h^{n}, w_{h}\right)+\tilde{b}_h(p_h^{1,n},p_h^{2,n},w_h)=\left(f_{h}^{n}, w_{h}\right)~~\forall w_{h} \in Q_{h}^k,~\text { for } n=1, \ldots, N,\vspace{3mm}\\
p_h^{i,0}=p^i_{h,0},~\phi_h^0=\phi_{h,0}.
\end{array}\right.
\end{equation}
\end{small}
\vspace{0mm}
Assume the following regularity properties hold,
\begin{align}\label{reassu3} &p^{i,n},~~p^{i,n}_t\in L^{\infty}(0,T;H^{k+1}(\ome)\cap L^{\infty}(\ome)),~~i=1,2,~~\notag\\
&\phi^n\in L^{\infty}(0,T;H^{k+2}(\ome)\cap W^{k+1,\infty}(\ome)),~~and~~\partial_t\phi^n\in L^{\infty}(0,T;W^{1,\infty}(\ome)).
\end{align}
We also suppose
\begin{align}\label{assuf3} f^n\in L^{\infty}(0,T;H^{k}(\ome))~~and~~F^{i,n}\in L^{\infty}(0,T;H^{k+1}(\ome)).
\end{align}
We present the error estimate of $\phi_h^n-\phi^n$ in the following lemma.
\begin{lemma}\label{theorem4.2}
Let $(\phi^n,p^{i,n})$ and $(\phi_h^n,p_h^{i,n})$ be the solutions of (\ref{weak formulation}) and (\ref{full}), respectively, and assume (\ref{reassu3})-(\ref{assuf3}) hold. 
Then for all $n=1,...,N$ there holds
\begin{equation}
\|\phi_h^n-\phi^n\|_{1}\leq C\Big(h^k+\sum\limits^{2}_{i=1}\|p_h^{i,n}-p^{i,n}\|_0\Big).\notag
\end{equation}
and
\begin{align}\label{piinfty2}\|\Pi_k^0\nabla\phi_h^n\|_{0,\infty,E}\leq Ch^{-1}\sum_{i=1}^2\|p_h^{i,n}-p^{i,n}\|_0.
\end{align}
\end{lemma}
\textbf{Proof.}~ The proof can be presented by repeating the arguments in Lemma \ref{theorem3.2}. $\hfill\Box$
\vspace{3mm}\\
\par Now, we proceed to estimate the error in the $L^2$ norm for $p^{i,n}$.\vspace{3mm}
\begin{theorem}\label{theorem4.1}
Let $(\phi^n,p^{i,n})$ and $(\phi_h^n,p_h^{i,n})$ be the solutions of (\ref{weak formulation}) and (\ref{full}), respectively. Suppose (\ref{reassu3})-(\ref{assuf3}) hold
 and $\tau$ is small enough. Then for all $n=1,...,N$, we have
\begin{eqnarray}
\sum\limits^{2}_{i=1}\|p_h^{i,n}-p^{i,n}\|_0+\|\phi_h^n-\phi^n\|_0\leq C(\tau+h^k).\notag
\end{eqnarray}
\end{theorem}
\textbf{Proof.}~ Set
\begin{eqnarray}
\begin{split}
p_h^{i,n}-p^{i,n}=(p_h^{i,n}-R_hp^{i,n})+(R_hp^{i,n}-p^{i,n})=:\upsilon^{i,n}+\varrho^{i,n},\notag
\end{split}
\end{eqnarray}
where $p^{i,n}=p^i(t_n),~~n=1,2,\cdots,N$. From Lemma \ref{rhl2err}, we obtain
\begin{eqnarray}
\|\varrho^{i,n}\|_0&\leq Ch^{k}.\notag
\end{eqnarray}
The estimate for $\upsilon^{i,n}$ requires more analysis. Denote by
\begin{align*} D_{\tau}u^n=\frac{u^n-u^{n-1}}{\tau}.
\end{align*}
For all $v_h\in Q_h^k$, from (\ref{weak formulation}) and (\ref{full}), it holds\vspace{-13mm}
\begin{spacing}{2.0}
\begin{align}\label{mh+ah}
&m_h(D_{\tau}\upsilon^{i,n},v_h)+a_h(\upsilon^{i,n},v_h)\notag\\
&=(F_h^{i,n},v_h)-b_{i,h}(p_h^{i,n},\phi_h^n,v_h)-m_h(D_{\tau}R_hp^{i,n},v_h)-a_h(R_hp^{i,n},v_h)\notag\\
&\;\;\;\;+a(p^{i,n},v_h)-a(p^{i,n},v_h)\notag\\
&=(F_h^{i,n},v_h)-b_{i,h}(p_h^{i,n},\phi_h^n,v_h)-m_h(D_{\tau}R_hp^{i,n},v_h)-a_h(R_hp^{i,n},v_h)\notag\\
&\;\;\;\;+a(p^{i,n},v_h)-(F^{i,n},v_h)+(p_t^{i,n},v_h)+b_i(p^{i,n},\phi^n,v_h)\notag\\
&=(F_h^{i,n}-F^{i,n},v_h)+\Big((p_t^{i,n},v_h)-m_h(D_{\tau}R_hp^{i,n},v_h)\Big)\notag\\
&\;\;\;\;+\big(a(p^{i,n},v_h)-a_h(R_hp^{i,n},v_h)\big)+\big(b_i(p^{i,n},\phi^n,v_h)-b_{i,h}(p_h^{i,n},\phi_h^n,v_h)\big)\notag\\
&:=H_{1}^n+H_{2}^n+H_{3}^n+H_{4}^n.
\end{align}
\end{spacing}
\vspace{-8mm}
\noindent From (\ref{sh1}), we have
\begin{eqnarray}\label{h1}
\begin{split}
H_{1}^n=(F_h^{i,n}-F^{i,n},v_h)=(\Pi_{k}^0F^{i,n}-F^{i,n},v_h)\leq Ch^{k+1}|F^{i,n}|_{k+1}\|v_h\|_0.
\end{split}
\end{eqnarray}
In order to bound the second term, adding and substracting suitable terms, we can obtain\vspace{-13mm}
\begin{spacing}{2.0}
\begin{align}\label{h2}
H_{2}^n&=(p_t^{i,n},v_h)-m_h\Bigg(D_{\tau}R_hp^{i,n},v_h\Bigg)\notag\\
&=\Big((p_t^{i,n},v_h)-(D_{\tau}p^{i,n},v_h)\Big)+\sum_E(D_{\tau}(p^{i,n}-\Pi_k^0p^{i,n}),v_h)_{E}+\Big(\sum_E(D_{\tau}\Pi_k^0p^{i,n},v_h)_E- m_h^E(D_{\tau}R_hp^{i,n},v_h)\Big)\notag\\
&:=H_{21}+H_{22}+H_{23}.
\end{align}
\end{spacing}
\vspace{-8mm}
The estimates for $H_{21}$ to $H_{23}$ can be determined as follows,
\begin{align*}H_{21}&\leq \|p_t^{i,n}-D_{\tau}p^{i,n}\|_0\|v_h\|_0\notag\\
&\leq \|\frac 1 2\tau\cdot \partial_{tt} p^i(x,\xi)\|_0\|v_h\|_0~~(t^{n-1}<\xi<t^n),\notag~~(by~~Taylor's~~ expansion)\\
&\leq C\tau \|v_h\|_0,
\end{align*}
\begin{align*} H_{22}&=\sum_E(D_{\tau}(p^{i,n}-\Pi_k^0p^{i,n}),v_h)\notag\\
&\leq \sum_E\frac 1 {\tau} \int_{t^{n-1}}^{t^n} \|\partial_t (p^{i,n}-\Pi_k^0p^{i,n})(s)\|_0 ds \|v_h\|_{0,E}\notag\\
&\leq C h^{k}\|v_h\|_0.
\end{align*}
\begin{align*} H_{23}&=\sum_E\big((D_{\tau}\Pi_k^0p^{i,n},v_h)_E- m_h^E(D_{\tau}R_hp^{i,n},v_h)\big)\notag\\
&=\sum_E\big(m_h^E(D_{\tau}\Pi_k^0p^{i,n}-D_{\tau}R_hp^{i,n},v_h)\big)~~(by~~(\ref{cons})) \notag\\
&\leq \sum_E\|D_{\tau}(\Pi_k^0p^{i,n}-p^{i,n})\|_{0,E}\|v_h\|_{0,E}+\|D_{\tau}(p^{i,n}-R_hp^{i,n})\|_0\|v_h\|_{}\notag\\
&\leq \sum_E\frac 1 {\tau} \int_{t^{n-1}}^{t^n} \|\partial_t (\Pi_k^0p^{i,n}-p^{i,n})(s)\|_0 ds \|v_h\|_{0,E}+
\int_{t^{n-1}}^{t^n} \|\partial_t (p^{i,n}-R_hp^{i,n})(s)\|_0 ds \|v_h\|_{0}\notag\\
&\leq C h^k\|v_h\|_0~~(by~~Lemma~~\ref{trhl2}).
\end{align*}
Hence
\begin{align}\label{h2}H_2^n=(p_t^{i,n},v_h)-m_h(D_{\tau}R_hp^{i,n},v_h)\leq C ({\tau}+h^k)\|v_h\|_0.
\end{align}
Using (\ref{continuity a}), the third term can be estimated a follows
\begin{eqnarray}\label{h3}
\begin{split}
H_{3}^n&=a(p^{i,n},v_h)-a_h(R_hp^{i,n},v_h)\notag\\
&=a(p^{i,n},v_h)-a_h(p^{i,n},v_h)+a_h(p^{i,n}-R_hp^{i,n},v_h)\notag\\
&\leq Ch^k|p^{i,n}|_{k+1}\|\nabla v_h\|_{0}.
\end{split}
\end{eqnarray}
From (\ref{piinfty2}) and Lemmas \ref{bi-bhi} and \ref{theorem4.2},  we can express the fourth term as follows
\begin{eqnarray}\label{h4}
\begin{split}
H_{4}^n&=b_i(p^{i,n},\phi^n,v_h)-b_{i,h}(p_h^{i,n},\phi_h^n,v_h)\\
&\leq C\Big(h^k+\sum\limits^{2}_{i=1}\|p_h^{i,n}-p^{i,n}\|_{0}+h^{-1}\sum\limits^{2}_{i=1}\|p_h^{i,n}-p^{i,n}\|_{0}^2\Big)
\|\nabla v_h\|_{0}.
\end{split}
\end{eqnarray}
 Collecting the estimation of (\ref{h1})-(\ref{h4}) into (\ref{mh+ah}) and taking $v_h=\upsilon^{i,n}$, we get
\begin{align}\label{mh}
m_h&\Big(\frac{\upsilon^{i,n}-\upsilon^{i,n-1}}{\tau},\upsilon^{i,n}\Big)+a_h(\upsilon^{i,n},\upsilon^{i,n})
=H_{1}^n+H_{2}^n+H_{3}^n+H_{4}^n\notag\\
&\leq C (\tau+h^k)\|\upsilon^{i,n}\|_0+C\Big(h^k+\sum\limits^{2}_{i=1}\|p_h^{i,n}-p^{i,n}\|_{0}
+h^{-1}\sum\limits^{2}_{i=1}\|p_h^{i,n}-p^{i,n}\|_{0}^2\Big)
\|\nabla\upsilon^{i,n}\|_{0}\notag\\
&\leq  C (\tau+ h^k+\sum\limits^{2}_{i=1}\|p_h^{i,n}-p^{i,n}\|_{0}+ h^{-1}\sum\limits^{2}_{i=1}\|p_h^{i,n}-p^{i,n}\|_{0}^2)\|\nabla\upsilon^{i,n}\|_{0}.
\end{align}
Note that
\begin{align*}\frac 1 {2\tau} (\|\upsilon^{i,n}\|_0^2-\|\upsilon^{i,n-1}\|_0^2)&=\frac 1 {\tau} \|\upsilon^{i,n}\|_0^2
-\frac 1 {2\tau}(\|\upsilon^{i,n-1}\|_0^2+\|\upsilon^{i,n}\|_0^2)\notag\\
&\leq \frac 1 {\tau} m_h(v^{i,n},v^{i,n})-\frac 1 {\tau}\|\upsilon^{i,n}\|_0\|\upsilon^{i,n-1}\|_0\notag\\
&\leq \frac 1 {\tau} m_h(v^{i,n},v^{i,n})-\frac 1 {\tau} m_h(v^{i,n-1},v^{i,n}).
\end{align*}
Combining the above with (\ref{mh}), we get
\begin{align}\label{5.7-1}
\frac 1 {2\tau} (\|\upsilon^{i,n}\|_0^2&-\|\upsilon^{i,n-1}\|_0^2)+\|\nabla\upsilon^{i,n}\|_0^2\leq
m_h\Big(\frac{\upsilon^{i,n}-\upsilon^{i,n-1}}{\tau},\upsilon^{i,n}\Big)+a_h(\upsilon^{i,n},\upsilon^{i,n})\notag\\
&\leq C (\tau+ h^k+\sum\limits^{2}_{i=1}\|p_h^{i,n}-p^{i,n}\|_{0}+ h^{-1}\sum\limits^{2}_{i=1}\|p_h^{i,n}-p^{i,n}\|_{0}^2)\|\nabla\upsilon^{i,n}\|_{0}.
\end{align}
Next, we shall use the mathematical induction to show
\begin{align}\label{l2up} \|\upsilon^{i,n}\|_{0}\leq C(\tau+h^k),~~n=0,1,\cdot,N.
\end{align}
First, there holds
\begin{align*} \|\upsilon^{i,0}\|_{0}&=\|p_h^{i,0}-\Pi_k^0p^{i,0}\|_0\leq \|p_h^{i,0}-p^{i,0}\|_0+\|p^{i,0}-\Pi_k^0p^{i,0}\|_0\notag\\
&\leq C(\tau+h^k).
\end{align*}
Assume
\begin{align}\label{assu2}\|\upsilon^{i,n}\|_{0}\leq C(\tau+h^k),~~\forall n=1,2,\cdots,J,~~1\leq J\leq N-1.
\end{align}
From (\ref{5.7-1}), it follows that
\begin{align}\label{5eqn2}
\|\upsilon^{i,n}\|_0^2&-\|\upsilon^{i,n-1}\|_0^2+2\tau\|\nabla\upsilon^{i,n}\|_0^2\notag\\
&\leq C (\tau^2+ \tau h^k+\tau\sum\limits^{2}_{i=1}\|p_h^{i,n}-p^{i,n}\|_{0}+ \tau h^{-1}\sum\limits^{2}_{i=1}\|p_h^{i,n}-p^{i,n}\|_{0}^2)\|\nabla\upsilon^{i,n}\|_{0}\notag\\
&\leq C\Big(\tau^4+\tau^2h^{2k}+\tau^2(\sum\limits^{2}_{i=1}\|p_h^{i,n}-p^{i,n}\|_{0})^2
+\tau^2h^{-2}(\sum\limits^{2}_{i=1}\|p_h^{i,n}-p^{i,n}\|_{0}^2)^2\Big)+
\epsilon\|\nabla\upsilon^{i,n}\|_{0}^2.
\end{align}
Note that
\begin{align*}\tau^2 h^{-2}(\sum\limits^{2}_{i=1}\|p_h^{i,n}-p^{i,n}\|_{0}^2)^2&\leq \tau^2 h^{-2}(\sum\limits^{2}_{i=1}(\|\upsilon^{i,n}\|_{0}^2+\|\rho^{i,n}\|_{0}^2)^2\notag\\
&\leq \tau^2 h^{-2}(\sum\limits^{2}_{i=1}(\|\upsilon^{i,n}\|_{0}^2)^2+\tau^2 h^{4k-2}\notag\\
&\leq \tau^2 h^{-2}(\tau+h^k)^2\sum\limits^{2}_{i=1}\|\upsilon^{i,n}\|_{0}^2+\tau^2 h^{4k-2}\notag~~(by~~(\ref{assu2}))\\
&\leq C\tau \sum\limits^{2}_{i=1}\|\upsilon^{i,n}\|_{0}^2+\tau^2 h^{4k-2},
\end{align*}
where we have used $\tau\leq Ch$. Inserting the above into (\ref{5eqn2}), it yields
\begin{align*}
\|\upsilon^{i,n}\|_0^2&-\|\upsilon^{i,n-1}\|_0^2
\leq C(\tau^4+\tau^2h^{2k}+\tau \sum\limits^{2}_{i=1}\|\upsilon^{i,n}\|_{0}^2),~~\forall n=1,2,\cdots,J,~~1\leq J\leq N-1.
\end{align*}
Then, we have
\begin{align*}\sum\limits^{2}_{i=1}(\|\upsilon^{i,n+1}\|_0^2&-\|\upsilon^{i,n}\|_0^2)\leq C(\tau^4+\tau^2h^{2k}+\tau \sum\limits^{2}_{i=1}\|\upsilon^{i,n+1}\|_{0}^2),~~n=0,1,\cdots,J,~~1\leq J\leq N-2.
\end{align*}
Summing up for the index $n=0,1,\cdots, J$, $1\leq J\leq N-1$, then we get
\begin{align*}\sum\limits^{2}_{i=1}(\|\upsilon^{i,J+1}\|_0^2-\|\upsilon^{i,0}\|_0^2)&\leq C\sum\limits^{J}_{n=0} (
\tau^4+\tau^2h^{2k}+\tau \sum\limits^{2}_{i=1}\|\upsilon^{i,n+1}\|_{0}^2).
\end{align*}
Hence
\begin{align*}\sum\limits^{2}_{i=1}\|\upsilon^{i,J+1}\|_0^2&\leq \sum\limits^{2}_{i=1}\|\upsilon^{i,0}\|_0^2+C\sum\limits^{J}_{n=0} (
\tau^4+\tau^2h^{2k}+\tau \sum\limits^{2}_{i=1}\|\upsilon^{i,n+1}\|_{0}^2)\notag\\
&\leq C (\tau^2+h^{2k}+\sum\limits^{J}_{n=0}(\tau \sum\limits^{2}_{i=1}\|\upsilon^{i,n+1}\|_{0}^2))
\end{align*}
By the discrete Gronwall's inequality, we get
\begin{align}\label{5eqn1}\sum\limits^{2}_{i=1}\|\upsilon^{i,J+1}\|_0^2\leq C(\tau^2+h^{2k}).
\end{align}
That is
\begin{align*}\|\upsilon^{i,J+1}\|_0\leq c(\tau+h^k),~~1\leq J\leq N-1.
\end{align*}
Hence (\ref{l2up}) holds for $n=0,1,\cdots,N$. We can easily obtain
\begin{align}\label{pih1}
\sum\limits^{2}_{i=1}\|p_h^{i,n}-p^{i,n}\|_0&\leq \sum_{i=1}^2( \|\upsilon^{i,n}\|_0+\|\rho^{i,n}\|_0)\notag\\
&\leq C(\tau+h^k).
\end{align}
Using the similar analysis for Theorem \ref{theorem3.3} and from (\ref{pih1}), we have
\begin{align}\|\phi_h^n-\phi^n\|_0\leq C(\tau+h^k).
\end{align}
This completes the proof.$\hfill\Box$
\vspace{3mm}

\setcounter{lemma}{0}
\setcounter{theorem}{0}
\setcounter{corollary}{0}
\setcounter{equation}{0}

\section{Numerical Results}
In this section, we report a numerical example to test the practical performance of the virtual element method for solving (\ref{pnp equation}). The implement of the numerical experiment is based on \cite{beirao2014hitchhiker}. All the computations are carried out in Fortran 90 on the computer with CPU-2.90GHz (Intel(R) Core (TM) i5-10400F), RAM-16GB.\\
\vspace{3mm}

\noindent\textbf{Example 5.1} Consider problem (\ref{pnp equation})-(\ref{boundary}) and the right-hand side functions are determined from the exact solution (cf. \cite{sun2016error})
\begin{equation}\label{exact s}
\left\{\begin{aligned}
&\phi(t,x,y)=(1-e^{-t}) \sin (\pi x) \sin (\pi y), \\
&p^{1}(t,x,y)=\sin {t} \sin (2 \pi x) \sin (2 \pi y), \\
&p^{2}(t,x,y)= \sin (2 t) \sin (3 \pi x) \sin (3 \pi y).
\end{aligned}\right.
\end{equation}
In the test we consider the time interval [0, 1] and the computational domain $\Omega=[0, 1] \times [0, 1]$ with $T=1$ and time step $\tau=h^2$. The rectangular domain is discretized with several different types of polygonal meshes, viz., triangle, square, non-convex polygons, mixed-polygon, random Voronoi and smooth Voronoi meshes, which are shown in (a)-(f) of Figures 1, respectively.

\begin{figure}[H]
    \begin{minipage}[t]{0.5\linewidth}
    \centering
    \includegraphics[height=5.5cm,width=5.5cm]{triangles.pdf}
    \caption*{(a) Triangular mesh}
    \label{fig:side:a}
    \end{minipage}
    \begin{minipage}[t]{0.5\linewidth}
    \centering
    \includegraphics[height=5.5cm,width=5.5cm]{square.pdf}
    \caption*{(b) Square mesh}
    \label{fig:side:b}
    \end{minipage}
\end{figure}

\begin{figure}[H]
 \begin{minipage}[t]{0.5\linewidth}
	\centering
	\includegraphics[height=5.5cm,width=5.5cm]{non-convex.pdf}
	\caption*{(c) Non-convex mesh}
	\label{fig:side:c}
\end{minipage}
\begin{minipage}[t]{0.5\linewidth}
	\centering
	\includegraphics[height=5.5cm,width=5.5cm]{mpolygon.pdf}
	\caption*{(d) Mixed-polygon mesh.}
	\label{fig:side:d}
\end{minipage}
\end{figure}

\begin{figure}[H]
    \begin{minipage}[t]{0.5\linewidth}
    \centering
\label{Fig.sub.1}
\includegraphics[height=5.5cm,width=5.5cm]{Voronoi.pdf}
     \caption*{(e) Random Voronoi mesh}
    \label{fig:side:e}
     \end{minipage}
    \begin{minipage}[t]{0.5\linewidth}
    \centering
\label{Fig.sub.2}
\includegraphics[height=5.5cm,width=5.5cm]{sVoronoi.pdf}
 \caption*{(f) Smooth Voronoi mesh}
    \label{fig:side:f}
    \end{minipage}
\caption{six types of polygonal meshes}
\end{figure}

\begin{figure}[H]
	\centering
	{
		\begin{minipage}{8cm}
			\centering
			\includegraphics[scale=0.5]{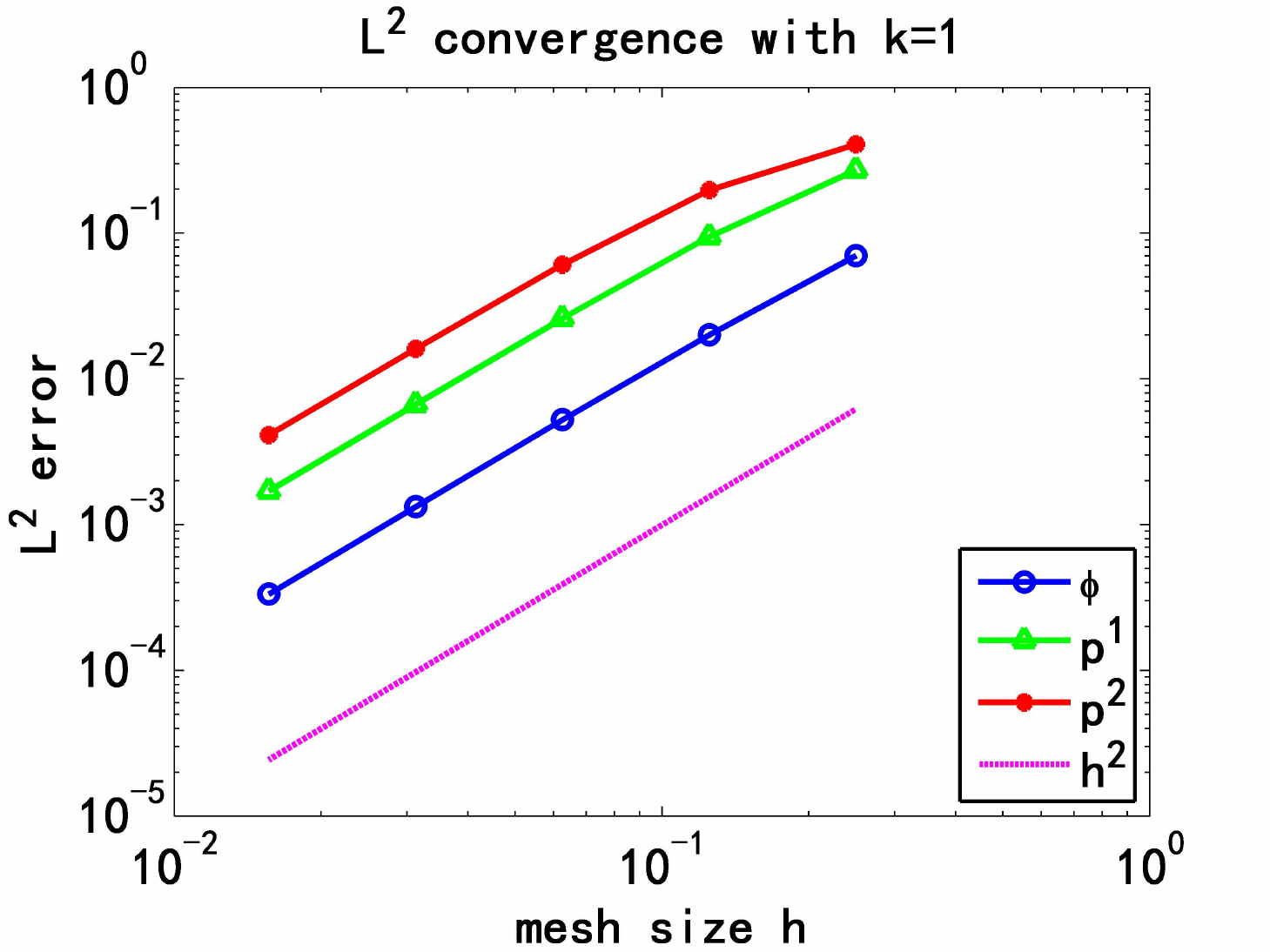}
		\end{minipage}
	}
	{
		\begin{minipage}{7cm}
			\centering
			\includegraphics[scale=0.5]{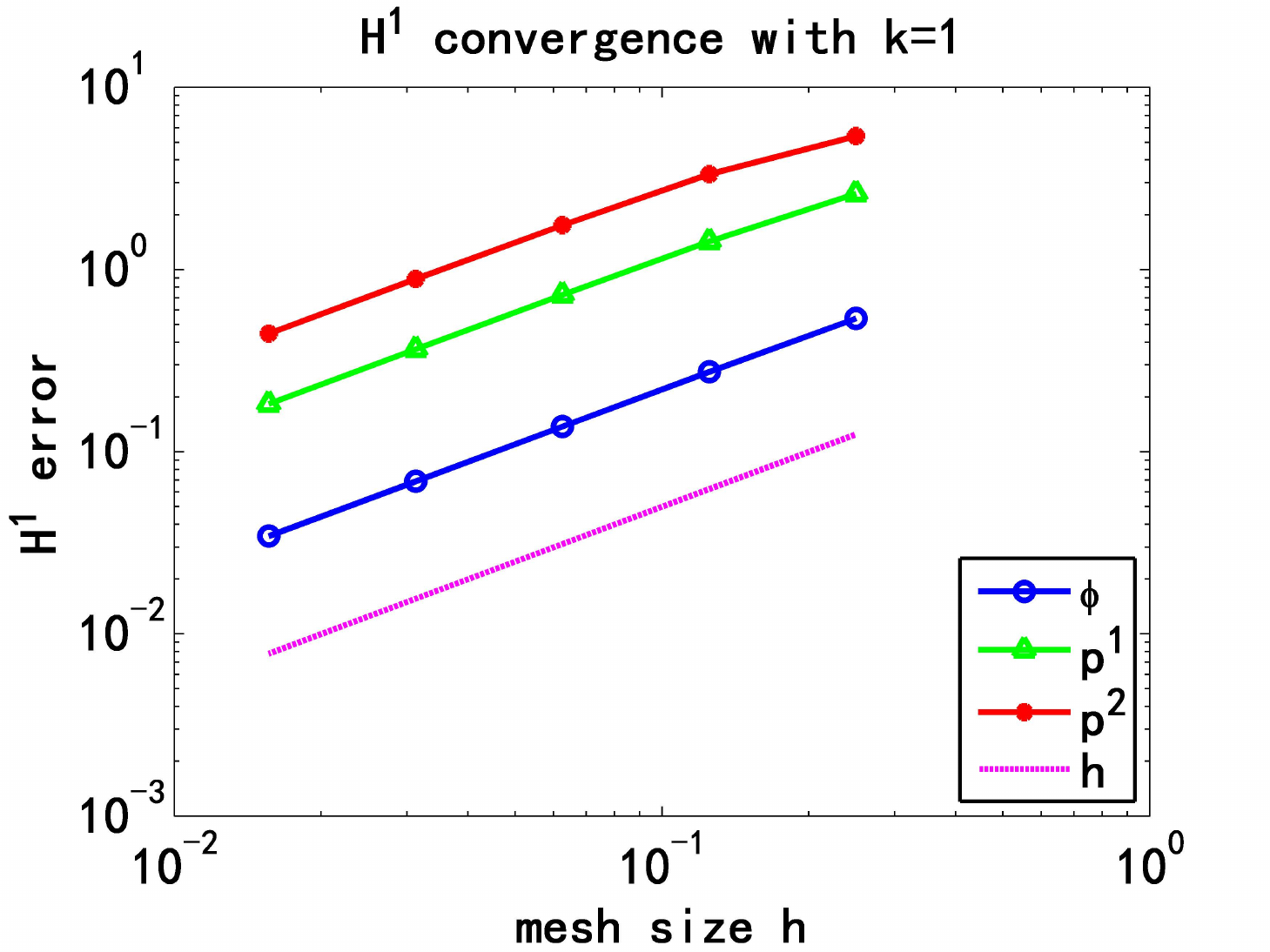}
		\end{minipage}
	}
	
	\caption{h-convergence on triangular mesh with $t=1.0$.}
	\label{tri}
\end{figure}

\begin{figure}[H]
	\centering
	{
		\begin{minipage}{8cm}
			\centering
			\includegraphics[scale=0.5]{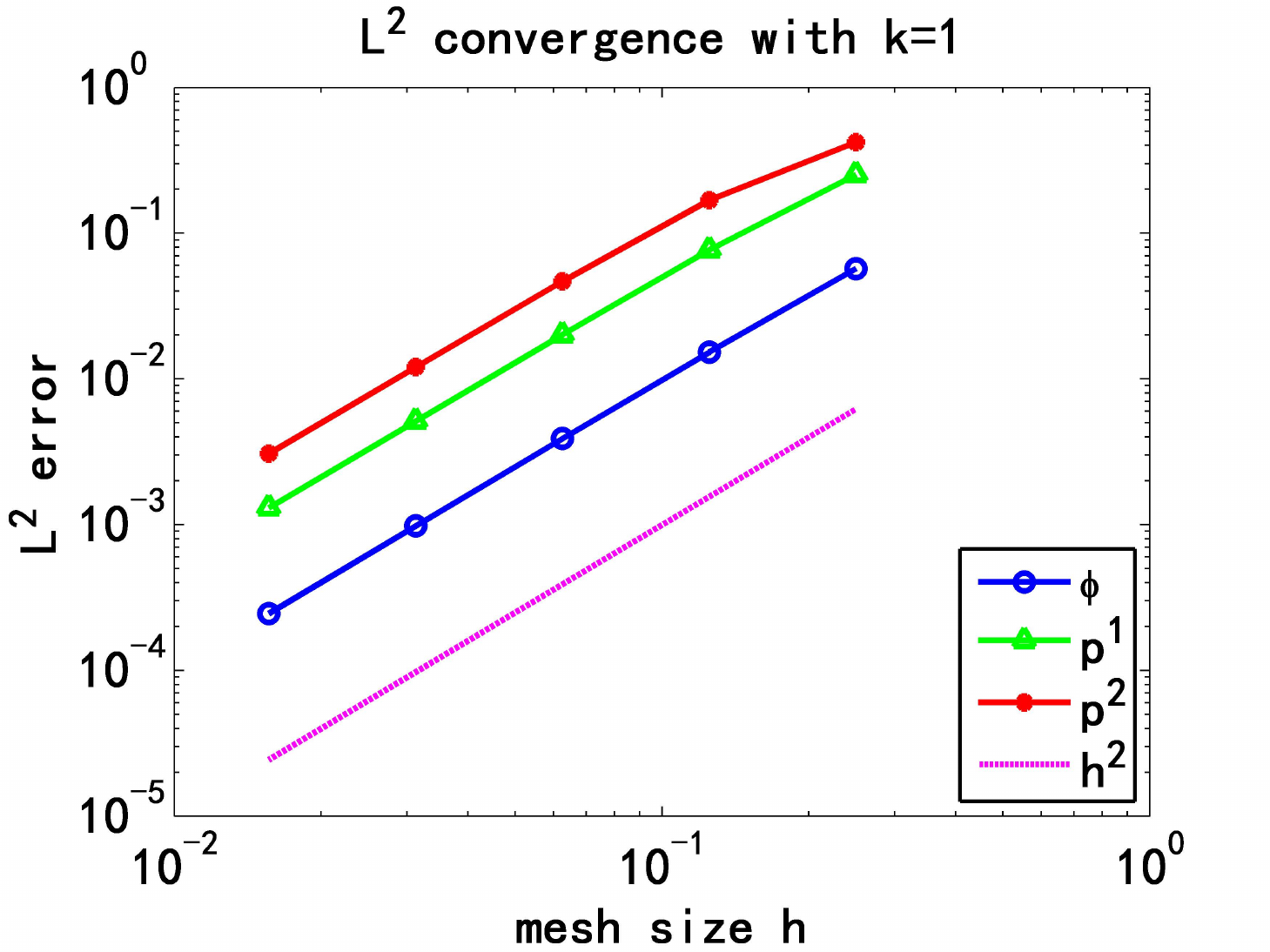}
		\end{minipage}
	}
	{
		\begin{minipage}{7cm}
			\centering
			\includegraphics[scale=0.5]{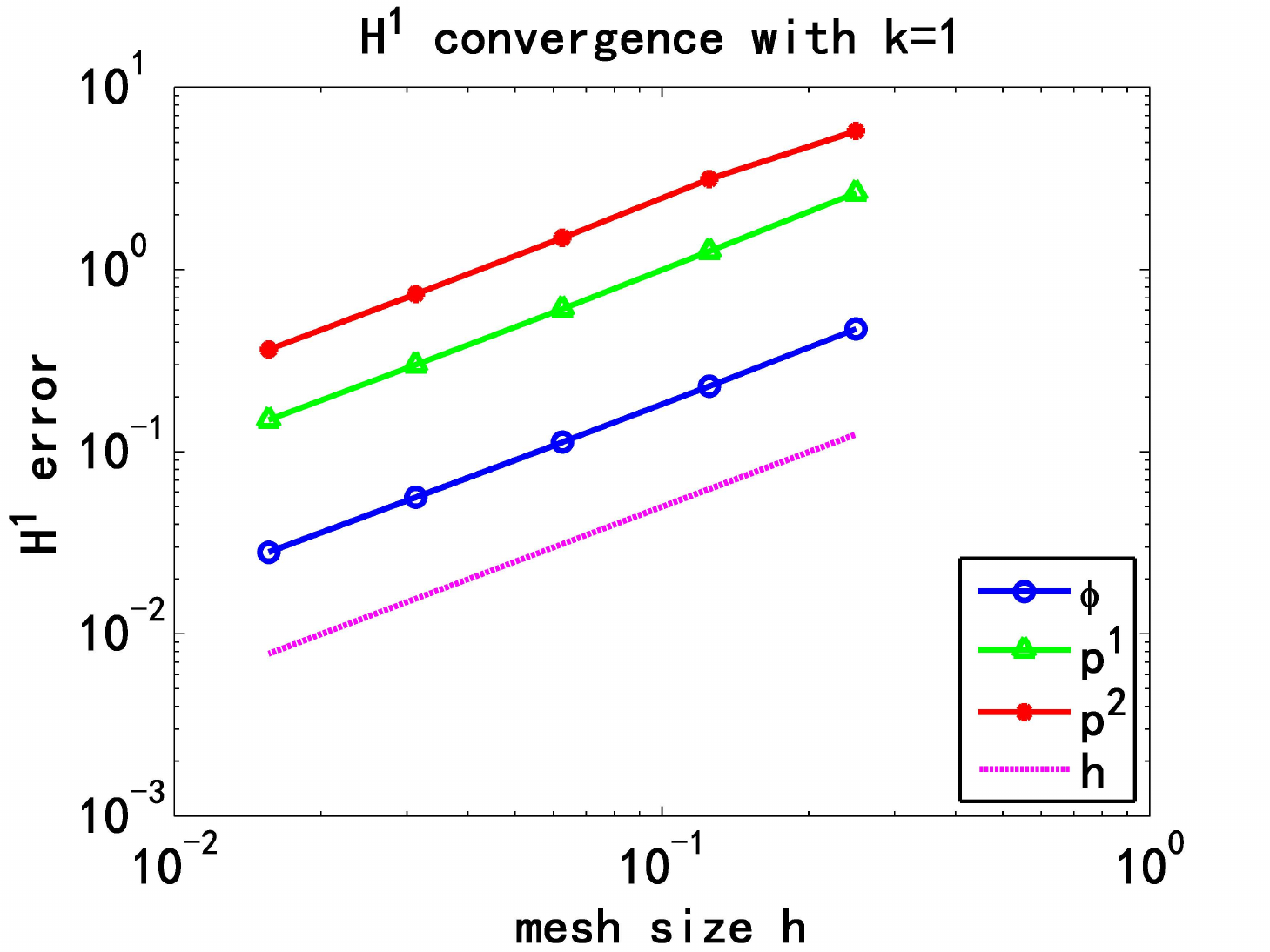}
		\end{minipage}
	}
	
	\caption{h-convergence on square mesh with $t=1.0$.}
	\label{squa}
\end{figure}

\begin{figure}[H]
	\centering
	{
		\begin{minipage}{8cm}
			\centering
			\includegraphics[scale=0.48]{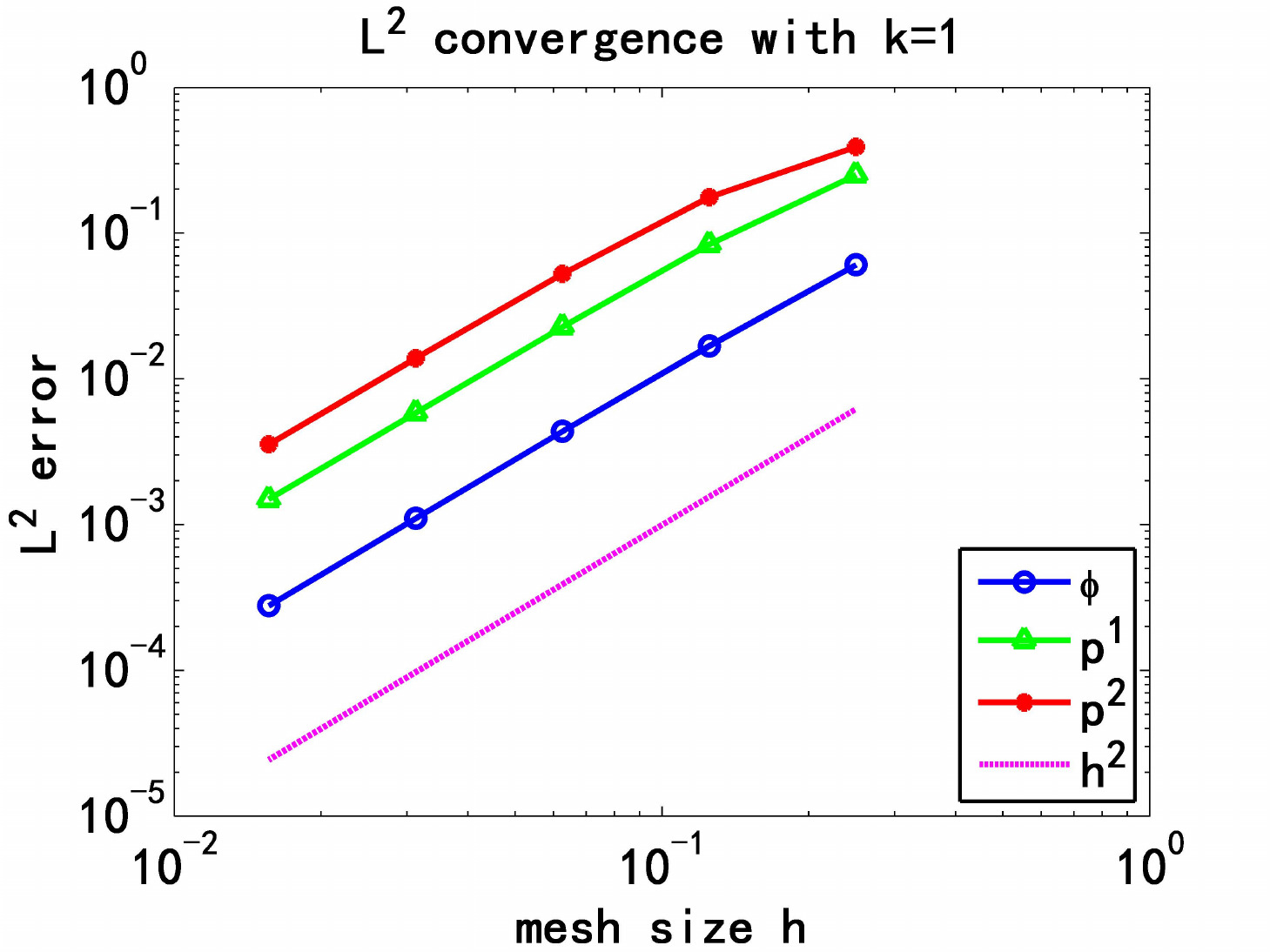}
		\end{minipage}
	}
	{
		\begin{minipage}{7cm}
			\centering
			\includegraphics[scale=0.48]{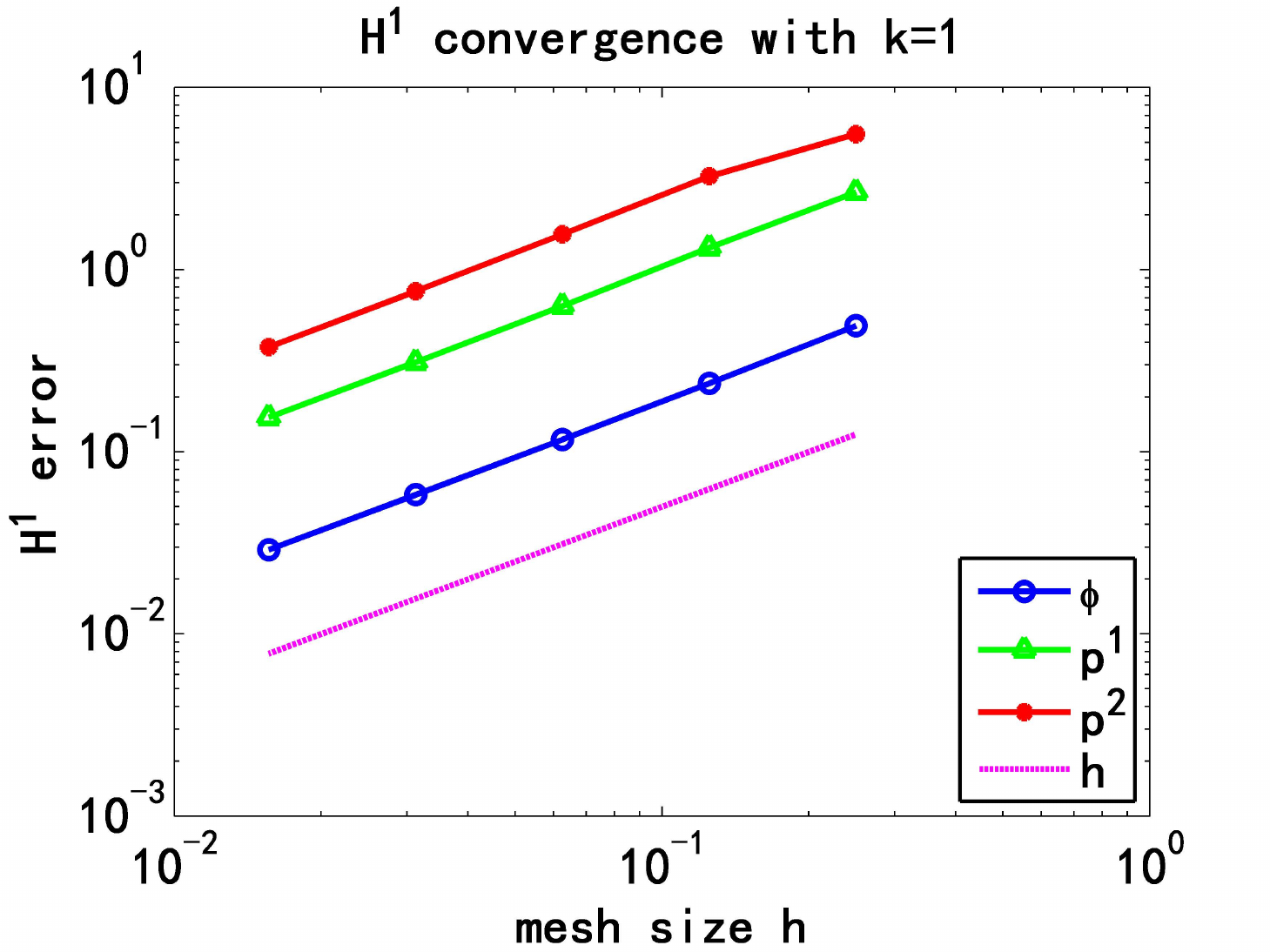}
		\end{minipage}
	}
	
	\caption{h-convergence on non-convex mesh with $t=1.0$. }
	\label{non}
\end{figure}

\begin{figure}[H]
	\centering
	{
		\begin{minipage}{8cm}
			\centering
			\includegraphics[scale=0.48]{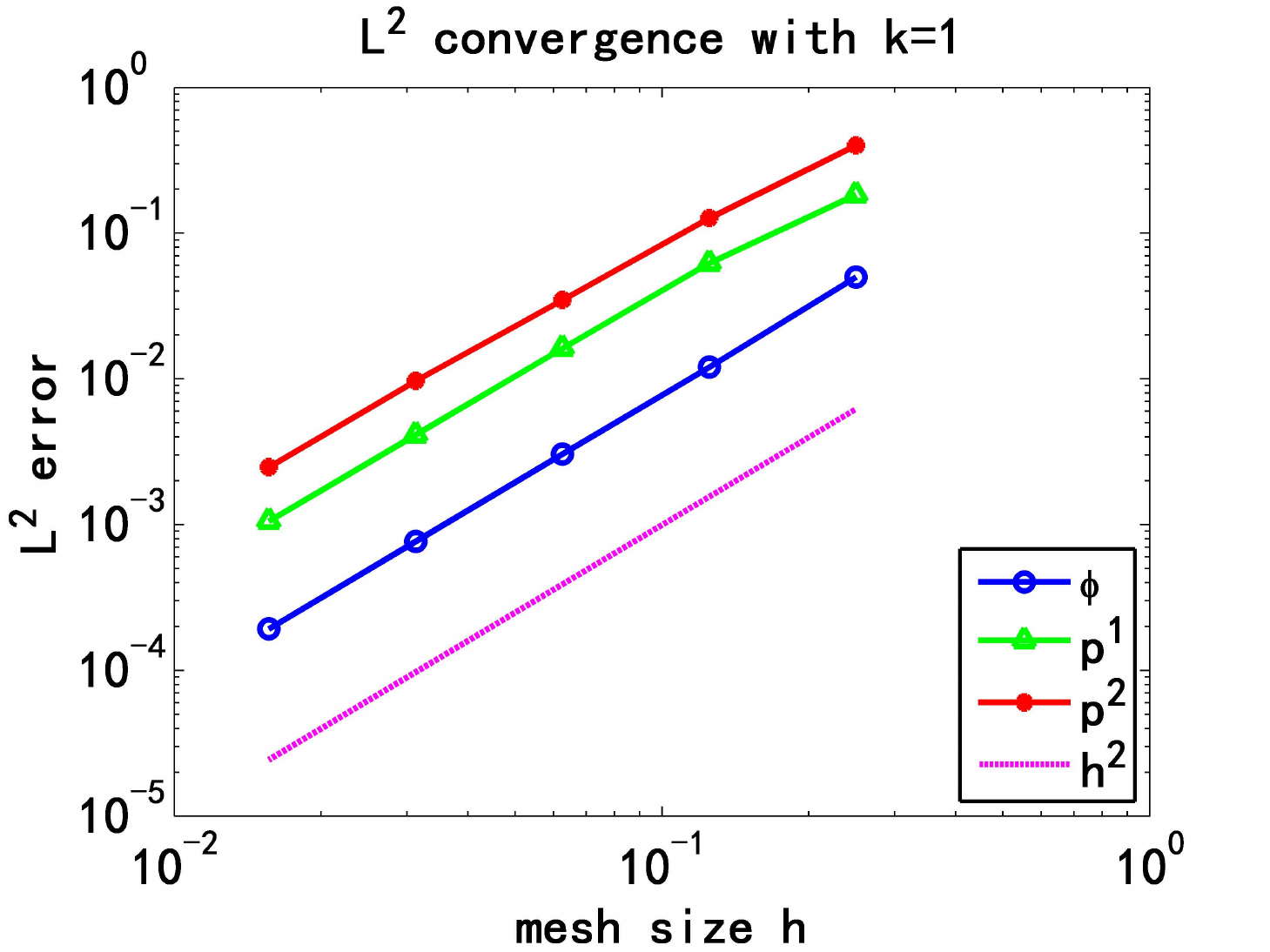}
		\end{minipage}
	}
	{
		\begin{minipage}{7cm}
			\centering
			\includegraphics[scale=0.48]{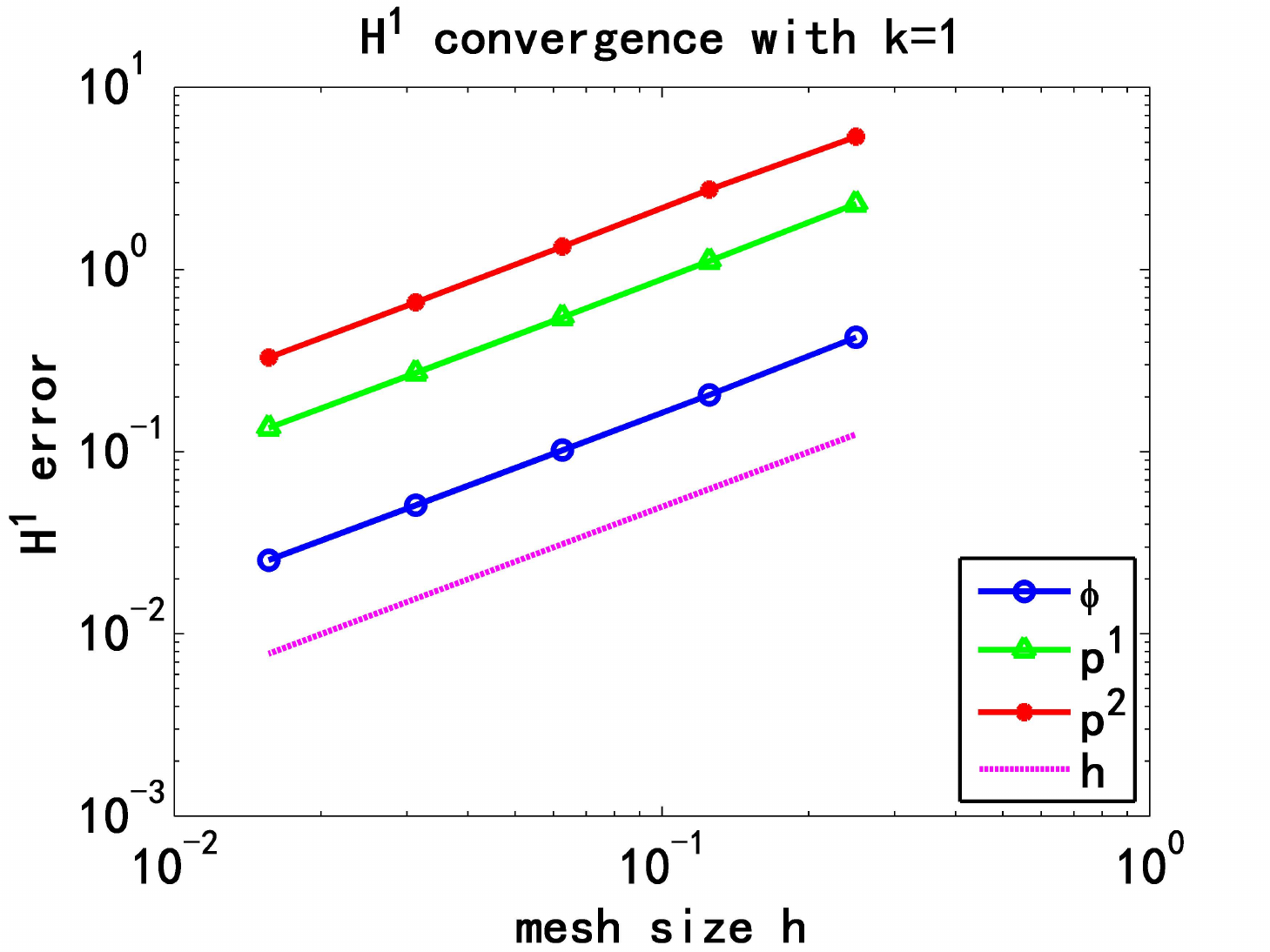}
		\end{minipage}
	}
	\caption{h-convergence on mixed-polygon mesh with $t=1.0$.}
	\label{mix}
\end{figure}

\begin{figure}[H]
	\centering
	{
		\begin{minipage}{8cm}
			\centering
			\includegraphics[scale=0.5]{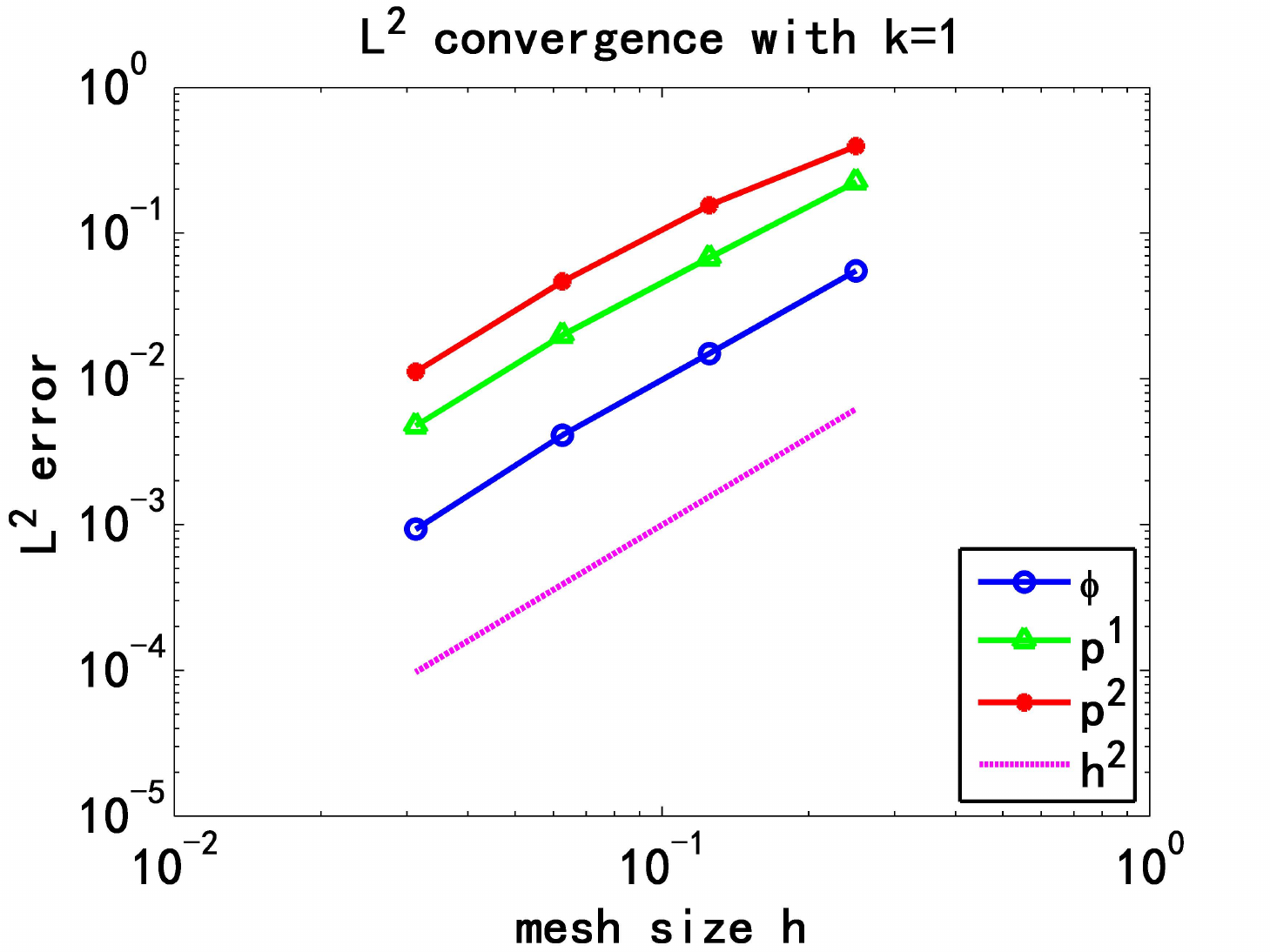}
		\end{minipage}
	}
	{
		\begin{minipage}{7cm}
			\centering
			\includegraphics[scale=0.5]{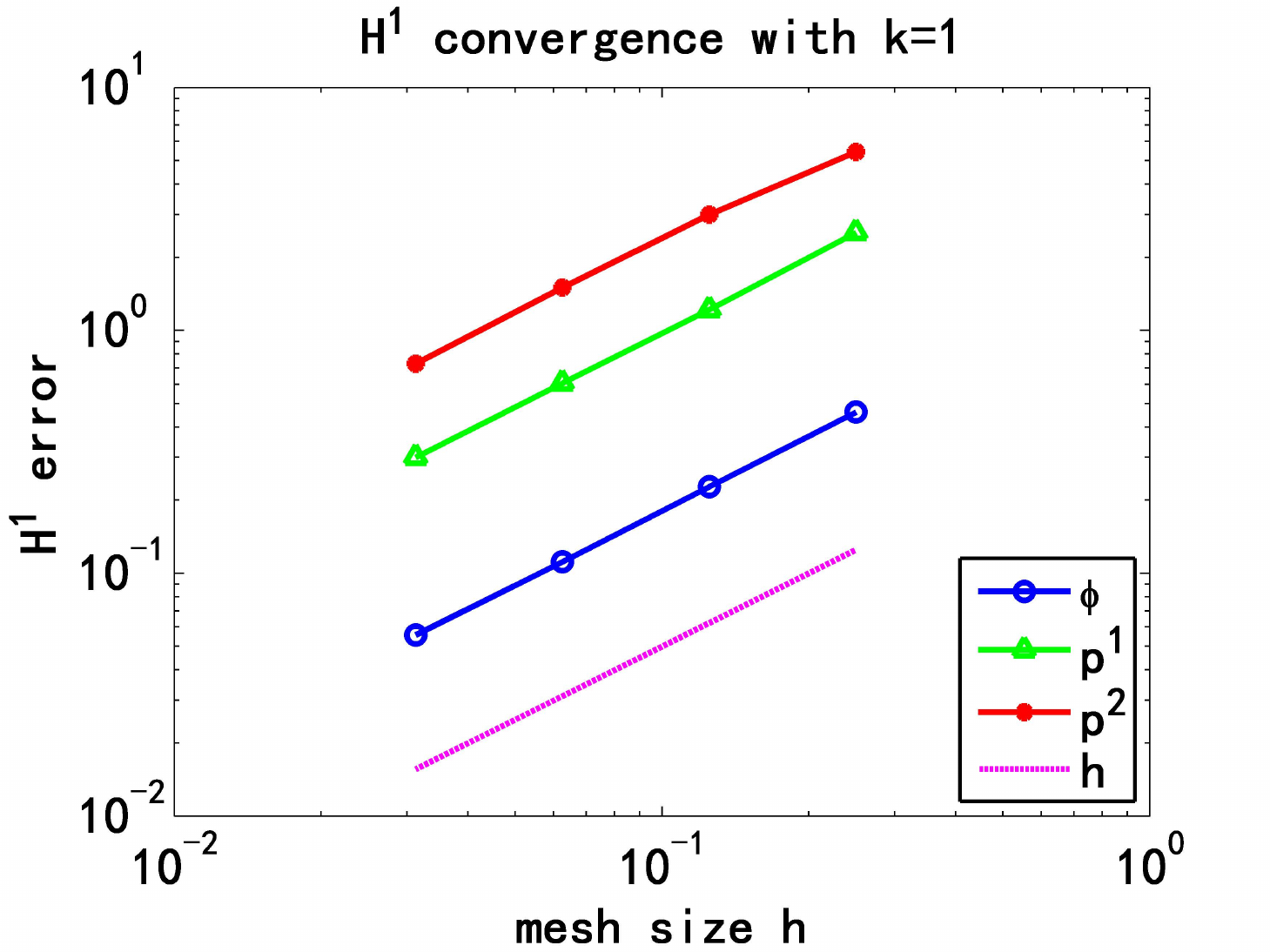}
		\end{minipage}
	}
	\caption{h-convergence on random Voronoi mesh with $t=1.0$.}
	\label{vor}
\end{figure}

\begin{figure}[H]
	\centering
	{
		\begin{minipage}{8cm}
			\centering
			\includegraphics[scale=0.5]{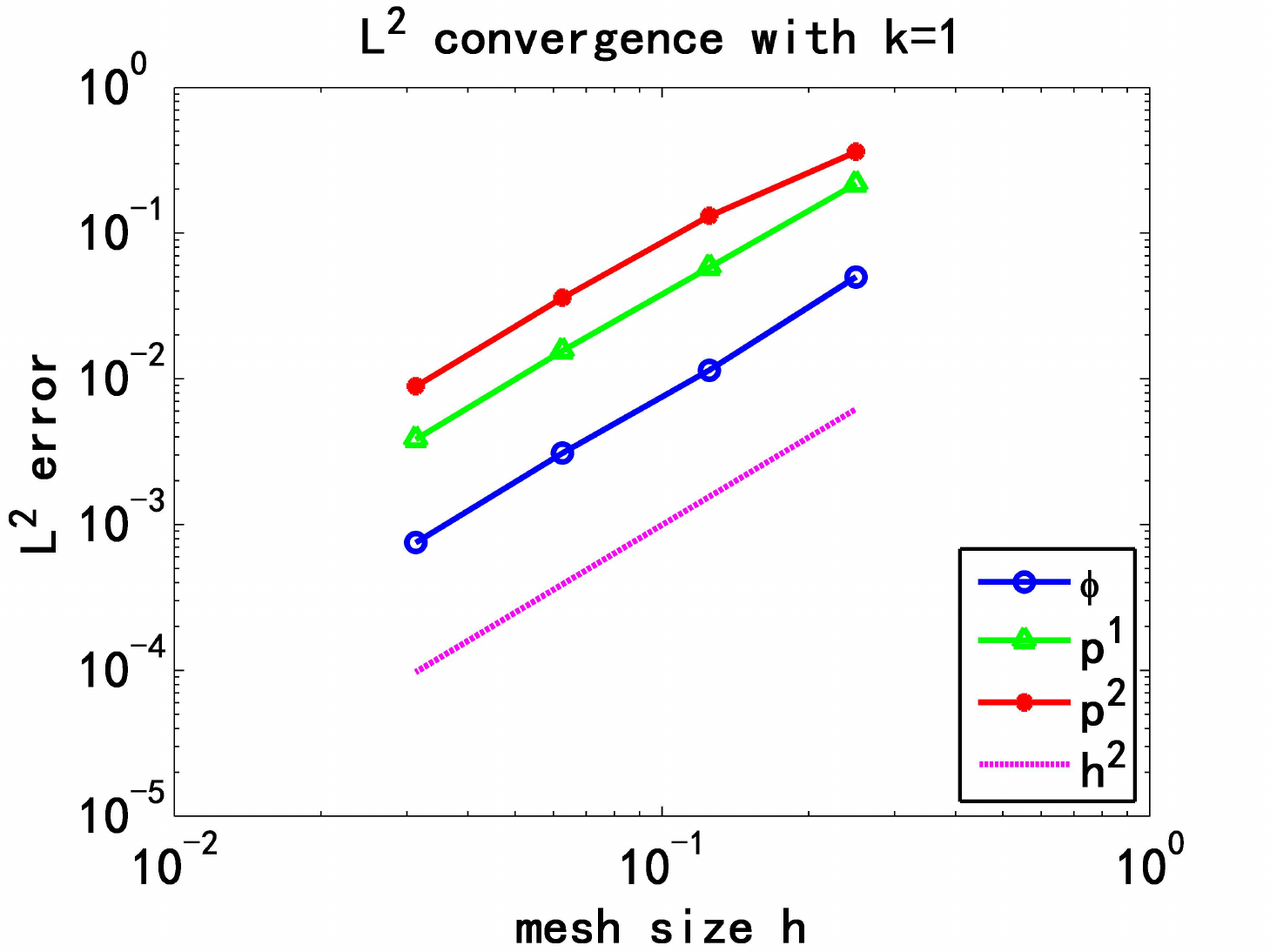}
		\end{minipage}
	}
	{
		\begin{minipage}{7cm}
			\centering
			\includegraphics[scale=0.5]{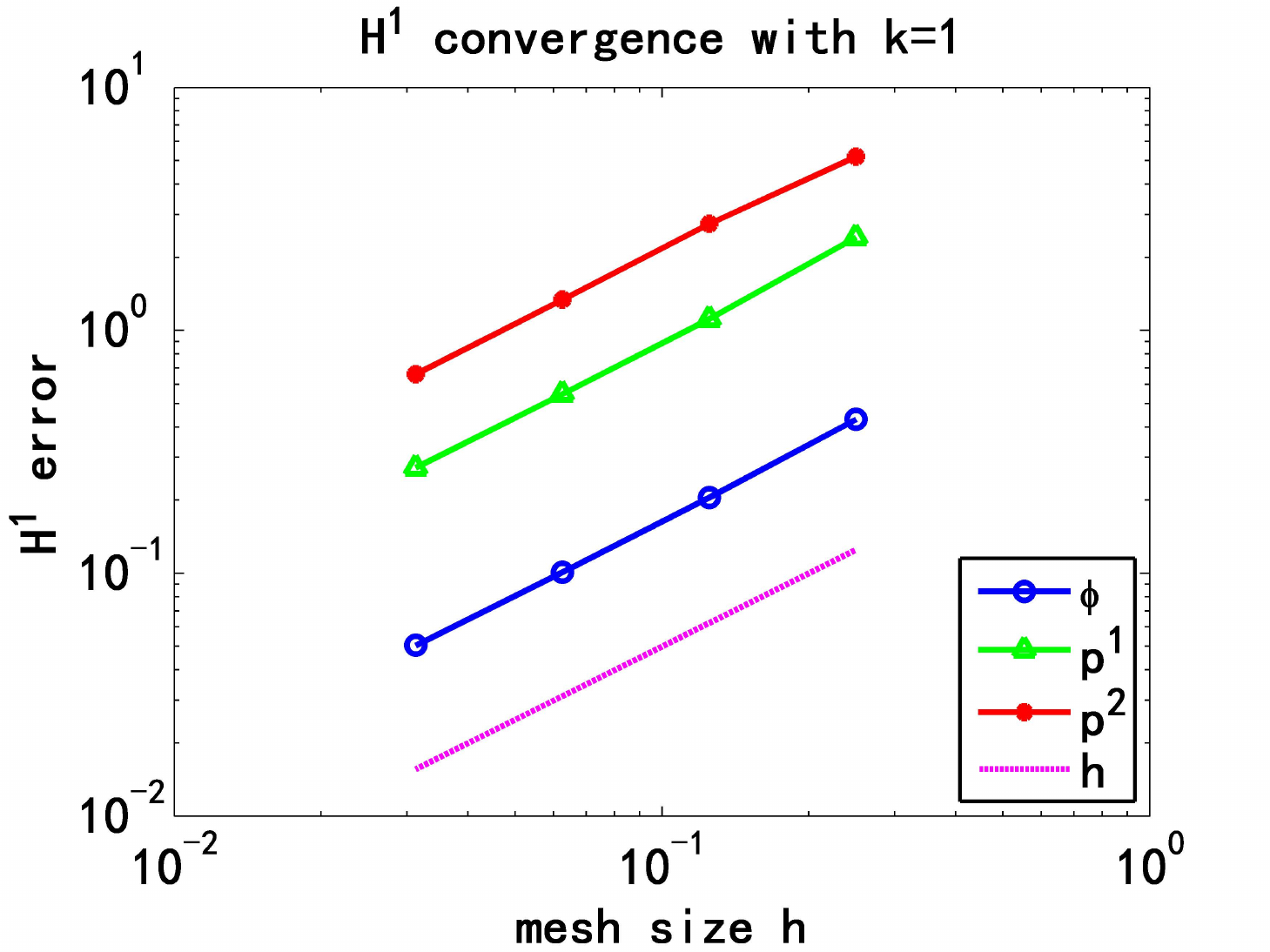}
		\end{minipage}
	}
	\caption{h-convergence on smooth Voronoi mesh with $t=1.0$.}
	\label{sm}
\end{figure}

 As the VEM solution can not be known explicitly inside the elements, the convergence of VEM is evaluated  through the relative $L^2$ norm and $H^{1}$ semi-norm  using the projection operator $\Pi_k^\nabla$ onto $\mathbb{P}_k$, that is
 \beas e_{L^2}:=\sqrt{\sum_{E\in \mathcal{T}^h} \|u-\Pi_k^{\triangledown}u_h\|_{0,E}^2}, ~~e_{H^1}:=\sqrt{\sum_{E\in \mathcal{T}^h} \|\nabla (u-\Pi_k^{\triangledown}u_h)\|_{0,E}^2},
 \eeas
\noindent where $u$ is corresponding the exact solution $\phi,~p^{1}~\text{or}~p^{2}$ in (\ref{exact s}), $u_h$ represents the VEM solution with the order $k=1$,  $\Pi_k^\nabla$ is defined by (\ref{3.7}). To display the convergence results in the figures, we shall use the mesh-size parameter $h$ which is measured in following ratios (cf. \cite{da2017high}).
\begin{equation}
h=\left(\frac{|\Omega|}{N_{E}}\right)^{1 / 2},
\end{equation}
where $N_E$ is the number of polygons. The convergences of the errors in $L^2$ and $H^1$ norms at $t=1.0$ are displayed in Figures \ref{tri}-\ref{sm}, which infers that the convergence orders approximate the  second order and first order in $L^2$ norm and $H^1$ norm, respectively. Figures \ref{tri}-\ref{sm} show that the numerical result matches with the theoretical  result established in Section 4.

\vspace{3mm}
\section{Conclusion}
In this paper, we study VEMs to approximate the solutions of the time-dependent PNP equations on general polygonal meshes. 
We derive the a priori error estimates for semidiscrete and fully discrete schemes, respectively. The numerical errors have been conducted to show that the convergence orders agree with the theoretical results well. This method is expected to be extended to more complex PNP models, for example, PNP equations for semiconductor devices and three-dimensional ion channel.\\

\noindent
{\bf Acknowledgement }The authors would like to thank Yang Liu and Zhiquan Wu for their valuable discussions on numerical experiments. Y. Yang was supported by the Guangxi Natural Science Foundation (2020GXNSFAA159098) and China NSF (NSFC 11561016). S. Shu was supported by the China NSF (NSFC 11971414). Y. Liu was supported by GUET Excellent Graduate Thesis Program (No. 2020YJSPYA02)

\bibliography{artpnp}
\bibliographystyle{abbrv}
\end{document}